\documentclass[12pt,twoside]{article}
\usepackage{a4wide}
\usepackage{amssymb}
\usepackage{amsthm}
\usepackage{amsmath}
\usepackage{amscd}
\usepackage[mathscr]{eucal}
\usepackage[all]{xy}
\usepackage{verbatim}
\numberwithin{equation}{section}

\theoremstyle{plain}
\newtheorem{theorem}{Theorem}[section]
\newtheorem{corollary}[theorem]{Corollary}
\newtheorem{lemma}[theorem]{Lemma}
\newtheorem{proposition}[theorem]{Proposition}

\theoremstyle{definition}
\newtheorem{definition}[theorem]{Definition}
\newtheorem{remark}[theorem]{Remark}
\newtheorem{example}[theorem]{Example}
\newtheorem*{notation}{Notation}

\theoremstyle{remark}

\newcommand{\R}{\mathbb{R}}
\newcommand{\A}{\mathbb{A}}
\newcommand{\Q}{\mathbb{Q}}
\newcommand{\Z}{\mathbb{Z}}

\newcommand{\C}{\mathbb{C}}
\newcommand{\h}{\mathbb{H}}
\renewcommand{\H}{\mathbb{H}}
\renewcommand{\P}{\mathbb{P}}

\newcommand{\zxz}[4]{\begin{pmatrix} #1 & #2 \\ #3 & #4 \end{pmatrix}}
\newcommand{\abcd}{\zxz{a}{b}{c}{d}}

\newcommand{\kzxz}[4]{\left(\begin{smallmatrix} #1 & #2 \\ #3 & #4\end{smallmatrix}\right) }
\newcommand{\kabcd}{\kzxz{a}{b}{c}{d}}

\newcommand{\calF}{\mathcal{F}}

\newcommand{\calH}{\mathcal{H}}

\newcommand{\calK}{\mathcal{K}}
\newcommand{\calL}{\mathcal{L}}
\newcommand{\calM}{\mathcal{M}}
\newcommand{\calN}{\mathcal{N}}
\newcommand{\calO}{\mathcal{O}}

\newcommand{\calS}{\mathcal{S}}
\newcommand{\calT}{\mathcal{T}}

\newcommand{\calX}{\mathcal{X}}

\newcommand{\fraka}{\mathfrak a}
\newcommand{\frakb}{\mathfrak b}
\newcommand{\frakc}{\mathfrak c}
\newcommand{\frakd}{\mathfrak d}

\newcommand{\frakl}{\mathfrak l}

\newcommand{\IC}{{\mathbb C}}

\newcommand{\IH}{{\mathbb H}}

\newcommand{\IP}{{\mathbb P}}
\newcommand{\IQ}{{\mathbb Q}}
\newcommand{\IR}{{\mathbb R}}

\newcommand{\IZ}{{\mathbb Z}}
\newcommand{\eps}{\varepsilon}
\newcommand{\bs}{\backslash}
\newcommand{\norm}{\operatorname{N}}
\newcommand{\vol}{\operatorname{vol}}
\newcommand{\tr}{\operatorname{tr}}

\newcommand{\Ad}{\operatorname{Ad}}

\newcommand{\Cl}{\operatorname{Cl}}
\newcommand{\Sl}{\operatorname{SL}}
\newcommand{\PSL}{\operatorname{PSL}}

\newcommand{\Spin}{\operatorname{Spin}}
\newcommand{\GSpin}{\operatorname{GSpin}}
\newcommand{\CG}{\operatorname{CG}} 

\newcommand{\Orth}{\operatorname{O}}
\newcommand{\Uni}{\operatorname{U}}

\newcommand{\Aut}{\operatorname{Aut}}
\newcommand{\Mat}{\operatorname{M}}

\newcommand{\End}{\operatorname{End}}

\newcommand{\Pet}{\text{\rm Pet}}
\newcommand{\Gr}{\operatorname{Gr}}

\newcommand{\SO}{\operatorname{SO}}
\newcommand{\Gal}{\operatorname{Gal}}

\newcommand{\CL}{\operatorname{Cl}}

\newcommand{\supp}{\operatorname{supp}}

\newcommand{\Pic}{\operatorname{Pic}}
\newcommand{\Div}{\operatorname{Div}}

\newcommand{\DN}{\operatorname{DN}}

\newcommand{\pr}{\operatorname{pr}}

\newcommand{\ord}{\operatorname{ord}}
\newcommand{\id}{\operatorname{id}}

\newcommand{\dv}{\operatorname{div}}

\newcommand{\CM}{\mathcal{CM}}
\newcommand{\ch}{\operatorname{CH}}
\newcommand{\rat}{\operatorname{Rat}}

\newcommand{\cc}{\operatorname{c}}

\newcommand{\darstell}[2]{#1 \mathcal{O}_F + #2 \mathcal{O}_F}
\long\def\symbolfootnote[#1]#2{\begingroup%
\def\thefootnote{\fnsymbol{footnote}}\footnote[#1]{#2}\endgroup}
\title{Hilbert modular forms and their applications}

\author{Jan Hendrik Bruinier}
\date{\today}

\begin{document}

\maketitle

\symbolfootnote[0]{ Mathematisches Institut, Universit\"at zu
K\"oln, Weyertal 86--90, D-50931 K\"oln, Germany }

\tableofcontents

\section*{Introduction}

The present notes contain the material of the lectures given by the author
at the summer school on ``Modular Forms and their Applications'' at
the Sophus Lie Conference Center in the summer of 2004.

We give an introduction to the theory of Hilbert modular forms and
some geometric and arithmetic applications. We tried to keep the
informal style of the lectures. In particular, we often do not work in
greatest possible generality, but rather consider a reasonable special
case, in which the main ideas of the theory become
clear. 

For a more comprehensive account to Hilbert modular varieties, we refer
to the books by Freitag \cite{Fr}, Garrett \cite{Ga}, van der Geer \cite{Ge}, and Goren \cite{Go}. 
We hope that the present text will be a useful
addition to these references.  

Hilbert modular surfaces can also be realized as modular varieties
corresponding to the orthogonal group of a rational quadratic space of
type $(2,2)$. This viewpoint leads to several interesting features of
these surfaces.  For instance, they come with a natural family of
divisors arising from embeddings of ``smaller'' orthogonal groups, the
so-called Hirzebruch-Zagier divisors. Their study led
to important discoveries and triggered generalizations in various
directions.  Moreover, the theta correspondence provides a source of
automorphic forms related to the geometry of Hirzebruch-Zagier divisors.

A more recent
development is the regularized theta lifting due to Borcherds, Harvey
and Moore, which yields to automorphic products and automorphic Green
functions.  The focus of the present text is on these topics,
highlighting the role of the orthogonal group.  We
added some background material on quadratic spaces and orthogonal
groups, to make the connection  explicit.

I thank G. van der Geer and D. Zagier for several interesting
conversations during the summer school at the Sophus Lie Conference
Center.  Moreover, I thank J. Funke for his helpful comments on a
first draft of this manuscript.

\section{Hilbert modular surfaces}

In this section we give a brief introduction to Hilbert modular
surfaces associated to real quadratic fields. For details we refer
to \cite{Fr}, \cite{Ga}, \cite{Ge}, \cite{Go}.

\subsection{The Hilbert modular group} \label{sect:1}

Let $d>1$ be a squarefree integer. Then $F=\IQ(\sqrt{d})$ is a real
quadratic field, which we view as a subfield of $\R$. The
discriminant of $F$ is
\begin{align}
D=\begin{cases} d,& \text{if $d\equiv 1\pmod{4}$,}\\
4d,& \text{if $d\equiv 2,3 \pmod{4}$.}
\end{cases}
\end{align}
We write $\calO_F$ for the ring of integers in $F$, so
\begin{align}
\calO_F= \begin{cases}
\Z+\frac{1+\sqrt{d}}{2}\Z,&  \text{if $d\equiv 1\pmod{4}$,}\\
\Z+\sqrt{d}\Z,& \text{if $d\equiv 2,3 \pmod{4}$.}
\end{cases}
\end{align}
The ring $\calO_F$ is a Dedekind domain, that is, it is a noetherian
integrally closed integral domain in which every non-zero prime
ideal is maximal.

We denote by $\calO_F^*$ the group of units in $\calO_F$. By the
Dirichlet unit theorem there is a unique unit $\eps_0>1$ such that
$\calO_F^*=\{ \pm 1 \}\times\{ \eps_0^n ;\; n\in\Z\}$. It is called
the fundamental unit of $F$. We write $x\mapsto x'$ for the
conjugation, $\norm(x)=xx'$ for the norm in $F$, and $\tr(x)=x+x'$
for the trace in $F$. The different of $F$ is denoted by $\frakd_F$.
Note that $\frakd_F=(\sqrt{D})$.

Recall that an (integral) ideal of $\calO_F$ is a
$\calO_F$-submodule of $\calO_F$. A fractional ideal of $F$ is a
finitely generated $\calO_F$-submodule of $F$.  Fractional ideals
form a group together with the ideal multiplication.  The neutral
element is $\calO_F$ and the inverse of a fractional ideal
$\fraka\subset F$ is
\[
\fraka^{-1}= \{ x\in F;\; x\fraka\subset\calO_F\}.
\]
Since $F$ is a quadratic extension of $\Q$, we have the useful
formula $\fraka^{-1}=\frac{1}{\norm(\fraka)}\fraka'$, where
$\fraka'$ is the conjugate of $\fraka$.  Two fractional ideals
$\fraka,\frakb$ are called equivalent, if there is a $r\in F$ such
that $\fraka=r\frakb$. The group of equivalence classes $\Cl(F)$ is
called the ideal class group of $F$. It is a finite abelian group.
Two fractional ideals $\fraka,\frakb$ are called equivalent in the
narrow sense, if there is a totally positive $r\in F$ such that
$\fraka=r\frakb$. The group of equivalence classes $\Cl^+(F)$ is
called the narrow ideal class group of $F$. It is equal to $\Cl(F)$,
if and only if $\eps_0$ has norm $-1$. Otherwise it is an extension
of degree $2$ of $\Cl(F)$. The (narrow) class number of $F$ is the
order of the (narrow) ideal class group. It measures how far
$\calO_F$ is from being a principal ideal domain.

If the class number of $F$ is greater than $1$, there are ideals
which cannot be generated by a single element. However, we have the
following fact, which holds in any Dedekind ring.

\begin{remark}\label{rem_ideal}
  If $\fraka\subset F$ is a fractional ideal, then there exist $\alpha
  ,\beta\in F$ such that $\fraka=\alpha\calO_F+\beta\calO_F$.  \hfill
  $\square$
\end{remark}

The group $\Sl_2(F)$ is embedded into $\Sl_2(\R)\times\Sl_2(\R)$ by
the two real embeddings of $F$. It acts on $\H\times\H$, where
$\H=\{\tau\in\C;\;\Im(\tau)>0\}$ is the complex upper half plane,
via fractional linear transformations,
\begin{align}
  \abcd z=\left( \frac{az_1 +b}{cz_1 +d}, \frac{a'z_2 +
      b'}{c'z_2 +d'} \right).
\end{align}
Here and throughout we use $z=(z_1,z_2)$ as a standard variable on
$\H^2$. If $z\in \H^2$ and $\kabcd\in \Sl_2(F)$, we write
\begin{align}
\norm(cz+d)=(cz_1+d)(c'z_2+d').
\end{align}

\begin{lemma}\label{lem:impart}
For $z\in \H^2$ and $\gamma=\kabcd\in \Sl_2(F)$ we have
\[
\Im(\gamma z) =  \frac{\Im(z)}{|\norm(cz+d)|^2}.
\]
\end{lemma}

\begin{proof}
This follows immediately from the analogous assertion in the
$1$-dimensional case.
\end{proof}

If $\fraka$ is a fractional ideal of $F$, we write
\begin{align}\label{def:fraka}
\Gamma(\calO_F\oplus\fraka)=\left\{ \abcd\in \Sl_2(F);\quad a,d\in \calO_F,\;b\in
  \fraka^{-1},\; c\in \fraka\right\}
\end{align}
for the {\em Hilbert modular group} corresponding to $\fraka$.
Moreover, we write
\begin{align}
\Gamma_F=\Gamma(\calO_F\oplus\calO_F)=\Sl_2(\calO_F).
\end{align}

Let $\Gamma\subset \Sl_2(F)$ be a subgroup which is commensurable with
$\Gamma_F$, i.e., $\Gamma\cap\Gamma_F$ has finite index in both,
$\Gamma$ and $\Gamma_F$. Then $\Gamma$ acts properly discontinuously
on $\H^2$, i.e., if $W\subset\H^2$ is compact, then $\{
\gamma\in\Gamma;\;\gamma W\cap W\neq\emptyset \}$ is finite (see
Corollary \ref{pd}). In particular, for any $a\in \H^2$, the
stabilizer $\Gamma_a=\{\gamma\in \Gamma;\; \gamma a = a\}$ is a finite
subgroup of $\Gamma$.  Let $\bar \Gamma_a$ be the image of $\Gamma_a$
in $\PSL_2(F)=\Sl_2(F)/\{\pm1\}$.  If $\#\bar\Gamma_a>1$ then $a$ is
called an {\em elliptic fixed point} for $\Gamma$ and $\#
\bar\Gamma_a$ is called the order of $a$. The order of $a$ only
depends of the $\Gamma$-class. Moreover, there are only finitely many
$\Gamma$-classes of elliptic fixed points. It can be shown that
$\Gamma$ always has a finite index subgroup which has no elliptic
fixed points.

The quotient
\begin{align}
Y(\Gamma)=\Gamma\backslash \H^2
\end{align}
is a normal complex surface. The singularities are given by the
elliptic fixed points. They are finite quotient singularities.

The surface $Y(\Gamma)$ is non-compact. It can be compactified by
adding a finite number of points, the cusps of $\Gamma$. They can be
described as follows.
The group $\Sl_2(F)$ also acts on $\P^1(F)=F\cup\{ \infty \}$ by
\[
\abcd\frac{\alpha}{\beta}=\frac{a\frac{\alpha}{\beta}+b}{c\frac{\alpha}{\beta}+d}=\frac{a\alpha+b\beta}{c\alpha+d\beta}.
\]
Notice that, since $\kabcd\infty = \frac{a}{c}$, the action of
$\Sl_2(F)$ is transitive. The $\Gamma$-classes of $\P^1(F)$ are
called the {\em cusps} of $\Gamma$.

\begin{lemma}
The map
\begin{eqnarray*}
\varphi :\Gamma_F\bs\IP^1(F) &\xymatrix{\ar[r] & }& \Cl(F), \\
(\alpha:\beta) &\xymatrix{ \ar@{|->}[r] & }&
\alpha\calO_F+\beta\calO_F,
\end{eqnarray*}
is bijective.
\end{lemma}

\begin{proof}
We first show that $\varphi$ is well-defined: It is clear that
$\varphi (\alpha:\beta)=\varphi (r\alpha:r\beta)$. Now let
$\kabcd\in\Gamma_F$, and let
$\frac{\gamma}{\delta}=\kabcd\frac{\alpha}{\beta}$. We need to show
that $\varphi (\gamma:\delta)=\varphi (\alpha:\beta)$. We have
\begin{align*}
\varphi\left(\gamma:\delta\right) &= \gamma \calO_F+ \delta\calO_F \\
&= (a\alpha + b\beta)\calO_F + (c\alpha + d\beta)\calO_F \\
&\subset \varphi\left(\alpha:\beta\right) .
\end{align*}
Interchanging the roles of $(\gamma:\delta)$ and $(\alpha:\beta)$,
we see
\begin{align*}
\varphi\left(\alpha:\beta \right)
&= (d\gamma - b\delta)\calO_F + (-c\gamma + a\delta)\calO_F \\
&\subset \varphi\left(\gamma:\delta\right).
\end{align*}
Consequently, $\varphi (\gamma:\delta)=\varphi (\alpha:\beta)$.

The surjectivity of 
$\varphi$ follows from Remark \ref{rem_ideal}.

Finally, we show that $\varphi$ is injective. Let $\fraka=\varphi
(\alpha:\beta)= \varphi (\gamma:\delta)$. Then $1 \in\calO_F
=\fraka\fraka^{-1} =\alpha\fraka^{-1} +\beta\fraka^{-1}$. So there
exist $\tilde{\alpha},\tilde{\beta}\in\fraka^{-1}$ such that $
1=\alpha\tilde{\beta} - \beta\tilde{\alpha}$. We find that
\begin{align*}
M_1: =
\zxz{\alpha}{\tilde{\alpha}}{\beta}{\tilde{\beta}}\in\zxz{\fraka}{\fraka^{-1}}{\fraka}{\fraka^{-1}}\cap\Sl_2(F),
\end{align*}
and $M_1\infty=(\alpha:\beta)$. In the same way we find
\begin{align*}
M_2:=\zxz{\gamma}{\tilde{\gamma}}{\delta}{\tilde{\delta}}\in\zxz{\fraka}{\fraka^{-1}}{\fraka}{\fraka^{-1}}\cap\Sl_2(F)
\end{align*}
such that $M_2\infty=(\gamma:\delta)$. Therefore we have
\begin{align*}
M_2 M_1^{-1}
\in\zxz{\fraka}{\fraka^{-1}}{\fraka}{\fraka^{-1}}\zxz{\fraka^{-1}}{\fraka^{-1}}{\fraka}{\fraka}
=\zxz{\calO_F}{\calO_F}{\calO_F}{\calO_F}.
\end{align*}
Hence, $M_2 M_1^{-1}\in\Gamma_F$ and $M_2
M_1^{-1}(\alpha:\beta)=(\gamma:\delta)$.  This concludes the proof
of the Lemma.
\end{proof}

\begin{corollary}
  The number of cusps of $\Gamma_F$ is equal to the class number
  $h(F)$ of $F$.  A subgroup $\Gamma\subset \Sl_2(F)$ which is
  commensurable with $\Gamma_F$ has finitely many cusps.
\end{corollary}

\begin{remark}
  Let $\Gamma_{\infty}\subset\Gamma$ be the stabilizer of $\infty$.
  Then there is a $\Z$-module $M\subset F$ of rank $2$ and a finite
  index subgroup $V\subset \calO_F^*$ acting on $M$ such that the
  group
\begin{align}
\label{def:gmv} G(M,V)=\left\{ \zxz{\eps}{\mu}{0}{\eps^{-1}};\;
\text{$\mu\in M$ and $\eps\in V$}\right\}
\end{align}
is contained in $\Gamma_\infty$ with finite index.
\end{remark}

\begin{example}
If $\fraka\subset F$ is a fractional ideal, then
\[
\Gamma(\calO_F\oplus\fraka)_\infty= \left\{ \zxz{\eps}{\mu}{0}{\eps^{-1}};\;
\mu\in
  \fraka^{-1},\, \eps\in \calO_F^*\right\}.
\]
\end{example}

\subsection{The Baily-Borel compactification}

We embed $\P^1(F)$ into $\P^1(\R)\times\P^1(\R)$ via the two real
embeddings of $F$. Then we may view $\P^1(F)$ as the set of rational
boundary points of $\H^2$ in the same way as $\P^1(\Q)$ is viewed as
the set of rational boundary points of $\H$. Here we consider
\begin{align}
(\IH^2)^{\ast}=\IH^2\cup\IP^1(F).
\end{align}
By introducing a suitable topology on $(\IH^2)^{\ast}$, the quotient
$\Gamma\bs (\IH^2)^{\ast}$ can be made into a compact Hausdorff
space. This leads to the Baily-Borel compactification of
$Y(\Gamma)$.

\begin{proposition}
\label{prop:bbtop} On $(\IH^2)^{\ast}$ there is a unique topology
with the following properties:
\begin{enumerate}
\item[(i)] 
The induced topology on $\IH^2$ agrees with the usual one.
\item[(ii)]  
$\IH^2$ is open and dense in $(\IH^2)^{\ast}$.
\item[(iii)] 
The sets $U_C\cup\infty$, where
\[
U_C=\left\{ (z_1,z_2)\in\IH^2 ;\;  \Im(z_1)\Im(z_2) >C \right\}
\]
for $C>0$, form a base of open neighborhoods of the cusp $\infty$.
\item[(iv)] 
If $\kappa\in\IP^1(F)$ and $\rho\in\Sl_2(F)$ with
$\rho \infty=\kappa$, then the sets
\[
\rho(U_C\cup\infty) \qquad (C>0)
\]
form a base of open neighborhoods of the cusp $\kappa$.
\end{enumerate}
\hfill $\square$
\end{proposition}

\begin{remark}
  The system of open neighborhoods of $\kappa$ defined by (iv) does not
  depend on the choice of $\rho$. The stabilizer $\Gamma_{\infty}$ of
  $\infty$ acts on $U_C$.  If
  $\gamma=\kzxz{\eps}{\mu}{0}{\eps^{-1}}\in \Gamma_{\infty}$, then
\[
\gamma z=(\eps^2 z_1+\eps \mu,\eps'{}^2 z_2+\eps' \mu') .
\]
\end{remark}

We consider the quotient
\begin{align}
X(\Gamma)=\Gamma\bs(\IH^2)^{\ast}.
\end{align}

\begin{theorem}
The quotient space $X(\Gamma)$, together with the quotient topology,
is a compact Hausdorff space. \hfill $\square$
\end{theorem}

\begin{proposition}
\label{prop:cuspinfty} For $C>0$ sufficiently large, the canonical
map
\[
\xymatrix{\Gamma_{\infty}\bs U_C\cup\infty\ar[r] & \Gamma\bs
(\IH^2)^{\ast} }
\]
is an open embedding. \hfill $\square$
\end{proposition}

The group $\Sl_2(F)$ acts by topological automorphisms on $(\H^2)^*$.
Hence, for $\rho\in \Sl_2(F)$, the natural map
\[
X(\Gamma) \longrightarrow X(\rho^{-1}\Gamma \rho),\quad z\mapsto
\rho^{-1} z
\]
is topological. If $\rho\infty =\kappa$, it takes the cusp $\kappa$
of $\Gamma$ to the cusp $\infty$ of $\rho^{-1} \Gamma \rho$. In that
way, local considerations near the cusps can often be reduced to
considerations at the cusp $\infty$ (for a conjugate group), for
which one can use Proposition \ref{prop:cuspinfty}.

We define a complex structure on $X(\Gamma)$ as follows. For
$V\subset X(\Gamma)$ open we let $U\subset (\H^2)^*$ be the inverse
image under the canonical projection $\pr:  (\H^2)^*\to X(\Gamma)$,
and let $U'$ be the inverse image in $\H^2$. We have the diagram
\[
\xymatrix{
\IH^2 \ar[r] & (\IH^2)^{*} \ar[r]^{\pr} & X(\Gamma)\\
U'\ar[u]\ar[r] & U \ar[u]\ar[r] & V \ar[u] }.
\]
We define $\calO_{X(\Gamma)}(V)$ to be the ring of continuous
functions $f:V\to \C$ such that $\pr^*(f)$ is holomorphic on $U'$.
This defines a sheaf $\calO_{X(\Gamma)}$ of rings on $X(\Gamma)$,
and the pair $(X(\Gamma), \calO_{X(\Gamma)})$ is a locally ringed
space.

\begin{theorem}[Baily-Borel]
The space $(X(\Gamma),\calO_{X(\Gamma)})$ is a normal complex space.
\hfill $\square$
\end{theorem}

The proof is based on a criterion of Baily and Cartan for the
continuation of complex structures, see \cite{Fr} p.~112.

In contrast to the case of modular curves the resulting normal
complex space $X(\Gamma)$ is not regular. There are finite quotient
singularities at the elliptic fixed points, and more seriously, the
cusps are highly singular points.  By the theory of Hironaka the
singularities can be resolved \cite{Hi}. A weak form of Hironaka's
result states that there exists a desingularization
\begin{align}
\label{desing} \pi:\widetilde{X}(\Gamma)\longrightarrow X(\Gamma),
\end{align}
where $\widetilde{X}(\Gamma)$ is a non-singular connected projective
variety such that $\pi$ induces a biholomorphic map
$\pi^{-1}(X(\Gamma)^{reg})\to X(\Gamma)^{reg}$.  Here
$X(\Gamma)^{reg}$ is the regular locus of $X(\Gamma)$.  One can
further require that the complement of $\pi^{-1}(X(\Gamma)^{reg})$
is a divisor with normal crossings.  The minimal resolution of
singularities was constructed by Hirzebruch \cite{Hirzebruch}.

It can be shown that there is an ample line bundle on $X(\Gamma)$,
the line bundle of modular forms (in sufficiently divisible weight).
Consequently, the space $X(\Gamma)$ carries the structure of a
projective algebraic variety over $\C$.  The surface $Y(\Gamma)$ is
a Zariski-open subvariety and therefore quasi-projective.

\begin{remark}
  The Hilbert modular surfaces $Y(\Gamma)$ often have a moduli
  interpretation, analogously to the fact that $\Sl_2(\Z)\bs\H$
  parametrizes isomorphism classes of elliptic curves over $\C$. It
  can be used to construct integral models of Hilbert modular
  surfaces.  For instance, $Y(\Gamma(\calO_F\oplus\fraka))$ is the coarse moduli
  space for isomorphism classes of triples $(A,\imath,m)$, where $A$
  is an abelian surface over $\C$, and $\imath:\calO_F\to \End(A)$ is
  an embedding of rings (real multiplication), and $m$ is an
  isomorphism from the polarization module of $A$ to
  $(\fraka\frakd_F)^{-1}$ respecting the positivities,
cf.~\cite{Go} Theorem 2.17. The variety
  $Y(\Gamma(\calO_F\oplus\fraka))$ can be interpreted as the complex points of a
  moduli stack over $\Z$. One can also
  construct toroidal compactifications and Baily-Borel compactifications
  over $\Z$, cf.~\cite{Ra}, \cite{DePa}, \cite{Ch}.
\end{remark}

\subsubsection{Siegel domains}

Here we recall the properties of Siegel domains for Hilbert modular
surfaces. They are nice substitutes for fundamental domains.

We write $(x_1,x_2)$ for the real part and $(y_1,y_2)$ for the
imaginary part of $(z_1,z_2)$. The top degree differential form
\begin{align}
d\mu=\frac{dx_1\,dy_1}{y_1^2}\,\frac{dx_2\,dy_2}{y_2^2}
\end{align}
on $\H^2$ is invariant under the action of $\Sl_2(\R)^2$. It defines
an invariant measure on $\H^2$, which is induced by the Haar measure
on $\Sl_2(\R)$.

\begin{definition}
A subset $S\subset\IH^2$ is called a {\em fundamental set} for
$\Gamma$, if
\[
\IH^2=\bigcup_{\gamma\in\Gamma} \gamma(S).
\]
\end{definition}

\begin{definition}
A fundamental set $S$ for $\Gamma$ is called a {\em fundamental
domain} for $\Gamma$, if
\begin{enumerate}
\item[(i)] 
$S$ is measurable.
\item[(ii)] 
There is a subset $N\subset S$ of measure zero,
such that for all $z, w\in S\bs N$ we have
\begin{align*}
z \sim_{\Gamma} w \  \Rightarrow z=w.
\end{align*}
\end{enumerate}
\end{definition}

\begin{remark}
  It can be shown that every measurable fundamental set contains a
  fundamental domain.
\end{remark}

Nice fundamental sets for the action of $\Gamma$ on $\IH^2$ are
given by {\em Siegel domains}: For a positive real number $t$ we put
\begin{align}
\label{siegeldom} \calS_t= \left\{ z\in\IH^2 ;\quad
\text{$|\,x_i|<t$ and $y_i >\tfrac{1}{t}$ for $i=1,2$} \right\}.
\end{align}

\begin{proposition}\label{siegel1}
For any fixed $t\in \R_{>0}$ there exist only finitely many
$\gamma\in \Gamma$ such that
\begin{align}
\label{siegelcond} \gamma \calS_t\cap \calS_t \neq \emptyset.
\end{align}
\end{proposition}

\begin{proof}
It is clear that there are only finitely many $\gamma=\kabcd\in
\Gamma$ with $c=0$ satisfying condition \eqref{siegelcond}.

On the other hand, assume that $\gamma\in \Gamma$ with $c\neq 0$,
and assume that there is a $z\in \gamma \calS_t\cap\calS_t$. Then we
have
\begin{align}
\label{in1}
\frac{y_1}{|c z_1 + d|^2}&>\frac{1}{t},\\
\label{in2} \frac{y_2}{|c' z_2 + d'|^2}&>\frac{1}{t}.
\end{align}
The first inequality implies that
\begin{align}
\label{in3} y_1 > \frac{1}{t}\left((cx_1+d)^2 + c^2 y_1^2\right)
\geq \frac{1}{t} c^2 y_1^2 > \frac{1}{t^2} c^2 y_1,
\end{align}
and therefore $|c|<t$. In the same way, inequality \eqref{in2}
implies that $|c'|<t$. Hence there are only finitely many
possibilities for $c$. For these, by  \eqref{in3} and its analogue
for $y_2$, the imaginary part $(y_1,y_2)$ is bounded, and there are
also just finitely many possibilities for $d$.

Moreover, replacing $\gamma$ by $\gamma^{-1}$ in the above argument,
we find that there are only finitely many possibilities for $a$. But
$a,c,d$ determine $\gamma$.
\end{proof}

\begin{corollary}\label{pd}
  The action of $\Gamma$ on $\H^2$ is properly discontinuous, that is,
  if $W\subset\H^2$ is compact, then $\{ \gamma\in\Gamma;\;\gamma
  W\cap W\neq\emptyset \}$ is finite.
\end{corollary}

\begin{proof}
This follows from Proposition \ref{siegel1} and the fact that
$\bigcup_{t\in \R_{>0}} \calS_t =\H^2$.
\end{proof}

\begin{theorem}\label{fundset}
Let $\kappa_1,\ldots,\kappa_r\in\IP^1(F)$ be a set of
representatives for the cusps of $\Gamma$, and let
$\rho_1,\ldots,\rho_r\in\Sl_2(F)$ such that
$\rho_j\,\infty=\kappa_j$. There is a $t>0$ such that
\[
\calS=\bigcup_{j=1}^{r} \rho_j\calS_t
\]
is a measurable fundamental set for $\Gamma$.
\end{theorem}

\begin{proof}
See \cite{Ga}, Chapter 1.6.
\end{proof}


\subsection{Hilbert modular forms}

Let $\Gamma\subset\Sl_2(F)$ be a subgroup which is commensurable
with $\Gamma_F$, and let $(k_1,k_2)\in\IZ^2$.

\begin{definition}
\label{hilbertmodular} A meromorphic function $f:\IH^2\to \IC$ is
called a meromorphic {\em Hilbert modular form} of weight
$(k_1,k_2)$ for $\Gamma$, if
\begin{align}
f(\gamma z)=(cz_1 +d)^{k_1} (c'z_2 +d')^{k_2} f(z) \label{formel}
\end{align}
for all $\gamma=\kabcd\in\Gamma$. If $k_1=k_2=:k$, then $f$ is said
to have {\em parallel weight}, and is simply called a meromorphic
Hilbert modular form of weight $k$. If $f$ is holomorphic on $\H^2$,
then it is called a holomorphic Hilbert modular form. Finally, a
Hilbert modular form $f$ is called {\em symmetric}, if
$f(z_1,z_2)=f(z_2,z_1)$.
\end{definition}

For a function $f: \IH^2 \to\IC$ and $\gamma=\kabcd\in\Sl_2(F)$ we
define the Petersson slash operator by
\begin{align*}
(f\mid_{k_1,k_2} \gamma)(z)=(cz_1 +d)^{-k_1} (c'z_2
+d')^{-k_2}f(\gamma z).
\end{align*}
The assignment $f\mapsto f\mid_{k_1,k_2} \gamma$ defines a right
action of $\Sl_2(F)$ on complex valued functions on $\H^2$. Using
it, we may rewrite condition \eqref{formel} as
\begin{align*}
f\mid_{k_1,k_2} \gamma = f, \qquad \gamma\in  \Gamma.
\end{align*}
If $k_1=k_2=:k$, we simply write $f\mid_k \gamma$ instead of
$f\mid_{k_1,k_2}\gamma$.

If $f$ is a holomorphic Hilbert modular form for $\Gamma$, it has a
Fourier expansion at the cusp $\infty$ of the following form. Let
$M\subset F$ be a $\Z$-module of rank $2$ and let $V\subset
\calO_F^*$ be a finite index subgroup $V\subset \calO_F^*$ acting on
$M$ such that the group
\begin{align}\label{def:gmv2}
G(M,V)=\left\{ \zxz{\eps}{\mu}{0}{\eps^{-1}};\; \text{$\mu\in M$ and
$\eps\in V$}\right\}
\end{align}
is contained in $\Gamma_\infty$ with finite index.
%
The transformation law \eqref{formel} for $\gamma\in G(M,V)$ implies
that
\begin{align*}
f(z+\mu)=f(z)
\end{align*}
for all $\mu\in M$. Therefore, $f$ has a normally convergent Fourier
expansion
\begin{align}\label{def:fexp}
f=\sum_{\nu\in M^{\vee}} a_{\nu} \,e(\tr(\nu z)),
\end{align}
where $e(w)=e^{2\pi iw}$, $\tr(\nu z)=\nu z_1+\nu'z_2$,  and
\begin{align}
M^{\vee}=\{ \lambda\in F ;\; \text{$\tr(\mu\lambda)\in\IZ$ for all
$\mu\in M$} \}
\end{align}
is the dual lattice of $M$ with respect to the trace form on $F$.
The Fourier coefficients $a_\nu$ are given by
\begin{align}\label{eq:coeff}
a_{\nu}=\frac{1}{\vol(\IR^2/M)} \int\limits_{\IR^2/M} f(z)
e(-\tr(\nu z))\,dx_1\,dx_2.
\end{align}

More generally, if $\kappa\in\IP^1(F)$ is any cusp of $\Gamma$, we
take $\rho\in\Sl_2(F)$ such that $\rho\infty=\kappa$, and consider
$f\mid_{k_1,k_2}\rho$ at $\infty$. The Fourier expansion of $f$ at
$\kappa$ is the expansion of $f\mid_{k_1,k_2}\rho$ at $\infty$ (it
depends on the choice of $\rho$).

It is a striking fact that, in contrast to the one-dimensional
situation, a holomorphic Hilbert modular form is automatically
holomorphic at the cusps by the G\"otzky-Koecher principle.

\begin{theorem}[G\"otzky-Koecher principle]\label{koecherprinc}
Let $f:\IH^2\rightarrow \IC$ be a holomorphic function satisfying
$f\mid_{k_1,k_2} \gamma=f$ for all $\gamma\in G(M,V)$ as in
\eqref{def:gmv2}. Then
\begin{enumerate}
\item[(i)] 
$a_{\eps^2\nu}=\eps^{k_1}\eps'{}^{k_2} a_{\nu} $ for all $\nu\in M^{\vee}$
and $\eps\in V$,
\item[(ii)] 
$a_{\nu}\neq 0\ \Rightarrow\ \nu=0$ or $\nu\gg 0$.
\end{enumerate}
\end{theorem}

\begin{proof}
(i) For $\eps\in V$ we have $\kzxz{\eps^{-1}}{0}{0}{\eps}\in
G(M,V)$.
The transformation law implies that
\begin{align*}
\eps^{-k_1}\eps'{}^{-k_2} \sum_{\nu\in M^{\vee}}
a_{\nu}\,e(\tr(\nu\eps^{-2}z))=\sum_{\nu\in M^{\vee}}
a_{\nu}\,e(\tr(\nu z)).
\end{align*}
Comparing Fourier coefficients, we obtain the first assertion.

(ii) Suppose that there  is a $\nu\in M^\vee$ such that $a_\nu\neq 0$
and such that $\nu<0$ or $\nu'<0$. Without loss of generality we
assume $\nu<0$. There is an $\eps\in V$ with $\eps>1$ and
$0<\eps'<1$ such that $\tr(\eps\nu)<0$. Then $\tr(\eps^{2n}\nu)$
goes to $-\infty$ for $n\to \infty$.

The series
\begin{align*}
\sum_{n\geq 1} a_{\nu\eps^{2n}}\,e(i\tr(\nu\eps^{2n}))
\end{align*}
is a subseries of the Fourier expansion of $f(z)$ at $z=(i,i)$ and
therefore converges absolutely. But by (i) we have
\begin{align*}
\sum_{n\geq 1}|a_{\nu\eps^{2n}}\,e(i\tr(\nu\eps^{2n}))| =
|a_{\nu}|\sum_{n\geq 1}\eps^{k_1}\eps'{}^{k_2}
e^{-2\pi\tr(\nu\eps^{2n})} \ \rightarrow \ \infty,
\end{align*}
contradicting the convergence.
\end{proof}

\begin{corollary}
A holomorphic Hilbert modular form for the group $\Gamma$ has a
Fourier expansion at the cusp $\infty$ of the form
\begin{align}
\label{hholmf}
f(z)=a_0 +\sum_{\substack{\nu\in M^{\vee} \\ \nu\gg 0}}
a_{\nu}\,e(\tr(\nu z)).
\end{align}
\end{corollary}

The constant term $a_0$ is called the value of $f$ at $\infty$. We
write $f(\infty)=a_0$. More generally, if $\kappa\in\IP^1(F)$ is
any cusp of $\Gamma$, we take $\rho\in\Sl_2(F)$ such that
$\rho\infty=\kappa$. We put
$f(\kappa)=(f\mid_{k_1,k_2}\rho)(\infty)$. If $(k_1,k_2)\neq (0,0)$,
the value $f(\kappa)$ of $f$ at $\kappa$ depends on the choice of $\rho$ (by a
non-zero factor).


\begin{definition}
A holomorphic Hilbert modular form $f$ is called a {\em cusp form},
if it vanishes at all cusps of $\Gamma$.
\end{definition}

\begin{proposition}\label{hcusp}
Let $f$ be a holomorphic 
modular form of weight $(k_1,k_2)$ for $\Gamma$.
If $k_1\neq k_2$, then $f$ is a cusp form.
\end{proposition}

\begin{proof}
  This follows from Theorem \ref{koecherprinc} (i), applied to the
  constant terms at the cusps.
\end{proof}

\begin{proposition}\label{hinv}
Let $f$ be a modular form of weight $(k_1,k_2)$ for $\Gamma$. Then the
function $h(z)=|f(z) y_1^{k_1/2}y_2^{k_2/2}|$ is $\Gamma$-invariant.
\end{proposition}

\begin{proof}
This follows from Lemma \ref{lem:impart}.
\end{proof}

\begin{proposition}\label{hbound}
  Let $f$ be a holomorphic modular form of weight $(k_1,k_2)$ for $\Gamma$, and
  let $h(z)=|f(z) y_1^{k_1/2}y_2^{k_2/2}|$ be the $\Gamma$-invariant
  function of Proposition \ref{hinv}.  
\begin{enumerate}
\item[(i)] If $f$ has parallel negative weight $k=k_1=k_2$, 
then $h$ 
attains a maximum on $\H^2$.
\item[(ii)]
If $f$ is a cusp form, then $h$ vanishes at the cusps and 
attains a maximum on $\H^2$.
\end{enumerate}
\end{proposition}

\begin{proof}
  We only prove the first statement, the second is similar.  By
  Proposition \ref{hinv}, it suffices to consider $h$ on a fundamental
  set for $\Gamma$.  In view of Theorem \ref{fundset}, it suffices to
  show that for any $\rho\in \Sl_2(F)$ and any $t\in \R_{>0}$, 
  the function $h(\rho z)$ is
  bounded and attains its maximum on the Siegel domain $\calS_t$.
  Using the Fourier expansion of $f$ at the cusp $\rho\infty$, we see that
\[
h(\rho z)=(y_1 y_2)^{k/2}a_0 + (y_1 y_2)^{k/2}\sum_{\substack{\nu\in M^{\vee} \\ \nu\gg 0}}
a_{\nu}\,e(\tr(\nu z))
\]
for a suitable 
rank $2$ lattice $M\subset F$.
Since $k<0$ we find that  $\lim_{y_1y_2\to \infty} h(\rho z)=0$ on $\calS_t$.
Consequently, $h(\rho z)$ is bounded and attains a maximum on $\calS_t$.
\end{proof}

\begin{proposition}\label{weightzero}
Let $f$ be a holomorphic modular form of weight $(k_1,k_2)$ for
$\Gamma$. Then $f$ vanishes identically unless $k_1,k_2$ are both
positive or $k_1=k_2=0$. In the latter case $f$ is constant.
\end{proposition}

\begin{proof}
  Let us first consider the case that $k_1=0$ and $k_2\neq 0$.  
  By Proposition \ref{hcusp}, $f$ is a cusp form.  The function
  $h(z)=y_2^{k_2/2} f(z)$ is holomorphic in $z_1$, and, by Proposition
  \ref{hbound}, $|h|$ attains a maximum on $\H^2$.
  According to the maximum principle, $h$ must be constant
as a function of $z_1$.  Hence,
\[
h(z_1,z_2)=h(z_1,\gamma' z_2)
\]
for all $\gamma\in \Gamma$.  Since $\{\gamma z_2;\;\gamma\in
\Gamma\}$ is dense in $\H$, the function $h$ must be constant on
$\H^2$. Because $h$ vanishes at the cusps, it must vanish
identically.
In the same way we see that $f=0$, if $k_2=0$ and $k_1\neq 0$.

Let us now consider the case that 
$k_1=k_2=0$. If $f$ is a cusp form,
then Proposition~\ref{hbound} implies that $|f|$ attains a maximum
on $\H^2$.  Hence, by the maximum principle, $f$ must be constant.
Since $f$ vanishes at $\infty$, we obtain that $f\equiv 0$. If $f$
is holomorphic (but not necessarily cuspidal), we consider the cusp
form
\begin{align*}
g(z):=\prod_{\kappa\in\Gamma\bs\IP^1(F)} (f(z)-f(\kappa)).
\end{align*}
We find that $g\equiv 0$, and therefore $f$ is constant.

Finally, assume that some $k_i$, say $k_1$, is negative. 
If $k_1\neq k_2$, then, by Proposition~\ref{hcusp}, $f$ is a cusp form.
If $k_1= k_2$, then $f$ has parallel negative weight.
In both cases, Proposition~\ref{hbound} implies that 
$h(z)=|f(z) y_1^{k_1/2}y_2^{k_2/2}|$ is bounded by a constant $C>0$ on $\H^2$.
We consider the Fourier expansion of $f$ at the cusp $\infty$ as in 
\eqref{def:fexp}.
The coefficients $a_\nu$ are given by \eqref{eq:coeff}. We find that
\begin{align*}
|a_\nu|&\leq \frac{1}{\vol(\IR^2/M)} \int\limits_{\IR^2/M} |f(z)
e(-\tr(\nu z))|\,dx_1\,dx_2\\
&\leq C y_1^{-k_1/2}y_2^{-k_2/2}e^{-2\pi \tr(\nu y)}.
\end{align*}
Taking the limit $y_1\to 0$, we see that $a_\nu$ vanishes for all 
$\nu\in M^\vee$, and therefore $f\equiv 0$.
%
%
\end{proof}

\begin{corollary}
Let $\widetilde{X}(\Gamma)\to X(\Gamma)$ be a desingularization as
in \eqref{desing}. Then any holomorphic $1$-form on
$\widetilde{X}(\Gamma)$ vanishes identically.
\end{corollary}

\begin{proof}
Let $\omega$ be a holomorphic $1$-form on $\widetilde{X}(\Gamma)$
and denote by $\eta$ its pullback to the regular locus
of $X(\Gamma)$. Viewing $\eta$ as a $\Gamma$-invariant holomorphic
$1$-form on $\H^2$, we may write
\[
\eta = f_1(z) \,dz_1 + f_2(z) \, dz_2,
\]
where $f_1$ and $f_2$ are holomorphic Hilbert modular forms of
weights $(2,0)$ and $(0,2)$, respectively. Hence, by Proposition
\ref{weightzero}, $\eta$ vanishes identically.
\end{proof}

In the same way one sees that any antiholomorphic $1$-form on
$\widetilde{X}(\Gamma)$ vanishes identically. Consequently,
\[
H^1(\widetilde{X}(\Gamma),\calO_{\widetilde{X}(\Gamma)})=0,
\]
that is, the surface $\widetilde{X}(\Gamma)$ is regular.  Using
Hodge theory we see that $H^1(\widetilde{X}(\Gamma),\C)=0$, so the
first Betti number vanishes. In particular, the interesting part of
the cohomology of a Hilbert modular surface is in degree $2$.

It can be shown that the Hilbert modular surfaces
$\widetilde{X}(\Gamma(\calO_F\oplus\fraka))$ corresponding to the groups
$\Gamma(\calO_F\oplus\fraka)$ are simply connected. This also implies the
vanishing of the first Betti number. However, there are also
examples of Hilbert modular surfaces which are not simply connected.
See \cite{Ge}, Chapter IV.6. (Recall that the fundamental group of a
complex surface is a birational invariant.)

For the rest of these notes we will only be considering Hilbert
modular forms of parallel weight $k$.

\begin{notation}
Let $k\in \Z$. We write $M_k(\Gamma)$ for the $\IC$-vector space of
holomorphic Hilbert modular forms of weight $k$ for the group
$\Gamma$, and denote by $S_k(\Gamma)$ the subspace  of cusp forms.
\end{notation}

The codimension of $S_k(\Gamma)$ in $M_k(\Gamma)$ is clearly
bounded by the number of cusps of $\Gamma$.

\begin{proposition}
Let $\widetilde{X}(\Gamma)\to X(\Gamma)$ be a desingularization.
\begin{enumerate}
\item[(i)]
Meromorphic Hilbert modular forms of weight $0$ for $\Gamma$, can be
identified with meromorphic functions on $\widetilde{X}(\Gamma)$.
\item[(ii)]
Meromorphic Hilbert modular forms of weight $2$ for $\Gamma$, can be
identified with meromorphic $2$-forms on $\widetilde{X}(\Gamma)$.
\item[(iii)]
Hilbert cusp forms of weight $2$ for $\Gamma$, can be identified
with holomorphic $2$-forms on $\widetilde{X}(\Gamma)$.
\end{enumerate}
\end{proposition}

\begin{proof}
The first two assertions are easy. For the third assertion we refer
to \cite{Fr}, Chapter II.4.
\end{proof}

In particular, the arithmetic genus of the surface
$\widetilde{X}(\Gamma)$, that is, the Euler characteristic  of the
structure sheaf,  is given by
\[
\chi(\calO_{\widetilde{X}(\Gamma)})=\sum_{p=0}^2 (-1)^p \dim
H^p(\widetilde{X}(\Gamma),
\calO_{\widetilde{X}(\Gamma)})=1+\dim(S_2(\Gamma)).
\]

Holomorphic Hilbert modular forms can be interpreted as sections of
the {\em sheaf $\calM_k(\Gamma)$ of modular forms}, which can be
defined as follows: If we write $\pr:\H^2\to Y(\Gamma)$ for the
canonical projection, then the sections over an open subset
$U\subset\Gamma\bs \H^2$ are holomorphic functions on $\pr^{-1}(U)$,
satisfying the transformation law \eqref{formel}. By the Koecher
principle, this sheaf on $Y(\Gamma)$ extends to $X(\Gamma)$. It is a
coherent $\calO_{X(\Gamma)}$-module.
%

Let $n(\Gamma)$ denote the least common multiple of the orders of
all elliptic fixed points for $\Gamma$. When $n(\Gamma)\mid k$, then
$\calM_k(\Gamma)$ is a line bundle.
One can show that this line bundle is ample and thereby prove that
$X(\Gamma)$ is algebraic. Notice that
$\calM_{nk}(\Gamma)=\calM_{k}(\Gamma)^{\otimes n}$ for any positive
integer $n$.

\subsection{$M_k(\Gamma)$ is finite dimensional}

In this section we show that $M_k(\Gamma)$ is finite dimensional.
The argument is based on the comparison of two different norms on
the space of cusp forms. It is a rather general principle and works
in a much more general setting (cf. \cite{Fr}, Chapter I.6).

We begin by defining the Petersson scalar product on $M_k(\Gamma)$.
The top degree differential form
$d\mu=\frac{dx_1\,dy_1}{y_1^2}\,\frac{dx_2\,dy_2}{y_2^2}$
on $\H^2$ is invariant under the action of $\Sl_2(\R)^2$.

\begin{definition}
Let $f,g\in M_k(\Gamma)$  such that the product $f g$ is a cusp
form. We define the Petersson scalar product of $f$ and $g$ by
\begin{align*}
\langle f,g \rangle = \int_{\calF} f(z) \overline{g(z)}(y_1y_2)^k\,
d\mu,
\end{align*}
where $\calF$ is a fundamental domain for $\Gamma$.
\end{definition}

\begin{lemma}
For $f,g$ as above the Petersson scalar product converges absolutely
and is independent of the choice of the fundamental domain.
\end{lemma}

\begin{proof}
Arguing as in Proposition \ref{hbound}, we see that
$f(z)\overline{g(z)}\,(y_1 y_2)^k$ is invariant under $\Gamma$ and
bounded on $\IH^2$. Hence, that the integral does not depend on the
choice of $\calF$ follows from the absolute convergence using the
theorem on dominated convergence for the Lebesgue integral. To prove
the absolute convergence, it suffices to show that
\begin{align*}
\int_{\calF}d\mu<\infty.
\end{align*}
In view of Proposition \ref{fundset}, it suffices to show that
\begin{align*}
\int_{\calS_t}d\mu <\infty
\end{align*}
for all $t>0$. This follows from the fact that
$
\int_{1/t}^{\infty}\frac{dy}{y^2}<\infty$.
\end{proof}

In particular, the Petersson scalar product defines a hermitian
scalar product on $S_k(\Gamma)$. We denote the corresponding
$L^2$-norm on $S_k(\Gamma)$ by
\begin{align}\label{l2norm}
\| f \|_2 :=\sqrt{\langle f,f \rangle}.
\end{align}
On the other hand we have the maximum norm on $S_k(\Gamma)$ which is
defined by
\begin{align}\label{maxnorm}
\| f\|_{\infty} = \max_{z\in\calF} \left( |f(z)| (y_1
y_2)^{k/2}\right).
\end{align}

\begin{lemma}\label{lem_norm}
There is a constant $A=A(\Gamma,k)>0$ such that
\begin{align*}
\| f \|_{\infty}\leq A \cdot \| f \|_2
\end{align*}
for all $f\in S_k(\Gamma)$.
\end{lemma}

\begin{proof}
The $L^2$-norm can be estimated by considering the Fourier
expansions of $f$ at the cusps of $\Gamma$ and using Siegel domains
(Proposition \ref{fundset}). See \cite{Fr}, Lemma 6.2 for details.
\end{proof}

\begin{theorem}
The vector space $M_k(\Gamma)$
is finite dimensional.
\end{theorem}

\begin{proof}
It suffices to show that $\dim S_k(\Gamma)<\infty$. Let
$f_1,\ldots,f_m$ be an orthonormal set with respect to the
Petersson scalar product, that is, $\langle f_i,f_j
\rangle=\delta_{ij}$. For an arbitrary linear combination
\begin{align*}
f=\sum_{j=1}^m c_j f_j
\end{align*}
with coefficients $c_j\in\IC$, we have $\|f \|_{\infty}\leq A \| f
\|_2$ by Lemma \ref{lem_norm}. Hence, for all $z\in\IH^2$ we have
\begin{align*}
\left| \sum_{j=1}^m c_j  f_j(z)(y_1  y_2)^{k/2} \right| &\leq A
\left( \sum_{j=1}^m |c_j|^2 \right)^{1/2}.
\end{align*}
We consider the inequality for $c_j=\overline{f_j(z)}$. We find
\begin{align*}
\sum_{j=1}^m | f_j(z)|^2 (y_1 y_2)^{k/2} &\leq A \left(\sum_{j=1}^m
| f_j(z) |^2 \right)^{1/2} .
\end{align*}
Dividing by the sum on the right hand side and taking the square we
obtain
\begin{align*}
\sum_{j=1}^m |\,f_j(z)|^2 (y_1\, y_2)^k & \leq A^2 .
\end{align*}
Integrating over $\calF$ we find that
\[
m \leq A^2\vol(\Gamma\bs\IH^2).
\]
This concludes the proof of the theorem.
\end{proof}

\subsection{Eisenstein series}

\label{sect:eis}

Here we define Eisenstein series for Hilbert modular groups. For
simplicity we only consider the full Hilbert modular group
$\Gamma_F=\Sl_2(\calO_F)$ of the real quadratic field $F$. We write
$\norm(x)$ for the norm of $x\in F$, and $\norm(\fraka)$ for the
norm of a fractional ideal $\fraka\subset F$.

Let $B\in\Cl(F)$ be an ideal class of $F$. There is a zeta function
associated with $B$ which is defined by
\begin{align}
\zeta_B(s)=\sum_{\substack{\frakc\in B \\ \frakc\subset\calO_F}}
\norm(\frakc)^{-s}.
\end{align}
Here $s$ is a complex variable, and the sum runs through all
integral ideals in the ideal class $B$. The series converges for
$\Re(s)>1$. It has a meromorphic continuation to the full complex
plane and the completed zeta function
\begin{align}
\Lambda_B(s)=D^{s/2} \pi^{-s} \Gamma(s/2)^2 \zeta_B(s)
\end{align}
satisfies the functional equation
\begin{align}\label{dedekindfunctional}
\Lambda_B(s)=\Lambda_{\frakd B^{-1}}(1-s)
\end{align} (see e.g.
\cite{Ne}, Chapter VII.5).  Here, $\frakd=\frakd_F=\sqrt{D}\calO_F$
is the different of $F$. The Dedekind zeta function $\zeta_F(s)$ of
$F$ is given by
\begin{align}\label{dede}
\zeta_F(s)=\sum_{B\in \Cl(F)} \zeta_B(s) =
\sum_{\substack{\frakc\subset\calO_F\\\text{integral ideal}}}
\norm(\frakc)^{-s}.
\end{align}

Let $\frakb$ be a fractional ideal in the ideal class $B$. The group
of units $\calO_F^*$ acts on $\frakb\times\frakb$ by
$(c,d)\mapsto(\eps c,\eps d)$ for $\eps\in \calO_F^*$. Recall that
for $\kabcd\in\Sl_2(F)$ and $z=(z_1,z_2)\in\IH^2$, we write
$\norm(cz+d)=(cz_1 +d)(c'z_2 +d')$.

\begin{definition} \label{def:gkb}
Let $k>2$ be an even integer. We define the Eisenstein series of
weight $k$ associated to $B$ by\footnote{The superscript at the summation sign means that the zero summand is omitted.} 
\begin{align*}
G_{k,B}(z)=\norm(\frakb)^k\sideset{}{'}\sum_{(c,d)\in
\calO_F^*\bs\frakb\times\frakb}\norm(cz+d)^{-k}.
\end{align*}
\end{definition}

The Eisenstein series $G_{k,B}$ does not depend on the choice of the
representative $\frakb$ of the ideal class $B$ and converges
uniformly absolutely in every Siegel domain $\calS_t$ ($t>0$).
Consequently, it defines an element of $M_k(\Gamma_F)$. The value at
the cusp $\infty$ is given by
\begin{align*}
G_{k,B}(\infty)&=\lim_{\substack{z\in\calS_t \\ y_1 y_2\rightarrow\infty}} G_{k,B}(z) \\
&= N(\frakb)^k\sideset{}{'}\sum_{d\in\calO_F^*\bs\frakb}\norm(d)^{-k} \\
&= \sideset{}{'}\sum_{d\in \calO_F^*\bs\frakb}\norm(d\frakb^{-1})^{-k} \\
&=\zeta_{B^{-1}}(k).
\end{align*}
%
If $\kappa\in\IP^1(F)$ is any cusp and
$\rho=\kzxz{\alpha}{\beta}{\gamma}{\delta}\in\Sl_2(F)$ with $\rho
\infty=\kappa$, then
\begin{align*}
  G_{k,B}(\kappa)&=\lim_{\substack{z\in\calS_t \\ y_1
      y_2\rightarrow\infty}} (G_{k,B}\mid_k\rho)
  (z)\\
  &=\norm(\fraka)^k \zeta_{[\fraka]B^{-1}}(k),
\end{align*}
where $\fraka=\calO_F\alpha+\calO_F\gamma$.

\begin{theorem} \label{thm:eissplit}
Let $k>2$ be an even integer. The Eisenstein series $G_{k,B}\in
M_k(\Gamma_F)$, where $B\in\Cl(F)$, are linearly independent. The
space $M_k(\Gamma_F)$ can be decomposed as a direct sum
\begin{align*}
M_k(\Gamma_F)=S_k(\Gamma_F)\oplus \bigoplus\limits_{B\in\Cl(F)}\IC
G_{k,B}.
\end{align*}
\end{theorem}

\begin{proof}
See \cite{Ge}, p.~21.
\end{proof}

\begin{remark}
One can define $G_{k,B}$ for $k>2$ odd in the same way. In this case
it is easily seen that $G_{k,B}\equiv 0$ if $\calO_F$ contains a
unit of negative norm. Moreover, the theorem also holds for $k=2$.
In this case one can define $G_{2,B}$ by analytic continuation using
the 'Hecke-trick'. In turns out that all Eisenstein series of weight
$2$ are holomorphic (in contrast to the case of elliptic modular
forms, where the constant term is sometimes non-holomorphic).
\end{remark}

The Fourier expansion of the Eisenstein series can be computed in
the same way as in the case of elliptic modular forms.

\begin{theorem}\label{thm:eisf}
Let $k\geq 2$ even. The Eisenstein series $G_{k,B}$ has the Fourier
expansion
\begin{align*}
G_{k,B}(z)=\zeta_{B^{-1}}(k)+\frac{(2\pi i)^{2k}}{(k-1)!^2} D^{1/2-k} \sum_{\substack{\nu\in\frakd^{-1} \\
\nu\gg 0}} \sigma_{k-1,\frakd B}(\frakd\nu)\,e(\tr(\nu z)).
\end{align*}
Here, for $A\in\Cl(F)$ and an integral ideal $\frakl\subset\calO_F$,
$\sigma_{s,A}(\frakl)$ denotes the divisor sum
\begin{align*}
\sigma_{s,A}(\frakl) =\sum_{\substack{\text{$\frakc\in A$ integral}
\\ \frakc \mid \frakl}} \norm(\frakc)^s.
\end{align*}
\hfill$\square$
\end{theorem}

Using the functional equation of $\zeta_B(s)$, we may write
\begin{align}
\label{normalisedfourier} G_{k,B}(z)=\frac{(2\pi
i)^{2k}}{(k-1)!^2}D^{1/2-k}\bigg[
\frac{1}{4}\zeta_{\frakd B}(1-k) +\sum_{\substack{\nu\in\frakd^{-1} \\
\nu\gg 0}} \sigma_{k-1,\frakd B}(\frakd\nu)\,e(\tr(\nu z))\bigg].
\end{align}

Recall that a Hilbert modular form $f$ of weight $k$ is called
symmetric, if $f(z_1,z_2)=f(z_2,z_1)$. It is easily seen that the
Eisenstein series $G_{k,B}$ are symmetric.

\subsubsection{Restriction to the diagonal} \label{sect:restriction}

If $f\in M_k(\Gamma_F)$ is a Hilbert modular form, we consider its
restriction to the diagonal $g(\tau)=f(\tau,\tau)$. Since the
elliptic modular group $\Sl_2(\IZ)$ is embedded diagonally into
$\Gamma_F=\Sl_2(\calO_F)$, the function $g$ has the transformation
behavior
\begin{align*}
g(\gamma \tau) = f(\gamma\tau,\gamma\tau)=(c\tau +d)^{2k}
f(\tau,\tau)
\end{align*}
for $\gamma=\kabcd\in\Sl_2(\IZ)$. Therefore $g(\tau)$ is an elliptic
modular form for $\Sl_2(\IZ)$ of weight $2k$. If $f$ has the Fourier
expansion
\begin{align*}
f(z)=a_0 + \sum_{\substack{\nu\in\frakd^{-1} \\ \nu\gg 0}} a_{\nu}\,
e(\tr(\nu z))
\end{align*}
at the cusp $\infty$, then $g$ has the expansion
\begin{align}
g(\tau)=a_0 +\sum_{n\geq 1}\sum_{\substack{\nu\in\frakd^{-1} \\
\nu\gg 0 \\ \tr(\nu)=n}} a_\nu \,e(n\tau).
\end{align}
The geometric interpretation is the following. The diagonal
embedding $\H\to \H^2$, $\tau\mapsto (\tau,\tau)$ induces a morphism
$\varphi:\Sl_2(\Z)\bs \H\to Y(\Gamma_F)$, which is birational onto
its image. If we view $f$ as a section of the line bundle of modular
forms of weight $k$ over $Y(\Gamma_F)$, then $g=\varphi^*(f)$ is the
pull-back.

We now consider the restriction of the Eisenstein series $G_{k,B}$.
Using the Fourier expansion \eqref{normalisedfourier}, we see that
\begin{align}\label{eisrestr}
\frac{1}{4}\zeta_{\frakd B}(1-k)+\sum_{n\geq 1}
\sum_{\substack{\nu\in\frakd^{-1} \\ \nu\gg 0 \\
\tr(\nu)=n}}\sigma_{k-1,\frakd B}(\frakd\nu)\,e(n\tau)
\end{align}
is a modular form for $\Sl_2(\IZ)$ of weight $2k$. As a first
consequence we see that the special values $\zeta_{\frakd B}(1-k)$
must be rational numbers. This follows from the fact that the
divisor sums $\sigma_{k-1,\frakd B}(\frakd\nu)$ are rational
(integers), and the fact that the spaces of elliptic modular forms
for $\Sl_2(\Z)$ have bases with rational Fourier coefficients. If
$k=2,4$, then \eqref{eisrestr} must be a multiple of the elliptic
Eisenstein series
\begin{align*}
E_{2k}(\tau)=-\frac{B_{2k}}{4k}+\sum_{n\geq
1}\sigma_{2k-1}(n)\,e(n\tau).
\end{align*}
Here $B_{2k}$ is the usual Bernoulli number and
$\sigma_s(n)=\sum_{d\mid n}d^s$. Comparing the first Fourier
coefficients we obtain a formula for $\zeta_{\frakd B}(1-k)$ due to
Siegel.

\begin{theorem}[Siegel]
If $k=2,4$, then
\begin{align*}
\zeta_B(1-k)=-\frac{B_{2k}}{k}\sum_{\substack{\nu\in\frakd^{-1} \\
\nu\gg 0 \\ \tr(\nu)=1}} \sigma_{k-1,B}(\frakd\nu).
\end{align*}
\hfill$\square$
\end{theorem}

By means of \eqref{dede}, and using $B_4=B_8=-1/30$, we obtain:

\begin{corollary}
The special values of the Dedekind zeta function of $F$ at $-1,-3$
are given by
\begin{align*}
\zeta_F(-1) &=\frac{1}{60}\sum_{\substack{x\in\IZ \\ x^2<D \\ x^2\equiv D\;(4)}}\sigma_1\left(\frac{D-x^2}{4}\right), \\
\zeta_F(-3) &=\frac{1}{120}\sum_{\substack{x\in\IZ \\ x^2<D \\
x^2\equiv D\;(4)}}\sigma_3\left(\frac{D-x^2}{4}\right).
\end{align*}
\hfill $\square$
\end{corollary}

We end this section with a table for these special values.

\begin{table}[h]
\caption{\label{table1} Special values of $\zeta_F(s)$}
\begin{tabular}{|r||c|c|c|c|c|c|c|c|c|c|c|c|c|c| }
\hline \rule[-3mm]{0mm}{8mm}
$D$ & 5 & 8 &   12 & 13 & 17 &  21 & 24 & 28 & 29 & 33& 37 & 40 & 41 & 44\\
\hline \rule[-3mm]{0mm}{8mm}
$\zeta_F(-1)$ & $\frac{1}{30}$ & $\frac{1}{12}$ & $\frac{1}{6}$ & $\frac{1}{6}$ & $\frac{1}{3}$ & $\frac{1}{3}$ &  $\frac{1}{2}$ & $\frac{2}{3}$ & $\frac{1}{2}$& 1&$\frac{5}{6}$ & $\frac{7}{6}$ & $\frac{4}{3}$& $\frac{7}{6}$\\
\hline \rule[-3mm]{0mm}{8mm}
$\zeta_F(-3)$ & $\frac{1}{60}$ & $\frac{11}{120}$& $\frac{23}{60}$& $\frac{29}{60}$& $\frac{41}{30}$& $\frac{77}{30}$&  $\frac{87}{20}$& $\frac{113}{15}$& $\frac{157}{20}$ & $\frac{141}{10}$ & $\frac{1129}{60}$ &$\frac{1577}{60}$& $\frac{448}{15}$&$\frac{2153}{60}$\\
\hline
\end{tabular}
\end{table}

\subsubsection{The example $\Q(\sqrt{5})$}

\label{sect:gundex}

Eisenstein series and restriction to the diagonal can be used to
determine the graded algebra of holomorphic Hilbert modular
forms in some cases where the discriminant of $F$ is small. Here we
illustrate this for $F=\Q(\sqrt{5})$.  The class number of $F$ is
$1$, and the fundamental unit of $\calO_F$ is given by
$\eps_0=\frac{1+\sqrt{5}}{2}$. The graded algebra of modular forms
for the group $\Gamma_F$ was determined by Gundlach \cite{Gu}, see
also \cite{Mu}. For further examples see \cite{Ge}, Chapter~8.

We denote by $g_k$ the Eisenstein series for $\Gamma_F$ of weight
$k$ normalized such that the constant term is $1$ (so
$g_k=\frac{1}{\zeta_{F}(k)} G_{k,\calO_F}$).

\begin{theorem}[Gundlach]
\label{gund1} The graded algebra $M_{2*}^{sym}(\Gamma_F)$ of
holomorphic symmetric Hil\-bert modular forms of even weight for
$\Gamma_F$ is the (weighted) polynomial ring $\C[g_2,g_6,g_{10}]$.
\end{theorem}

Often it is more convenient to replace the generators $g_6$ and
$g_{10}$ by the cusp forms
\begin{align*}
s_6 &= 67\cdot(2^5\cdot 3^3\cdot 5^2)^{-1}\cdot(g_2^3-g_6),\\
s_{10}&= (2^{10}\cdot 3^5\cdot 5^5 \cdot 7)^{-1}\cdot (2^2\cdot
3\cdot 7\cdot 4231 \cdot g_2^5-5\cdot 67\cdot 2293\cdot g_2^2\cdot
g_6 +412751\cdot g_{10}).
\end{align*}
Then Gundlach's result can be restated as
\begin{align}
\label{gund2} M_{2*}^{sym}(\Gamma_F)=\C[g_2,s_6,s_{10}].
\end{align}
Notice that $g_2,s_6,s_{10}$ all have rational integral and coprime
Fourier coefficients.

The key idea for the proof is to show that there is a ``square
root'' for $s_{10}$. Gundlach constructed a cusp form $s_5$ of
weight $5$ for $\Gamma_F$ as a product of $10$ theta constants.
(Later, in Section \ref{sect:borcherdsex}, we will construct it as
the Borcherds lift $\Psi_{1}$.) One can show that $s_5$ is
anti-symmetric, that is, $s_5(z_1,z_2)=-s_5(z_2,z_1)$. Hence it has
to vanish on the diagonal. It turns out that the divisor of $s_5$ is
given by the image of the diagonal in $Y(\Gamma_F)$. Moreover,
$s_5^2=s_{10}$.

\begin{proof}[Proof of Theorem \ref{gund1}]
It is clear that the restriction of $g_2$ to the diagonal is the
Eisenstein series of weight $4$ for $\Sl_2(\Z)$, normalized such
that the constant term is $1$. A quick computation shows that the
restriction of $s_6$ is the delta function, the unique normalized
cusp form of weight 12 for $\Sl_2(\Z)$. Consequently the
restrictions of $g_2$ and  $s_6$ generate the algebra of modular
forms for $\Sl_2(\Z)$ of weight divisible by $4$.

Suppose that $f$ is a  symmetric Hilbert modular form of even
weight $k$ for $\Gamma_F$. Then the restriction to the diagonal of
$f$ is a modular form for $\Sl_2(\Z)$ of weight divisible by $4$ and
therefore a polynomial in the restrictions of $g_2$ and $s_6$. Hence
there is a polynomial $P\in \C[X,Y]$ such that
\[
f_1=f-P(g_2,s_6)
\]
vanishes on the diagonal. Therefore $f_1/s_5$ is a holomorphic
Hilbert modular form for $\Gamma_F$. It is anti-symmetric and
therefore vanishes on the diagonal. Consequently, $f_1/s_5^2\in
M_{k-10}(\Gamma_F)$ is symmetric. Now the assertion follows by
induction on the weight.
\end{proof}

To get the full algebra $M_{*}(\Gamma_F)$ of Hilbert modular forms
for $\Gamma_F$ one needs in addition the existence of a symmetric
Hilbert cusp form $s_{15}$ of weight $15$. Gundlach constructed it
as a product of differences of Eisenstein series of weight $1$ for a
principal congruence subgroup of $\Gamma_F$. (In Section
\ref{sect:borcherdsex} we will construct it as the Borcherds lift
$\Psi_{5}$.)

\begin{theorem}[Gundlach]
\label{gund3} The graded algebra $M_{*}(\Gamma_F)$ of Hilbert
modular forms for $\Gamma_F$ is generated by $g_2, s_5, s_{6},
s_{15}$. The anti-symmetric cusp form $s_5$ and the symmetric cusp
form $s_{15}$ satisfy the following relations over
$M_{2*}^{sym}(\Gamma_F)=\C[g_2,s_6,s_{10}]$:
\begin{align*}
s_5^2&=s_{10},\\
s_{15}^2&=5^5 \cdot s_{10}^3-2^{-1}\cdot 5^3 \cdot g_2^2 s_6
s_{10}^2 + 2^{-4} \cdot g_2^5 s_{10}^2
+2^{-1}\cdot 3^2\cdot 5^2 \cdot g_2 s_6^3 s_{10}\\
&\phantom{=}{}-2^{-3}\cdot  g_2^4 s_6^2 s_{10} -2\cdot 3^3 \cdot
s_6^5 +2^{-4} \cdot g_2^3 s_6^4.
\end{align*}
\hfill $\square$
\end{theorem}

\subsection{The $L$-function of a Hilbert modular form}

\label{sect:L}

In this section we briefly discuss how one can attach an
$L$-function to a Hilbert modular form. First, one needs to know
that the coefficients have  polynomial growth.

\begin{proposition}[Hecke estimate]
Let $f=\sum_{\nu}a_{\nu}\,e(\tr(\nu z))\in M_k(\Gamma)$.
\begin{enumerate}
\item[(i)] Then $a_{\nu}=O(\norm(\nu)^k)$ for $\norm(\nu)\rightarrow\infty$.
\item[(ii)] If $f$ is   a cusp form, we have the stronger estimate
$a_{\nu}=O(\norm(\nu)^{k/2})$ for $\norm(\nu)\rightarrow\infty$.
\end{enumerate}
\end{proposition}

\begin{proof}
Here we only carry out the proof for cusp forms. For non-cuspidal
modular forms one has to slightly modify the argument. (For the
group $\Gamma_F$ one can also use Theorems~\ref{thm:eissplit} and
\ref{thm:eisf}.) According to \eqref{eq:coeff} we have
\begin{align*}
a_{\nu}=\frac{1}{\vol(\IR^2/M)} \int\limits_{\IR^2/M} f(z)
e(-\tr(\nu z))\,dx_1\,dx_2.
\end{align*}
By Proposition \ref{hbound} we know that $|f(z) (y_1 y_2)^{k/2}|$ is
bounded on $\H^2$. Hence there is a constant $C>0 $ such that
\[
|a_\nu|\leq C \int\limits_{\IR^2/M} (y_1 y_2)^{-k/2} e^{-2\pi
\tr(\nu y)}\,dx_1\,dx_2
\]
for all $y\in (\R_{>0})^2$. Choosing
$y=(\frac{1}{\nu},\frac{1}{\nu'})$, we see that
\[
|a_\nu|\leq C \vol(\IR^2/M)\norm(\nu)^{k/2}.
\]
This proves the proposition.
\end{proof}

For the rest of this section we only consider the full Hilbert
modular group $\Gamma_F$. Let $f\in M_k(\Gamma_F)$,
and denote the Fourier expansion by
\begin{align*}
f=a_0 +\sum_{\substack{\nu\in\frakd^{-1} \\ \nu\gg 0}}a_{\nu}\,
e(\tr(\nu z)).
\end{align*}
Let $\calO_F^{*,+}$ be the group of totally positive units of
$\calO_F$, and let $U=\{\eps^2;\; \eps\in \calO_F^{*,+}\}$. Then $U$
has index $2$ in the cyclic group $\calO_F^{*,+}$. We have
$a_{\nu}=a_{\eps\nu}$ for all $\nu\in \frakd^{-1}$ and all $\eps\in
U$.
%

\begin{definition}
We define an $L$-series associated to $f$ by
\begin{align*}
L(f,s)&=\sum_{\substack{\nu\in\frakd^{-1}/ U\\ \nu\gg 0}}
a_{\nu}\norm(\nu\frakd)^{-s}.
\end{align*}
\end{definition}

\begin{example}
For the Eisenstein series $G_{k,B}$ (see Definition \ref{def:gkb})
one easily checks that
\begin{align*}
L( G_{k,B},s)= 2 \frac{(2\pi
i)^{2k}}{(k-1)!^2}D^{1/2-k}\zeta_{\frakd^{-1}B^{-1}}(s)\,\zeta_{\frakd
B}(s+1-k).
\end{align*}
\end{example}

The functional equation \eqref{dedekindfunctional} of the partial
Dedekind zeta function $\zeta_{B}(s)$ implies that $ L( G_{k,B},s)$
has a meromorphic continuation and satisfies a functional equation
relating $s$ and $k-s$. Therefore it is reasonable to expect similar
properties for the $L$-functions of cusp forms as well.

\begin{theorem}\label{thm:L}
Let $f\in S_k(\Gamma_F)$. The completed $L$-function
\begin{align*}
\Lambda(f,s)=D^s(2\pi)^{-2s}\Gamma(s)^2 L(f,s)
\end{align*}
has a holomorphic continuation to $\IC$, is entire and bounded in
vertical strips, and satisfies the functional equation
\begin{align*}
\Lambda(f,s)=(-1)^{k}\Lambda(f,k-s).
\end{align*}
\end{theorem}

\begin{proof}
Using the Euler integral for the Gamma function, we see that
\[
(2\pi)^{-2s}\Gamma(s)^2 \norm(\nu)^{-s}
=\int\limits_0^\infty\int\limits_0^\infty e^{-2\pi \tr(\nu y)} (y_1
y_2)^{s} \,\frac{dy_1}{y_1}\,\frac{dy_2}{y_2} .
\]
Hence, by unfolding we find that
\[
\Lambda(f,s)=\int\limits_{(\IR_{>0})^2/ U} f(iy) (y_1  y_2)^s\,
\frac{dy_1}{y_1}\,\frac{dy_2}{y_2}.
\]
We split up the integral into an integral over $y_1 y_2>1$ and a
second integral over $y_1 y_2<1$. The modularity of $f$ implies that
$f(i(\frac{1}{y_1},\frac{1}{y_2}))= (-1)^k (y_1 y_2)^k f(iy)$. Hence
the second integral can be transformed into an integral over $y_1
y_2>1 $ as well. We find that
\[
\Lambda(f,s)=\int\limits_{\substack{(\IR_{>0})^2/ U \\ y_1 y_2>1}}
f(iy) (y_1 y_2)^s\, \frac{dy_1}{y_1}\,\frac{dy_2}{y_2} +(-1)^k
\int\limits_{\substack{(\IR_{>0})^2/ U\\ y_1 y_2 >1}} f(iy) (y_1
y_2)^{k-s}\, \frac{dy_1}{y_1}\,\frac{dy_2}{y_2}.
\]
This integral representation converges for all $s\in \C$ and defines
the holomorphic continuation of $\Lambda(f,s)$. Moreover, the
functional equation is obvious now.
\end{proof}

We now suppose that $\calO_F$ contains a unit of negative norm. Then
$M_k(\Gamma_F)=\{0\}$ for $k$ odd. So we further suppose that $k$ is
even. Then the Fourier coefficients of $f\in M_k(\Gamma_F)$ satisfy
$a_{\nu}=a_{\eps\nu}$ for all $\nu\in \frakd^{-1}$ and all $\eps\in
\calO_F^{*,+}$. Thus, $a_{\nu}$ only depends on the ideal
$(\nu\frakd)\subset\calO_F$, and we write $a((\nu\frakd))=a_{\nu}$.
Then we may rewrite the $L$-function of $f$ in the form
\[
L(f,s)=\sum_{\substack{\fraka\subset\calO_ F\\ \text{principal
ideal}}} a(\fraka)\norm(\fraka)^{-s}.
\]
So this $L$-function is analogous to the zeta function $\zeta_B(s)$
associated to an  ideal class $B$ of $F$ (here the unit class). It
is natural to associate more general $L$-functions to $f$, for
instance, an $L$-function where one sums over all integral ideals of
$F$ analogous to the full Dedekind zeta function of $F$. To this end
it is more convenient to view Hilbert modular forms as automorphic
functions on $(\operatorname{Res}_{F/\Q}\Sl_2)(\A)$, where
$\operatorname{Res}_{F/\Q}$ denotes the Weil restriction of scalars,
and $\A$ denotes the ring of adeles of $\Q$ (cf. \cite{Ga}).


\section{The orthogonal group $\Orth(2,n)$}

An important property of Hilbert modular surfaces is that they can
also be regarded as modular varieties associated to the orthogonal
group of a suitable rational quadratic space $V$ of type $(2,2)$.
There is an accidental isomorphism
$\operatorname{Res}_{F/\Q}\Sl_2\cong \Spin_V$ of algebraic groups
over $\Q$. Modular varieties for orthogonal groups $\Orth(2,n)$ come
with natural families of special algebraic cycles on them arising
from embeddings of ``smaller'' orthogonal groups. They provide a
rich source of extra structure and can be used to study geometric
questions. In the $\Orth(2,2)$-case of Hilbert modular surfaces
these special cycles lead to Hirzebruch-Zagier divisors (codimension
$1$) and certain CM-points (codimension $2$).

To put things in the right context, in this section we study
quadratic spaces and modular forms for orthogonal groups in slightly
greater generality than needed for the application to Hilbert
modular surfaces. However, we hope that this will rather clarify
things than complicate them. For a detailed  account of the theory
of quadratic forms and orthogonal groups we refer to \cite{Ki},
\cite{Sch}.

\subsection{Quadratic forms}

Let $R$ be a commutative ring with unity $1$. We write $R^*$ for the
group of units in $R$. Let $M$ be a finitely generated
$R$-module.  A {\em quadratic form} on $M$ is a mapping $Q:M\to R$
such that
\begin{enumerate}
\item[(i)]
$Q(rx)=r^2Q(x)$ for all $r\in R$ and all $x\in M$,
\item[(ii)]
$B(x,y):=Q(x+y)-Q(x)-Q(y)$ is a bilinear form.
\end{enumerate}
The first condition follows from the second if $2$ is invertible in
$R$.  Then we have $Q(x)=\frac{1}{2} B(x,x)$.  The pair $(M,Q)$ is
called a quadratic module over $R$. If $R$ is a field, we frequently
say space instead of module. Two elements $x,y\in M$ are called
  orthogonal if $B(x,y)=0$. If $A\subset M$ is a subset, we denote
the orthogonal complement by
\begin{align}
A^\perp=\{ x\in M;\; \text{$B(x,y)=0$ for all $y\in A$}\}.
\end{align}
It is a submodule of $M$. For $x\in M$ we briefly write $x^\perp$
instead of $\{x\}^\perp$.  The quadratic module $M$ is called {\em
  non-degenerate}, if $M^\perp=\{0\}$.
%
A non-zero vector $x\in M$ is called {\em isotropic} if $Q(x)=0$,
and {\em anisotropic}, if $Q(x)\neq 0$, respectively.

Let $(M,Q)$ and $(M',Q')$ be quadratic modules over $R$. An
$R$-linear map $\sigma:M\to M'$ is called an {\em isometry}, if
$\sigma$ is injective and
\[
Q'(\sigma(x)) = Q(x)
\]
for all $x\in M$.  If $\sigma$ is also surjective then $M$ and $N$
are called isometric.  The {\em orthogonal group} of $M$,
\begin{align}
\Orth_M=\{ \sigma\in \Aut(M);\quad \text{$\sigma$ isometry}\},
\end{align}
is the group of all isometries from $M$ onto itself. The special
orthogonal group is the subgroup
\begin{align}
\SO_M=\{ \sigma\in \Orth_M;\quad \det(\sigma)=1\}.
\end{align}

Important examples of isometries are given by {\em reflections}. For
an element $x\in M$ with $Q(x)\in R^*$ we define $\tau_x:M\to M$ by
\begin{align}
\tau_x(y)=y-B(y,x)Q(x)^{-1} x, \qquad y\in M.
\end{align}
Then $\tau_x$ is an isometry and satisfies
\begin{align*}
\tau_x(x)&=-x\\
\tau_x(y)&=y, \quad \text{for $y\in x^\perp$,}\\
\tau_x^2&=\id .
\end{align*}
So $\tau_x$ is the reflection in the hyperplane $x^\perp$.

Further examples of isometries are given by {\em Eichler elements}.
Let $u\in M$ be isotropic and $v\in M$ with $B(u,v)=0$. We define
$E_{u,v}:M\to M$ by
\begin{align}
E_{u,v}(y)=y+B(y,u)v-B(y,v)u-B(y,u)Q(v)u, \qquad y\in M.
\end{align}
One easily checks that $E_{u,v}$ is an isometry and
\begin{align*}
E_{u,v}(u)&=u,\\
E_{u,v}(v)&=v-2Q(v)u, \\
E_{u,v_1}E_{u,v_2}&=E_{u,v_1+v_2} \quad \text{for $v_1,v_2\in
u^\perp$.}
\end{align*}

If $M$ is free, and $v_1,\dots,v_n$ is a basis of $M$, we have the
corresponding {\em Gram matrix} $S=(B(v_i,v_j))_{i,j}$. The class of
$\det (S)$ in $R^*/(R^*)^2$ is independent of the choice of the
basis. It is called the discriminant of $M$ and is denoted by
$d(M)$.
Note that if $v_1,\dots, v_n$ is an orthogonal basis of $M$, we have
\begin{align}
d(M)=2^nQ(v_1)\cdots Q(v_n).
\end{align}

For the rest of this subsection, let $M$ be a quadratic space of
dimension $n$ over a field $k$ of characteristic $\neq 2$.
The space is non-degenerate, if and only if $d(M)\neq 0$. If $M$ is
non-degenerate and $v_1\in M$ is an anisotropic vector,
there exist anisotropic vectors $v_2,\dots, v_n\in M$ such that
$v_1,\dots, v_n$ is an orthogonal basis of $M$. Consequently, the
reflection corresponding to $v_1$ satisfies $\det(\tau_{v_1})=-1$.

\begin{theorem}\label{thm:refl}
  Let $M$ be a regular quadratic space over a field $k$ of
  characteristic $\neq 2$. Then $\Orth_M$ is generated by
  reflections. Moreover,
$\SO_M$ is the subgroup of elements of $\Orth_M$ which can be
written as a product of an even number of reflections. \hfill
$\square$
\end{theorem}

\begin{example}\label{quadreal}
Let $p,q$ be non-negative integers. We denote by $\R^{p,q}$ the
quadratic space over $\R$ given by $\R^{p+q}$ with the quadratic
form
\[
Q(x)=x_1^2+\dots +x_{p}^2 - x_{p+1}^2\dots -x_{p+q}^2.
\]
If $V$ is a finite dimensional quadratic space over $\R$, then there
exist non-negative integers $p,q$ such that $V$ is isometric to
$\R^{p,q}$. The pair $(p,q)$ is uniquely determined by $V$ and is
called the {\em type} of $V$. Moreover, $p-q$ is called the {\em
  signature} of $V$. The orthogonal group of $\R^{p,q}$ is also
denoted by $\Orth(p,q)$.
\end{example}

\subsection{The Clifford algebra}

As before, let $R$ be a commutative ring with unity $1$, and let
$(V,Q)$ be a finitely generated quadratic module over $R$. If $A$ is
any $R$-algebra, we write $Z(A)$ for the center of $A$.

We consider the tensor algebra
\[
T_V= \bigoplus_{m=0}^\infty V^{\otimes m} = R\oplus V\oplus
(V\otimes_R V) \oplus \dots
\]
of $V$. Let $I_V\subset T_V$ be the two-sided ideal
generated by $v\otimes v-Q(v)$ for $v\in V$. The {\em Clifford
algebra} of $V$ is defined by
\begin{align}
C_V=T_V/I_V.
\end{align}
Observe that $R$ and $V$ are embedded into $C_V$ via the canonical
maps.  For simplicity, the element of $C_V$ represented by
$v_1\otimes \dots \otimes v_m$ (where $v_i\in V$) is denoted by $v_1
\cdots v_m$. By definition, we have for $u,v\in V\subset C_V$:
\begin{align*}
v^2&=Q(v),\\
uv+vu&=B(u,v).
\end{align*}
In particular, $uv=-vu$ if and only if $u$ and $v$ are orthogonal. The
Clifford algebra has the following universal property.

\begin{proposition}\label{univ}
Let $f:V\to A$ be an $R$-linear map to an $R$-algebra $A$ with unity
$1_A$ such that $f(v)^2=Q(v)1_A$ for all $v\in V$. Then there exists
a unique $R$-algebra homomorphism $C_V\to A$ such that the following
diagram commutes:
\[
\xymatrix{ V \ar[dr]\ar[r]& C_V \ar[d]\\
&A}.
 \]
\end{proposition}

The universal property implies that an isometry $\varphi:V\to V'$ of
quadratic spaces over $R$ induces a unique $R$-algebra homomorphism
$\tilde \varphi: C_V\to C_{V'}$ compatible with the natural
inclusions. Therefore, the assignment $V\mapsto C_V$ defines a
functor from the category of quadratic spaces over $R$ with
isometries as morphisms to the category of associative $R$-algebras
with unity.

Moreover, if we fix a quadratic space $(V,Q)$ over $R$, for any
commutative $R$-algebra $S$ with unity we can consider the extension
of scalars $V(S)=V\otimes_R S$, which is a quadratic module over $S$
in a natural way (the quadratic form being defined by $Q(v\otimes
s)=s^2Q(v)$). In the same way, we consider $C_V(S)=C_V\otimes_R S$.
One easily checks that $C_V(S)=C_{V(S)}$. So taking the Clifford
algebra commutes with extension of scalars.

\begin{example}
We denote by $C_{p,q}$ the Clifford algebra of the real quadratic
space $\R^{p,q}$ of Example \ref{quadreal}. For small $p,q$ we have
\[
C_{0,0}=\R,\quad C_{1,0}=\R\oplus\R,\quad C_{0,1}=\C, \quad
C_{2,0}=\Mat_2(\R), \quad C_{1,1}=\Mat_2(\R),\quad C_{0,2}=\H.
\]
Here $\H$ denotes the Hamilton quaternion algebra.
(This follows from Examples \ref{ex:n1} and \ref{ex:n2} below). 
\end{example}

Now assume that $V$ is free. If $v_1,\dots,v_n$ is a basis of $V$,
then these vectors generate $C_V$ as an $R$-algebra. The elements
\[
v_{i_1}\cdots v_{i_r}\qquad (\text{$1\leq i_1<\dots<i_r\leq n$ and
$0\leq r\leq n$})
\]
form a basis of $C_V$. In particular, $C_V$ is a free $R$-module of
rank $2^n$. Observe that for the trivial quadratic form $Q\equiv 0$
the Clifford algebra $C_V$ is simply the Grassmann algebra of $V$.

We write $C^0_V$ for the $R$-subalgebra of $C_V$ generated by
products of an even number of basis vectors of $V$, and  $C^1_V$ for
the $R$-submodule of $C_V$ generated by products of an odd number of
basis vectors of $V$. This definition is meaningful, since the
defining relations of $C_V$ only involve products of an even number
of basis vectors. We obtain a decomposition
\[
C_V=C^0_V\oplus C^1_V,
\]
which is a $\Z/2\Z$-grading on $C_V$. The subalgebra $C_V^0$ is
called the {\em even Clifford algebra} of $V$ (or the second
Clifford algebra of $V$).

Multiplication by $-1$ defines an isometry of $V$. By Proposition
\ref{univ} it induces an algebra automorphism
\begin{align}
J:C_V\longrightarrow C_V,
\end{align}
called the {\em canonical automorphism}. If $2$ is invertible in
$R$, then the even Clifford algebra can be characterized by
\[
C_V^0=\{ v\in C_V;\; J(v)=v\}.
\]

There is a second involution on $C_V$, which is an
anti-automorphism. It is called the {\em canonical involution} on
$C_V$ and is defined by
\begin{align}
{}^t:C_V\longrightarrow C_V,\quad  (x_1\otimes\dots \otimes  x_m)^t=
x_m\otimes\dots \otimes  x_1 .
\end{align}
It reduces to the identity on $R\oplus V$. It is used to define the
{\em Clifford norm} on $C_V$ by
\begin{align}
\norm: C_V \longrightarrow C_V,\quad \norm(x)=x^t x.
\end{align}
For $x\in V$ we have $\norm(x)=Q(x)$. So the norm map extends the
quadratic form on $V$. Note that the Clifford norm is in general
{\em not} multiplicative on $C_V$.

For the rest of this subsection, let $k$ be a field of
characteristic $\ne 2$, and let $(V,Q)$ be a non-degenerate
quadratic space over $k$ of dimension $n$. Moreover, let $v_1,\dots,
v_n$ be an orthogonal basis of $V$. We put
\[
\delta=v_1\cdots v_n\in C_V.
\]

\begin{remark}
\label{rem:disc} When $n$ is even, we have
\[
\delta v_i = -v_i \delta,\quad \delta^2 = (-1)^{n/2}2^{-n}d(V)\in
k^*/(k^*)^2.
\]
When $n$ is odd, we have
\[
\delta v_i = v_i \delta,\quad \delta^2 = (-1)^{(n-1)/2}2^{-n}d(V)\in
k^*/(k^*)^2.
\]
\hfill $\square$
\end{remark}


\begin{theorem}
\label{thm:center} The center of $C_V$ is given by
\[
Z(C_V)=\begin{cases}k &\text{if $n$ is even,}\\
k+k\delta&\text{if $n$ is odd.}
\end{cases}
\]
The center of $C_V^0$ is given by
\[
Z(C_V^0)=\begin{cases}k+k\delta&\text{if $n$ is even,}\\
k&\text{if $n$ is odd.}
\end{cases}
\]
\hfill $\square$
\end{theorem}

Let $A$ be a ring with unity such that $k\subset Z(A)$.
Recall that $A$ is called a {\em quaternion algebra} over $k$, if it has a basis $\{1,x_1,x_2,x_3\}$  as a $k$-vector space such that
\[
x_1^2=\alpha,\quad x_2^2=\beta,\quad x_3=x_1x_2 = -x_2 x_1
\]
for some $\alpha,\beta\in k^*$. Then it is denoted by $(\alpha,\beta)$.
The parameters $\alpha,\beta$ determine $A$ up to $k$-algebra isomorphism.
It is easily seen that $k=Z(A)$.
The conjugation in $A$ is defined by
\[
x=a_0+a_1 x_1 + a_2 x_2 + a_3 x_3 \mapsto \bar x = a_0-a_1 x_1 - a_2 x_2 - a_3 x_3
\]
for $x\in A$. The norm is defined by $\norm(x)=x\bar{x}\in k$.  A
quaternion algebra over $k$ is either isomorphic to $\Mat_2(k)$ or it
is a division algebra.  For more details we refer to \cite{Ki}
Chapter~1.5.  We end this section by giving some examples of Clifford
algebras associated to low dimensional quadratic spaces (see
\cite{Ki}, p.~28).

\begin{example}\label{ex:n1}
If $n=1$ then $C_V=k+k\delta$ and $\delta^2=d(V)/2$. As a $k$-algebra,
we have
\[
C_V\cong k[X]/(X^2-d(V)/2).
\]
When $d(V)/2$ is not a square in $k$, this is a quadratic field
extension of $k$. When $d(V)/2$ is a square in $k$, then $C_V\cong
k\oplus k$.
\end{example}

\begin{example}\label{ex:n2}
Suppose that $n=2$ and that $V$ has the orthogonal basis $v_1,v_2$
with $Q(v_i)=q_i\in k^*$. Then $C_V=k\oplus kv_1\oplus kv_2\oplus
kv_1 v_2$ is isomorphic to the quaternion algebra $(q_1,q_2)$ over
$k$. Moreover, $C_V^0\cong k[X]/(X^2 + d(V))$.
\end{example}

\begin{example}\label{ex:n3}
Suppose that $n=3$ and that $V$ has the orthogonal basis $v_1,v_2, v_3$
with $Q(v_i)=q_i\in k^*$. Then $C^0_V=k\oplus kv_1v_2\oplus kv_2v_3\oplus
kv_1 v_3$ is isomorphic to the quaternion algebra $(-q_1q_2,-q_2q_3)$ over
$k$. The
conjugation in the quaternion algebra is identified with the main
involution of the Clifford algebra, and the norm with the Clifford
norm.
\end{example}

\begin{example} \label{ex:n4}
Suppose that $n=4$ and that $V$ has the orthogonal basis
$v_1,v_2,v_3,v_4$ with $Q(v_i)=q_i\in k^*$.  Then the center $Z$ of
the even Clifford algebra $C_V^0$ is $k+k\delta$, and we have
\[
C^0_V=Z+Zv_1 v_2 + Zv_2 v_3 + Z v_1 v_3.
\]
Since $(v_1v_2)^2=-q_1 q_2$, $(v_2v_3)^2=-q_2 q_3$, and
$(v_1v_2)(v_2v_3)=q_2 (v_1 v_3)$, the algebra $C^0_V$ is isomorphic
to the quaternion algebra $(-q_1 q_2,-q_2 q_3)$ over $Z$.  The
conjugation in the quaternion algebra is identified with the main
involution of the Clifford algebra, and the norm with the Clifford
norm.
\end{example}

\subsection{The Spin group}
As before, let $R$ be a commutative ring with unity $1$, and let
$(V,Q)$ be a finitely generated quadratic module over $R$. The {\em
Clifford group} $\CG_V$ of $V$ is defined by
\begin{align}
\CG_V=\{ x\in C_V;\quad \text{$x$ invertible and $x V J(x)^{-1} =
V$}\}.
\end{align}
By definition, every $x\in \CG_V$ defines an automorphism $\alpha_x$
of $V$ by
\[
\alpha_x(v)=x v J(x)^{-1}. \] We obtain a linear representation
$\alpha :\CG_V\to \Aut_R(V)$, $x\mapsto \alpha_x$, called the {\em
vector representation}.

It is easily seen that the involution $x\mapsto x^t$ takes $\CG_V$
to itself. Consequently, if $x\in \CG_V$, then the Clifford norm
$\norm(x)$ belongs to $\CG_V$ as well.

\begin{lemma}\label{lem:normtarget}
The Clifford norm induces a group homomorphism $\CG\to Z(C_V)^*\cap
C^0_V$.
\end{lemma}

\begin{proof}
Let $x\in \CG_V$. We first show that $\norm(x)$ acts trivially on
$V$ via $\alpha$. Let $v\in V$. Since $w:=\alpha_x(v)\in V$, we have
$J(w)^t=-w$. This implies
\[
(x^t)^{-1} v J(x^t) = x v J(x)^{-1}
\]
and therefore
\[
\norm(x) v J(\norm(x))^{-1} =\norm(x) v \norm(x)^{-1}=  v.
\]
Since $V$ generates $C_V$ as an algebra, we see that $\norm(x)\in
Z(C_V)$. Moreover, it is clear that $\norm(x)$ is invertible and
contained in the even part of the Clifford algebra.

Now we see that for $x,y\in \CG_V$ we have
\[
\norm(xy)=(xy)^t(xy)=y^t(x^tx)y=(x^tx)(y^ty).
\]
This concludes the proof of the lemma.
\end{proof}

\begin{lemma}
For $x\in \CG_V$ the automorphism $\alpha_x\in \Aut_R(V)$ is an
isometry.
\end{lemma}

\begin{proof}
Let $v\in V$. Since $w=\alpha_x(v)\in V$, we have
\begin{align*}
Q(w)&=\norm(w)\\
&=J(x^{-1})^t v^t x^t x v J(x^{-1})\\
&=Q(v).
\end{align*}
This shows that $\alpha_x$ is an isometry.
\end{proof}

Consequently, the vector representation defines a homomorphism
\begin{align}\label{vectorrep}
\alpha:\CG_V\to \Orth_V.
\end{align}
Moreover, if $x\in \CG_V\cap V$, then $Q(x)\in R^*$ and $ \alpha_x$
is equal to the reflection $\tau_x$ in the hyperplane $x^\perp$.

\begin{definition}
We define the {\em general Spin group} $\GSpin_V$  and the {\em Spin
group} $\Spin_V$ of $V$ by
\begin{align*}
\GSpin_V&=\CG_V\cap C_V^0,\\
 \Spin_V&=\{ x\in \GSpin_V;\; \norm(x)=1\}.
\end{align*}
\end{definition}

For the rest of this section we assume that $R=k$ is a field of
characteristic $\neq 2$. We briefly discuss the structure of the
Clifford and the Spin group.

In this case, by Theorem \ref{thm:refl}, the vector representation
\eqref{vectorrep} is surjective onto $\Orth_V$. Moreover, the kernel is given by
$k^*$ (see \cite{Sch}, Chapter 9.3). Hence $\CG_V$ and $\GSpin_V$
are central extensions of $\Orth_V$ and $\SO_V$, respectively,
\begin{gather*}
\xymatrix{ 1 \ar[r] & k^*\ar[r] & \CG_V \ar[r] & \Orth_V\ar[r] &1},\\
\xymatrix{ 1 \ar[r] & k^*\ar[r] & \GSpin_V \ar[r] & \SO_V\ar[r] &1}.
\end{gather*}

According to Lemma \ref{lem:normtarget} and Theorem
\ref{thm:center}, the Clifford norm defines a homomorphism $\CG_V\to
k^*$. It induces a homomorphism
\begin{align}
\theta:\Orth_V\longrightarrow  k^*/(k^*)^2,
\end{align}
called the {\em spinor norm}. It is characterized by the property
that for the reflection $\tau_v$ corresponding to an anisotropic
vector $v\in V$ we have
\[
\theta(\tau_v)= Q(v)\in k^*/(k^*)^2.
\]
We obtain the exact sequence
\begin{gather*}
\xymatrix{ 1 \ar[r] & \{\pm 1\}\ar[r] & \Spin_V \ar[r]^\alpha &
\SO_V\ar[r]^\theta & k^*/(k^*)^2}.
\end{gather*}

The groups $\CG_V$, $\GSpin$, and $\Spin_V$ can be viewed as  the
groups of $k$-valued points of an affine algebraic group over $k$.
If $A$ is a commutative $k$-algebra with unity, then the group of
$A$-valued points of $\CG_V$ is $\CG_V(A)=\CG_{V(A)}$, and
analogously for the other groups.

The following lemma will be useful in Section
\ref{sect:hilbert-orth}.

\begin{lemma}
\label{lem:nleq4} Assume that $\dim(V)\leq 4$. Then
\begin{align*}
\GSpin_V &= \{ x\in C_V^0;\quad \norm(x)\in k^*\},\\
\Spin_V &= \{ x\in C_V^0;\quad \norm(x)=1\}.
\end{align*}
\end{lemma}

\begin{proof}
It is clear that the left hand sides are contained in the right hand
sides.

Conversely, let $x\in C_V^0$ with $\norm(x)\in k^*$. Then $x$ is
invertible, because $y=x^t \norm(x)^{-1}\in  C_V^0$ is inverse to
$x$. Hence, it suffices to show that $x V x^{-1} \subset V$.

Let $v\in V$. It is clear that $w:=x v x^{-1} \in C^1_V$. The
assumption $\dim(V)\leq 4$ implies that
\[
V=\{ g\in C^1_V;\quad g^t=g\}.
\]
Therefore it suffices to show that $w=w^t$. Since $\norm(x)\in k^*$,
we have $\norm(x)v\norm(x)^{-1}=v$. This implies that
\[
x v x^{-1} = (x^t)^{-1}v x^t
\]
and therefore $w=w^t$.
\end{proof}

\subsubsection{Quadratic spaces in dimension four}

We now consider the special cases that $(V,Q)$ is a rational
quadratic space of dimension $4$ over the field $k$.  Let
$v_1,v_2,v_3,v_4$ be an orthogonal basis of $V$ and put
$q_i=Q(v_i)\in k^*$.  By means of Example \ref{ex:n4} and Lemma
\ref{lem:nleq4}, we see that $\Spin_V$ is the group of norm $1$
elements in the quaternion algebra $(-q_1 q_2,-q_2 q_3)$ over
$Z:=Z(C^0_V)=k+k\delta$, where $\delta=v_1v_2v_3v_4$.

We would like to describe the vector representation of $\Spin_V$
intrinsically in terms of $C_V^0$. This can be done by identifying
$V$ with an isometric copy $\tilde V$ in $C_V^0$. (Note that by
definition $V\not\subset C_V^0$.)  The vector representation on $V$
translates to a ``twisted'' vector representation on $\tilde V$. We
partly follow \cite{KR} \S0.

\begin{lemma}
\label{lem:conj} Let $v_0\in V$ with $q_0=Q(v_0)\neq 0$, and denote
by $\sigma=\Ad(v_0)$  the adjoint automorphism of $C_V^0$ associated
to $v_0$, i.e., $x^\sigma= v_0 x v_0^{-1}$ for $x\in C_V^0$. Then
$\delta^\sigma=-\delta$ and the fixed algebra of $\sigma$ in $C_V^0$
is a quaternion algebra $B_0$ over $k$ such that $C_V^0=B_0\otimes_k
Z$.
\end{lemma}

\begin{proof}
See \cite{KR} Lemma 0.2.
\end{proof}

In particular, on the center $Z$ of $C^0_V$, the automorphism
$\sigma$ agrees with the conjugation in $Z/k$. Let
\begin{align}
\label{def:tildev} \tilde V&= \{ x\in C_V^0;\quad x^t=x^\sigma\}.
\end{align}
This is a quadratic space over $k$ together with the quadratic form
\begin{align}
\tilde Q(x)=q_0\cdot x^\sigma x = q_0 \cdot \norm(x).
\end{align}
The group $\Spin_V$ acts on $\tilde V$ by
\begin{align}
x\mapsto \tilde\alpha_g(x):=g x g^{-\sigma},
\end{align}
for $x\in \tilde V$ and $g\in \Spin_V$. The quadratic form $\tilde
Q$ is preserved under this action:
\begin{align}
\tilde Q(g x g^{-\sigma})=q_0\cdot (g x g^{-\sigma})^t (g x
g^{-\sigma}) =q_0\cdot x^t x = \tilde Q(x).
\end{align}

\begin{lemma}
\label{lem:viso} The assignment $x\mapsto x\cdot v_0$ defines an
isometry of quadratic spaces
\[
(\tilde V, \tilde Q)\longrightarrow (V,Q),
\]
which is compatible with the actions of $\Spin_V$.
\end{lemma}

\begin{proof}
See \cite{KR} Lemma 0.3.
\end{proof}

\subsection{Rational quadratic spaces of type $(2,n)$.}

\label{sect:orth}

Let $(V,Q)$ be a non-degenerate quadratic space over $\Q$. Then
$V(\R)=V\otimes_\Q\R$ is isometric to $\R^{p,q}$ for a unique pair
of non-negative integers $(p,q)$, called the type of $V$. If
$K\subset \Orth_V(\R)$ is a maximal compact subgroup, then $
\Orth_V(\R)/K$ is a symmetric space. It is hermitian, i.e., has a
complex structure, if and only if $p=2$ or $q=2$. Since this is the
case of interest to us, throughout this subsection we assume that
$V$ has type $(2,n)$. We discuss several realizations of the
corresponding hermitian symmetric domain. We frequently write
$(\cdot, \cdot)$ for the bilinear form $B( \cdot, \cdot)$.

\subsubsection{The Grassmannian model}

We consider the Grassmannian of $2$-dimensional subspaces of $V(\R)$
on which the quadratic form is positive definite,
\[
\Gr(V)=\{ v\subset V(\IR);\;\quad \text{$\dim v=2$ and $Q|_v >0$}
\}.
\]
By Witt's theorem, $\Orth_V(\R)$ acts transitively on $\Gr(V)$. If
$v_0\in \Gr(V)$ is fixed, then the stabilizer $K$ of $v_0$ is a
maximal compact subgroup of  $\Orth_V(\R)$, and $K \cong
\Orth(2)\times \Orth(n)$. Thus $\Gr(V)\cong \Orth_V(\R)/K$ is a
realization of the hermitian symmetric space. The Grassmannian model
has the advantage that it provides an easy description of
$\Orth_V(\R)/K$, but unfortunately we do not see the complex
structure.

\subsubsection{The projective model}

We consider the complexification $V(\C)$ of $V$ and the
corresponding projective space
\begin{align}
P(V(\IC)) &= (V(\IC)\bs \{0\} ) /\IC^{\ast}.
\end{align}
The zero quadric
\begin{align}
\label{def:N} \calN &=\{ [Z]\in P(V(\IC)) ;\quad \text{$(Z,Z)=0$}\}
\end{align}
is a closed algebraic subvariety. The subset
\begin{align}
\calK &=\{ [Z]\in P(V(\IC)) ;\quad \text{$(Z,Z)=0$,
$(Z,\bar{Z})>0$}\}
\end{align}
of the zero quadric is a complex manifold of dimension $n$
consisting of two connected components. The orthogonal group
$\Orth_V(\R)$ acts transitively on $\calK$. The subgroup
$\Orth^+_V(\R)$ of elements whose spinor norm equals the determinant
preserves the components of $\calK$, whereas
$\Orth_V(\R)\bs\Orth^+_V(\R)$ interchanges them. We choose one fixed
component of $\calK$ and denote it by $\calK^+$. If $Z\in V(\C)$ we
write $Z=X+iY$ with $X,Y\in V(\R)$ for the decomposition into real
and imaginary part.

\begin{lemma}\label{lem:grk}
The assignment $[Z]\mapsto v(Z):=\R X+\R Y$ defines a real analytic
isomorphism $\calK^+\to \Gr(V)$.
\end{lemma}

\begin{proof}
If $Z \in V(\C)$,  then the condition $[Z]\in \calK$ is equivalent
to
\begin{equation}\label{ob}
X \perp Y,\quad\text{and} \quad (X,X)=(Y,Y)>0.
\end{equation}
But this means that $X$ and $Y$ span a two dimensional positive
definite subspace of $V(\R)$ and thereby define an element of
$\Gr(V)$. Conversely for a given $v\in \Gr(V)$ we may choose a
(suitably oriented) orthogonal basis $X$, $Y$ as in \eqref{ob} and
obtain a unique $[Z]=[X+iY]\in \calK^+$. (Then $[\bar{Z}]\in\calK$
corresponds to the same $v \in \Gr(V)$.) We get a real analytic
isomorphism between  $\calK^+$ and $\Gr(V)$.
\end{proof}

The advantage of the projective model is that it comes with a
natural complex structure.
However, it is not the direct analogue of the upper half plane, the
standard model for the hermitian symmetric space for $\Sl_2(\R)$.

\subsubsection{The tube domain model}

We may realize $\calK^+$ as a tube domain in the following way.
Suppose that $e_1 \in V$ is a non-zero isotropic vector and $e_2\in
V$ with $(e_1,e_2)=1$.  We define a rational quadratic subspace
$W\subset V$ by $W=V\cap e_1^\perp \cap e_2^\perp$.  Then $W$ is
Lorentzian, that is, has type $(1,n-1)$ and
\[
V = W \oplus \Q e_2 \oplus\Q e_1.
\]
If $Z\in V(\C)$ and $Z=z+ae_2+be_1$ with $z\in W(\C)$ and $a,b\in
\C$, we briefly write $Z=(z,a,b)$. We consider the tube domain
\begin{align}
\label{td1} \calH=\{ z\in W(\C);\quad Q(\Im( z))>0\}.
\end{align}

\begin{lemma}
The assignment
\begin{align}\label{td2}
z \mapsto \psi(z):=[(z,1,-Q(z)-Q(e_2))]
\end{align}
defines a biholomorphic map $\psi:\calH\to \calK$.
\end{lemma}

\begin{proof}
One easily checks that if $z\in\calH$ then $\psi(z)\in \calK$.
Conversely assume that $[Z]\in \calK$ with $Z=X+iY$. From the fact
that $X$, $Y$ span a two dimensional positive definite subspace of
$V(\R)$ it follows that $(Z,e_1)\neq 0$. Thus $[Z]$ has a unique
representative of the form $(z,1,b)$. The condition $Q(Z)=0$ implies
that $b=-Q(z)-Q(e_2)$, and thereby $[Z]=[(z,1,-Q(z)-Q(e_2))]$.
Moreover, from $(Z,\bar{Z})>0$ one easily deduces $Q(\Im(z))>0$. We may
infer that the map $\psi$ is biholomorphic.
\end{proof}

The domain $\calH\subset W(\C)\cong \C^n$ has two components
corresponding to the two cones of positive norm vectors in the
Lorentzian space $W(\R)$.  We denote by $\calH^+$ the component
which is mapped to $\calK^+$ under the above isomorphism.  It can be
viewed as a generalized upper half plane.  The group $\Orth_V^+(\R)$
acts transitively on it.  In the $\Orth(2,1)$ case the domain
$\calH^+$ can be identified with the usual upper half plane $\H$.
For $\Orth(2,2)$ it can be identified with $\H^2$ as we shall see
below. However, a disadvantage of the tube domain model is that the
action of $\Orth_V^+(\R)$ is not linear anymore.

\subsubsection{Lattices}

As before, let $(V,Q)$ be a non-degenerate quadratic space over $\Q$
of type $(2,n)$.

\begin{definition}
A {\em lattice} in $V$ is a $\Z$-module $L\subset V$ such that
$V=L\otimes_\Z\Q$.
\end{definition}

A lattice $L\subset V$ is called {\em integral} if the bilinear form
is integral valued on $L$, that is, $(x,y)\in \Z$ for all $x,y\in
L$. A lattice is called {\em even} if the quadratic form is integral
valued on $L$, that is, $Q(x)\in \Z$ for all $x\in L$. So an even
lattice is a free quadratic module over $\Z$ of finite rank. Clearly
every even lattice is integral.


The dual lattice $L^\vee$
is defined by
\[
L^\vee=\{x\in V;\; \text{$(x,y)\in \Z$ for all $y\in L$}\}.
\]
The lattice $L$ is integral if and only if $L\subset L^\vee$. In
this case the quotient $L^\vee/L$ is a finite abelian group, called
the {\em discriminant group}. If $S$ is the Gram matrix
corresponding to a lattice basis of $L$, we have
\[
|L^\vee/L|=|\det(S)|.
\]

For the rest of this section we assume that $L\subset V$ is an even
lattice.  The orthogonal group $\Orth_L$ is a discrete subgroup of
$\Orth_V(\R)\cong \Orth(2,n)$.  Let
\begin{align}\label{gammacond}
\Gamma\subset  \Orth_L\cap \Orth_V^+(\R)
\end{align}
be a subgroup of finite index. Then $\Gamma$ acts properly
discontinuously on $\Gr(V)$, $\calK^+$, and $\calH^+$. We consider
the quotient
\begin{align}
Y(\Gamma)=\Gamma \bs\calH^+
\end{align}
similarly as in the construction of  Hilbert modular surfaces in
Section \ref{sect:1}.  It is a normal complex space, which is
compact if and only if $V$ is anisotropic.

If $Y(\Gamma)$ is non-compact, it can be compactified by adding
rational boundary components (see e.g.~\cite{BrFr}). These boundary
components are most easily described in the projective model
$\calK^+$.  The boundary points of $\calK^+$ in the zero quadric
$\calN$ correspond to non-trivial  isotropic subspaces of $V(\R)$.

Let $F\subset V(\R)$ be an isotropic line. Then $F$ represents a
boundary point of $\calK^+$. A boundary point of this type is called
{\em special}, otherwise {\em generic}. A set consisting of one
special boundary point is called a zero-dimensional boundary
component.

Let $F\subset V(\R)$ be a two-dimensional isotropic subspace.  The
set of all generic boundary points of $\calK^+$ which can be
represented by an element of $F(\C)$ is called the one-dimensional
boundary component attached to $F$.

By a boundary component we understand a one- or two-dimensional
boundary component. One can show (see \cite{BrFr}, Section 2):

\begin{lemma}
There is a bijective  correspondence between boundary components of
$\calK^+$ in the zero quadric $\calN$ and non-zero  iso\-tro\-pic
subspaces $F\subset V(\R)$ of the corresponding dimension. The
boundary of $\calK^+$ is the disjoint union of the boundary
components. \hfill $\square$
\end{lemma}

A boundary component is called rational if the corresponding
isotropic space $F$ is defined over $\Q$. The union of $\calK^+$
with all rational boundary components is denoted by $(\calK^+)^*$.
The rational orthogonal group $\Orth_V(\Q)\cap\Orth_V^+(\R)$ acts on
$ (\calK^+)^*$. By the theory of Baily-Borel, the quotient
\[
X(\Gamma)=(\calK^+)^*/\Gamma
\]
together with the Baily-Borel topology is a compact Hausdorff space.
There is a natural complex structure on $X(\Gamma)$ as a normal
complex space.  Moreover, using modular forms, one can construct an
ample line bundle on $X(\Gamma)$. Therefore, $X(\Gamma)$ is
projective algebraic. It is called the modular variety associated to
$\Gamma$. Using the theory of canonical models, one can show that
$X(\Gamma)$ is actually defined over a number field (see
\cite{Ku:Duke}).

\subsubsection{Heegner divisors}
Let $\Gamma$ be as above, see \eqref{gammacond}. In order to
understand the geometry of $X(\Gamma)$, we study special divisors on
this variety, obtained from embeddings of modular varieties
corresponding to quadratic subspaces of $V$.

Let $\lambda\in L^\vee$ with $Q(\lambda)<0$. Then the orthogonal
complement $V_{\lambda}=\lambda^{\perp}\subset V$ is a rational
quadratic space of type $(2,n-1)$. Moreover, the orthogonal
complement of $\lambda$ in $\calK^+$,
\[
H_\lambda=\{[Z]\in \calK^+;\; (Z,\lambda)=0\},
\]
is an analytic divisor. It is the hermitian symmetric domain
corresponding to $(V_{\lambda},Q|_{V_{\lambda}})$. Let us briefly
look at the description of $H_\lambda$ in the tube domain model
$\calH^+$ using the above notation. We write $\lambda =
\lambda_W+ae_2+be_1$ with $\lambda_W\in W$ and $a,b\in \Q $. Then
\[
H_\lambda\cong\{z\in \calH^+;\; aQ(z)-(z,\lambda_W)-aq(e_2)-b=0\}
\]
is given by a quadratic equation in the coordinates of $\calH^+$.
(Therefore it is sometimes called a rational quadratic divisor.)

If $\beta\in L^\vee/L$ is fixed and $m$ is a fixed negative rational
number, then
\begin{align}
H(\beta,m)=\sum_{\substack{\lambda\in \beta+L\\Q(\lambda)=m}}
H_\lambda
\end{align}
defines an analytic divisor on $\calK^+$ called the {\em Heegner
  divisor} of discriminant $(\beta,m)$.  If $\Gamma$ acts trivially on
$L^\vee/L$, then, by Chow's lemma, this divisor descends to an
algebraic divisor on $Y(\Gamma)$ (denoted in the same way). By
\cite{Ku:Duke}, it is defined over a number field.  Here we mainly
consider the composite Heegner divisor
\begin{align}
H(m)=\frac{1}{2}\sum_{\beta\in L^\vee/L} H(\beta,m) =
\sum_{\substack{\lambda\in L^\vee/\{\pm 1\}\\Q(\lambda)=m}}
H_\lambda .
\end{align}
It is $\Gamma$-invariant and descends to an algebraic divisor on
$Y(\Gamma)$.  Hence, $Y(\Gamma)$ comes with a natural family of
algebraic divisors indexed by negative rational numbers (with
denominators bounded by the level of $L$). The existence of such a
family is special for orthogonal and unitary groups.

\subsection{Modular forms for $\Orth(2,n)$}

Let $V$, $L$, $\Gamma$ be as above. We write
\[
\widetilde{\calK}^+=\{ Z\in V(\IC)\bs \{0\};\quad [Z]\in \calK^+\}
\]
for the cone over $\calK^+$.

\begin{definition}
\label{orthmodular} Let $k\in \Z$, and let $\chi$ be character of
$\Gamma$.  A meromorphic
function $F$ on $\widetilde{\calH}$ is called a meromorphic {\em
modular form} of weight $k$ and character $\chi$ for the group
$\Gamma$, if
\begin{enumerate}
\item[(i)] $F$ is homogeneous of degree $-k$, i.e., $F(c Z)= c^{-k}F(Z)$ for any $c\in \C\bs\{0\}$;
\item[(ii)] $F$ is invariant under $\Gamma$, i.e., $F(g Z)=\chi(g) F(Z)$ for any $g\in \Gamma$;
\item[(iii)] $F$ is meromorphic  at the boundary.
\end{enumerate}
If $f$ is actually holomorphic on $\widetilde{\calK}^+$ and at the
boundary, it is called a holomorphic modular form.
\end{definition}

By the Koecher principle the boundary condition is automatically
fulfilled if the Witt rank of $V$, that is, the dimension of a
maximal isotropic subspace, is smaller than $n$. (Note that because
of the signature the Witt rank of $L$ is always $\leq 2$.)


\subsection{The Siegel theta function}

\label{sect:siegeltheta}

Examples of modular forms on orthogonal groups can be constructed
using Eisenstein series similarly as in Section \ref{sect:eis}.
However, we do not discuss this. Here we consider a rather different
source of modular forms, the so called theta lifting. The groups
$\Sl_2(\R)$ and $\Orth(2,n)$ form a dual reductive pair in the sense
of Howe \cite{Ho}. Hence, Howe duality gives rise to a
correspondence between automorphic representations for the two groups. 
Often one can realize this correspondence as a lifting from automorphic forms 
on one group to the other, by integrating
against certain kernel functions given by theta functions.

Let $V$, $L$, $\Gamma$ be as above and assume that $n$ is even so
that $\dim (V)$ is even. Let
\begin{align*}
N=\min \{ a\in\Z_{>0}; \quad \text{$aQ(\lambda)\in\IZ$ for all
$\lambda\in L^\vee$ }\}
\end{align*}
be the {\em level} of $L$. We modify the discriminant of $L$ by a
sign and consider
\[
\Delta=(-1)^{\frac{n+2}{2}}\det S,
\]
where $S$ is the Gram matrix for a lattice basis of $L$. One can
show that $\Delta\equiv 1,0\pmod{4}$. Therefore
$\chi_{\Delta}=\left( \frac{\Delta}{\cdot}\right)$ is a quadratic
Dirichlet character modulo $N$.

For $\lambda\in V(\R)$ and $v\in \Gr(V)$ we have a unique
decomposition $\lambda=\lambda_v+\lambda_{v^\perp}$, where $
\lambda_v$ and $ \lambda_{v^\perp}$ are the orthogonal projections
of $\lambda$ to $v$ and $v^\perp$, respectively.  The positive
definite quadratic form
\[
Q_v(\lambda)=Q(\lambda_v)-Q(\lambda_{v^\perp})
\]
on V is called the majorant associated to $v$. If $Z\in
\widetilde{\calK}^+$, we briefly write $\lambda_Z$ and $Q_Z$ instead
of $\lambda_{v(Z)}$ and $Q_{v(Z)}$, where $v(Z)$ is the positive
definite plane corresponding to $Z$ via Lemma~\ref{lem:grk}.

\begin{definition}
\label{def:siegeltheta} Let $r\in\Z_{\geq 0}$. The Siegel theta function
of weight $r$ of the lattice $L$ is defined by
\begin{align*}
\Theta_r(\tau,Z)&=v^{n/2}\sum_{\lambda\in
L^\vee}\frac{(\lambda,Z)^r}{(Z,\bar{Z})^r}
e\big(Q(\lambda_Z)N\tau+Q(\lambda_{Z^{\perp}})N\bar{\tau}\big)\\
&=v^{n/2}\sum_{\lambda\in
L^\vee}\frac{(\lambda,Z)^r}{(Z,\bar{Z})^r}e\big(
Q(\lambda)Nu+Q_Z(\lambda)Niv\big),
\end{align*}
for $\tau=u+iv\in\IH$ and $Z\in\widetilde{\calK}^+$. Here
$e(w)=e^{2\pi iw}$ as usual (see e.g.~\cite{Bo2}, \cite{Od},
\cite{RS}).
\end{definition}

Because of the rapid decay of the exponential term
$e(Q_Z(\lambda)Niv)$, the series converges normally on
$\IH\times\widetilde{\calK}^+$. It defines a real analytic function,
which is non-holomorphic in both variables, $\tau$ and $Z$.  Using
the Poisson summation formula, or the theory of the Weil
representation, one can show that as a function in $\tau$, the
Siegel theta function satisfies
\begin{align}
\Theta_r(\gamma\tau,Z)=\chi_{\Delta}(d)(c\tau+d)^{r+\frac{2-n}{2}}\Theta_{r}(\tau,Z)
\end{align}
for all $\gamma=\kabcd\in\Gamma_0(N)$, where
\begin{align}
\Gamma_0(N)=\left\{ \abcd\in\Sl_2 (\IZ);\quad c\equiv 0\pmod{N}
\right\}.
\end{align}
Moreover, in the variable $Z$, the function
$\overline{\Theta_r(\tau,Z)}$ transforms as a modular form of weight
$r$ for $\Gamma$. This follows by direct inspection.

We may use the Siegel theta function as an integral kernel to lift
elliptic modular forms for $\Gamma_0(N)$ to modular forms on the
orthogonal group. More precisely, let $f\in
S_k(\Gamma_0(N),\chi_{\Delta})$ be a cusp form for $\Gamma_0(N)$
with character $\chi_{\Delta}$ of weight $k=r+\frac{2-n}{2}$. We
define the theta lift $\Phi_0(Z,f)$ of $f$ by the integral
\begin{align}
\label{def:thetaint}
 \Phi_0 (Z,f)=\int\limits_{\calF} f(\tau)
\overline{\Theta_r(\tau,Z)} v^k \,\frac{du\,dv}{v^2},
\end{align}
where $\calF$ denotes a fundamental domain for $\Gamma_0(N)$.

\begin{theorem}\label{theo_thetalift}
The theta lift $\Phi_0 (Z,f)$ of $f$ is a holomorphic modular form
of weight $r=k-\frac{2-n}{2}$ for the orthogonal group $\Gamma$.
\end{theorem}

\begin{proof}
The transformation properties of the Siegel theta function
immediately imply that  $\Phi_0(Z,f)$ transforms as a modular form
of weight $r$ for the group $\Gamma$. However, it is not clear at
all, that $\Phi_0(Z,f)$ is holomorphic. This can be proved by
computing the Fourier expansion. For details we refer to
e.g.~\cite{Bo2} Theorem 14.3, \cite{Od} Section 5, Theorem 2, or
\cite{RS}.
\end{proof}

\begin{remark}
  The linear map $f\mapsto \Phi_0(Z,f)$ often has a non-trivial kernel.
  The question when it vanishes is related to the vanishing of a
  special value of the standard $L$-function of $f$ \cite{Ral}.
  Therefore it can
  be rather difficult. However, in many cases it is also possible to
  obtain non-vanishing results by computing the Fourier expansion of
  the lift.
\end{remark}

\subsection{The Hilbert modular group as an orthogonal group}

\label{sect:hilbert-orth}

In this section we discuss the accidental isomorphism relating the
Hilbert modular group to an orthogonal group of type $(2,2)$ in
detail.
The Heegner divisors of the previous section give rise to certain
algebraic divisors on Hilbert modular surfaces, known as
Hirzebruch-Zagier divisors \cite{HZ}.

Let $d\in \Q^*$ be not a square, and put $F=\Q(\sqrt{d})$. We
consider the four dimensional $\Q$-vector space
\[
V=\Q\oplus \Q\oplus F
\]
together with the quadratic form $Q(a,b,\nu)=\nu\nu'-ab$, where
$\nu\mapsto \nu'$ denotes the conjugation in $F$.  So $(V,Q)$ is a
rational quadratic space of type $(2,2)$ if $d>0$ and of type
$(3,1)$ if $d<0$. We consider the orthogonal basis
\begin{align*}
v_1&= (1,1,0),  & v_3&=(0,0,1),\\
v_2&= (1,-1,0), & v_4&=(0,0,\sqrt{d}).
\end{align*}
Then $\delta=v_1v_2v_3v_4$ satisfies $\delta^2=d$. According to
Remark \ref{rem:disc} and Theorem \ref{thm:center}, the center
$Z(C^0_V)$ of the even Clifford algebra of $V$ is given by
$Z(C^0_V)=\Q+\Q\delta\cong F$. Moreover, in view of Example
\ref{ex:n4},
\[
C^0_V=Z+Zv_1 v_2 + Zv_2 v_3 + Z v_1 v_3
\]
is isomorphic to the split quaternion algebra $\Mat_2(F)$ over $F$.
This isomorphism can be realized by the assignment
\begin{align*}
1&\mapsto \zxz{1}{0}{0}{1}, & v_2v_3&\mapsto  \zxz{0}{1}{-1}{0},\\
v_1v_2&\mapsto  \zxz{1}{0}{0}{-1}, & v_1v_3&\mapsto
\zxz{0}{1}{1}{0}.
\end{align*}
The canonical involution on $C^0_V$ corresponds to the conjugation
\[
\abcd^*=\zxz{d}{-b}{-c}{a}
\]
in $\Mat_2(F)$. The Clifford norm corresponds to the determinant.
Hence, by Lemma \ref{lem:nleq4}, $\Spin_V$ can be identified with
$\Sl_2(F)$. As algebraic groups over $\Q$ we have $\Spin_V\cong
\operatorname{Res}_{F/\Q}\Sl_2$.  Consequently, the group
$\Gamma_F=\Sl_2(\calO_F)$ and commensurable groups can be viewed as
arithmetic subgroups of $\Spin_V$. For instance, using \eqref{viso2}
below, it is easily seen that $\Gamma_F=\Spin_L$, where $L$ denotes
the lattice $\Z\oplus \Z \oplus \calO_F\subset V$.

We now describe the vector representation explicitly using
Lemmas~\ref{lem:conj} and \ref{lem:viso}.  Let $\sigma=\Ad(v_1)$ be
the adjoint automorphism of $C_V^0$ associated to the basis vector
$v_1$, i.e., $x^\sigma= v_1 x v_1^{-1}$ for $x\in C_V^0$. Then
$\delta^\sigma=-\delta$, and on the center $F$ of $C^0_V$, the
automorphism $\sigma$ agrees with the conjugation in $F/\Q$. On
$\Mat_2(F)$ the action of $\sigma$ is given by
\[
\abcd\mapsto \abcd^\sigma=\zxz{d'}{-c'}{-b'}{a'}.
\]
As in \eqref{def:tildev} let
\begin{align*}
\tilde V&= \{ X\in \Mat_2(F);\quad X^*=X^\sigma\}\\
&= \left\{ X\in \Mat_2(F);\quad X^t=X'\right\}\\
&= \left\{ \zxz{a}{\nu'}{\nu}{b};\quad
   \text{$a,b\in \Q$ and $\nu\in F$}\right\}.
\end{align*}
This is a rational quadratic space together with the quadratic form
\[
\tilde Q(X)=- X^\sigma\cdot X =  -\det(X).
\]
The corresponding bilinear form is
\[
\tilde B(X_1,X_2)= -\tr (X_1\cdot X_2^*),
\]
for $X_1,X_2\in \tilde V$. The group $\Sl_2(F)\cong\Spin_V$ acts
isometrically on $\tilde V$ by
\begin{align}
\label{twistedspin} x\mapsto g.X:=g X g^{-\sigma}=g X {g'}^t,
\end{align}
for $X\in \tilde V$ and $g\in \Sl_2(F)$. A computation shows that in
the present case the isometry of quadratic spaces $\tilde V\to V$,
$X\mapsto X\cdot v_1$, of Lemma $\ref{lem:viso}$ is given by
\begin{align}
\label{viso2} \zxz{a}{\nu'}{\nu}{b}\mapsto (a,b,\nu).
\end{align}

Throughout the rest of this section we work with $\tilde V$ and the
twisted vector representation \eqref{twistedspin}. We assume that $d$
is positive so that $F$ is real quadratic. We now describe the
hermitian symmetric space corresponding to $\Orth_{\tilde V}$ as in
Section \ref{sect:orth}.

The two real embeddings $x\mapsto (x,x')\in \R^2$ induce an
embedding $\tilde V \to \Mat_2(\R)$. Hence we have $\tilde
V(\C)=\Mat_2(\C)$ and
\begin{align*}
\calK &=\{ [Z]\in P(\Mat_2(\C)) ;\quad \text{$\det (Z)=0$,
$-\tr(Z\bar{Z}^*)>0$}\}.
\end{align*}
We consider the isotropic vectors $e_1=\kzxz{-1}{0}{0}{0}$ and
$e_2=\kzxz{0}{0}{0}{1}$ in $\tilde V$, and the orthogonal complement
$W=\tilde{V}\cap e_1^\perp \cap e_2^\perp$. For $z=(z_1,z_2)\in
\C^2\cong W(\C)$ we put
\begin{align*}
M(z)=\zxz{z_1 z_2}{z_1}{z_2}{1} \in \Mat_2(\IC).
\end{align*}
Then $M(z)$ lies in the zero quadric, and $[M(z)]\in \calK$ if and
only if $\Im(z_1)\Im(z_2)>0$. Consequently, we may identify $\H^2$
with $\calH^+$. If we denote by $\calK^+$ the corresponding
component of $\calK$, we obtain a biholomorphic map
\begin{align}
\label{hilbertorth} \H^2\longrightarrow \calK^+, \quad z\mapsto
[M(z)].
\end{align}
It commutes with the actions of $\Sl_2(F)$, where the action on
$\calK^+$ is given by \eqref{twistedspin}. More precisely, in the
cone $\widetilde{\calK}^+$ we have
\begin{align}
\gamma.M(z)=\norm(cz+d) M(\gamma z)
\end{align}
for $\gamma=\kabcd\in \Sl_2(F)$. This implies that modular forms of
weight $k$ in the sense of Definition \ref{orthmodular} can be
identified with Hilbert modular forms of parallel weight $k$ in the
sense of Definition \ref{hilbertmodular}.

We consider in $V$ the lattice
\begin{align}\label{def:L}
L=\Z\oplus \Z \oplus \calO_F \cong \left\{ \zxz{a}{\nu'}{\nu}{b}\in
\tilde V; \quad\text{$a,b\in\IZ$ and $\nu\in\calO_F$ }\right\}.
\end{align}
The dual lattice of $L$ is
\begin{align}
L^\vee=\Z\oplus \Z \oplus \frakd_F^{-1} \cong \left\{
\zxz{a}{\nu'}{\nu}{b}\in \tilde V; \quad\text{$a,b\in\IZ$ and
$\nu\in\frakd_F^{-1}$ }\right\}.
\end{align}
The discriminant group is given by
$L^\vee/L\cong \calO_F /\frakd_F$.

\begin{proposition}
Under the isomorphism $\Spin_V\cong \Sl_2(F)$, the subgroup
$\Spin_L$ is identified with $\Gamma_F$. \hfill $\square$
\end{proposition}

The map \eqref{hilbertorth} induces an isomorphism of modular
varieties $Y(\Gamma_F)\to Y(\Spin_L)$.

\begin{remark}
More generally, let $\fraka$ be a fractional ideal of $F$ and put
$A=\norm(\fraka)$. We may consider the lattices
\begin{align*}
L(\fraka)&=
\left\{ \zxz{a}{\nu'}{\nu}{Ab}\in
\tilde V; \quad\text{$a,b\in\IZ$ and $\nu\in\fraka$ }\right\},\\
L^\vee(\fraka)&=
\left\{ \zxz{a}{\nu'}{\nu}{Ab}\in \tilde V; \quad\text{$a,b\in\IZ$
and $\nu\in\fraka\frakd_F^{-1}$ }\right\}.
\end{align*}
Observe that $L(\fraka)$ is $A$-integral (that is, the bilinear form
has values in $A\Z$), and $L^\vee(\fraka)$ is the $A\Z$-dual of
$L(\fraka)$. The group $\Gamma(\calO_F\oplus\fraka)\subset \Sl_2(F)$ defined in
\eqref{def:fraka} preserves these lattices.
\end{remark}

\subsubsection{Hirzebruch-Zagier divisors}

In view of the above discussion, the construction of Heegner
divisors provides a natural family of algebraic divisors on a
Hilbert modular surface, in this case known as {\em
Hirzebruch-Zagier divisors} \cite{HZ}.

If $A=\kzxz{a}{\lambda'}{\lambda}{b}\in \tilde V$ and
$z=(z_1,z_2)\in\H^2$, then
\[
(M(z),A)=-\tr(M(z)\cdot A^*) = -bz_1 z_2 +\lambda z_1+\lambda' z_2
-a.
\]
The zero locus of the right hand side defines an analytic divisor
on $\H^2$.

\begin{definition}
Let $m$ be a positive integer. The Hirzebruch-Zagier divisor $T_m$
of discriminant $m$ is defined as the Heegner divisor $H(-m/D)$ for
the lattice $L\subset V$, i.e.,
\begin{align*}
T_m
&=\sum_{\substack{(a,b,\lambda)\in L^\vee/\{\pm 1\} \\
ab-\lambda\lambda'=m/D}} \left\{ (z_1,z_2)\in\H^2;\quad
az_1z_2+\lambda z_1+\lambda' z_2 +b=0 \right\}.
\end{align*}
\end{definition}

It defines an algebraic divisor on the Hilbert modular surface
$Y(\Gamma_F)$. Here the multiplicities of all irreducible components
are $1$. (There is no ramification in codimension $1$.) By taking
the closure, we also obtain a divisor on $X(\Gamma_F)$. We will
denote these divisors by $T_m$ as well, since it will be clear from
the context where they are considered.

\begin{remark}
When $m$ is not a square modulo $D$, then $T_m=\emptyset$.
\end{remark}

\begin{example}
The divisor $T_1$ on $X(\Gamma_F)$ can be identified with the image
of the modular curve $X(1)=\overline{\Sl_2(\Z)\bs \H}$ under the
diagonal embedding considered in Section \ref{sect:restriction}.
\end{example}

\section{Additive and multiplicative liftings}

Let $F\subset \R$ be the real quadratic field of discriminant $D$.
Let $(V,Q)$ be the corresponding rational quadratic space of type
$(2,2)$ as in Section~\ref{sect:hilbert-orth}, and let $L\subset V$
be the even lattice \eqref{def:L}.
The corresponding Siegel theta function $\Theta_k(\tau,z)$ in weight
$k$ is modular in both variables $\tau$ and $z$: As a function of
$\tau$, $\Theta_k(\tau,z)$ is a non-holomorphic modular form of
weight $k$ for the group $\Gamma_0(D)$ with character $\chi_D=\left(
\frac{D}{\cdot}\right)$. As a function in $z$,
$\overline{\Theta_k(\tau,z)}$ is a non-holomorphic modular form of
weight $k$ for the Hilbert modular group $\Gamma_F$. For a cusp form
$f\in S_k(D,\chi_D)$ of weight $k$ for $\Gamma_0(D)$ with character
$\chi_D$, we may consider the theta integral $\Phi_0(z,f)$ as in
\eqref{def:thetaint}. By means of Theorem \ref{theo_thetalift} we
find that $\Phi_0(z,f)$ defines a Hilbert cusp form of weight $k$
for the group $\Gamma_F$ (which may vanish identically). Similar
constructions can be done for the Hilbert modular groups
$\Gamma(\calO_F\oplus\fraka)$ and for their congruence subgroups.

\subsection{The Doi-Naganuma lift}

\label{sect:DN}

In the following we discuss the theta lift in more detail. To keep
the exposition simple, we assume that $D=p$ is a prime and
$F=\Q(\sqrt{p})$. We consider the full Hilbert modular group
$\Gamma_F$.

For explicit computations it is convenient to modify the theta
lifting a bit. Let $M_k(p,\chi_p)$ denote the space of holomorphic
modular forms of weight $k$ for $\Gamma_0(p)$ and $\chi_p$. Since
this space is trivial when $k$ is odd, we assume that $k$ is even.
A function $f\in M_k(p,\chi_p)$ has a Fourier expansion
\begin{align*}
f(\tau)=\sum_{n\geq 0}c(n)q^n,
\end{align*}
where $q=e^{2\pi i\tau}$ as usual. We define the ``plus'' and
``minus'' subspace of $M_k(p,\chi_p)$ by
\begin{align}
\label{pm-space} M_k^{\pm}(p,\chi_p)=\{ f\in M_k(p,\chi_p);´\quad
\chi_p(n)=\mp 1\ \Rightarrow\ c(n)=0\},
\end{align}
and write $S_k^{\pm}(p,\chi_p)$ for the subspace of cusp forms.

Examples of modular forms in $M_k^{\pm}(p,\chi_p)$ can be
constructed by means of Eisenstein series. Recall that there are the
two Eisenstein series
\begin{align}\label{defg}
G_k(\tau) &=1+\frac{2}{ L(1-k,\chi_p)} \sum_{n=1}^\infty \sum_{d\mid n} d^{k-1}\chi_p(d) q^n,\\
\label{defh} H_k(\tau) &=\sum_{n=1}^\infty \sum_{d\mid n}
d^{k-1}\chi_p(n/d) q^n
\end{align}
in $M_{k}(p,\chi_p)$ (cf.~\cite{He} Werke p.~818), the former
corresponding to the cusp $\infty$, the latter corresponding to the
cusp $0$. The linear combination
\begin{equation}\label{eis}
E_k^\pm = 1+\sum_{n\geq 1} B_k^\pm(n)q^n= 1+\frac{2}{L(1-k,\chi_p)}
\sum_{n\geq 1} \sum_{d\mid n} d^{k-1}\left( \chi_p(d) \pm
\chi_p(n/d)\right) q^n
\end{equation}
belongs to $M_{k}^{\pm}(p,\chi_p)$.

\begin{proposition}[Hecke]
The space $M_k(p,\chi_p)$ decomposes into the direct sum
\begin{align*}
M_k(p,\chi_p)=M_k^+(p,\chi_p)\oplus M_k^-(p,\chi_p).
\end{align*}
\item Moreover,
\begin{align*}
M_{k}^{\pm}(p,\chi_p) = \C E_k^\pm \oplus S_{k}^{\pm}(p,\chi_p).
\end{align*}
\hfill $\square$
\end{proposition}

Modular forms in the plus space behave in many ways similarly as
modular forms on the full elliptic modular group $\Sl_2(\Z)$. In
fact, Theorem~5 of \cite{BB} states that $M_{k}^{\pm}(p,\chi_p)$ is
isomorphic to the space of vector-valued modular forms of
weight $k$ for $\Sl_2(\Z)$ transforming with the Weil representation
of $L^\vee/L$.

\begin{notation}
For a formal Laurent series $\sum c(n)q^n\in\IC((q))$  we put
\begin{align}
\label{tildecn}
\tilde{c}(n)=\begin{cases}  c(n), & \text{ if } p\nmid n, \\
2c(n), & \text{ if }p\mid n.
\end{cases}
\end{align}
\end{notation}

\begin{proposition}
\label{pairing} Let $f=\sum c(n)q^n\in M_k^\pm(p,\chi_p) $ and
$g=\sum b(n)q^n\in M_{k'}^\pm(p,\chi_p) $. Then
\begin{align*}
\langle f,g \rangle=\sum_{n\in\IZ}\sum_{m\in\IZ}
\tilde{c}(m)b(pn-m)q^n
\end{align*}
is a modular form of weight $k+k'$ for $\Sl_2(\IZ)$. The assignment
$(f,g)\mapsto \langle f,g \rangle $ defines a bilinear pairing.
\end{proposition}

\begin{proof}
This can be proved by interpreting modular forms in the plus space
as vector valued modular forms for $\Sl_2(\Z)$, see \cite{BB}.
\end{proof}

\begin{remark}
Proposition \ref{pairing} implies some amusing identities of divisor sums 
arising from the equalities $\langle E_k^+,E_k^+\rangle =E_{2k}$ for $k=2,4$. Here $E_{2k}$ denotes the Eisenstein series of weight $2k$ for $\Sl_2(\Z)$ normalized such that the constant term is $1$.
\end{remark}

Note that the statement of Proposition \ref{pairing} does not depend
on the holomorphicity of $f$. An analogous result holds for
non-holomorphic modular forms. For instance, the complex conjugate of
the Siegel theta function $\overline{\Theta_k(\tau,Z)}$ of the lattice
$L$ satisfies the plus space condition. This follows from Definition
\ref{def:siegeltheta}, since for $(a,b,\lambda)\in L^\vee$ we have
\[
-pQ(a,b,\lambda)=p(ab-\lambda\lambda')\equiv \square \pmod{p}.
\]

\
\begin{definition}
For $f\in M_k^+(p,\chi_p)$ we define the (modified) theta lift by
the integral
\begin{align*}
\Phi(z,f)=\int\limits_{\Sl_2(\IZ)\bs\IH}\langle
f(\tau),\overline{\Theta_k(\tau,Z)}\rangle v^k\, \frac{du\,dv}{v^2}.
\end{align*}
\end{definition}

The integral converges absolutely if $f$ is a cusp form. If $f$ is
not cuspidal, the integral has to be regularized (see \cite{Bo2}).
By computing the Fourier expansion of $\Phi(z,f)$, the following
theorem can be proved (cf.~\cite{Za}, \cite{Bo2} Theorem~14.3).

\begin{theorem}
\label{theo:ZKB} Let $f=\sum_n c(n)q^n\in M_k^+(p,\chi_p)$. The
theta lift $\Phi(z,f)$ has the following properties.
\begin{enumerate}
\item[(i)] $\Phi(z,f)$ is a Hilbert modular form of weight $k$ for
$\Gamma_F$.
\item[(ii)] It has the Fourier expansion
\begin{align*}
\Phi(z,f)=-\frac{B_k}{2k}\tilde{c}(0)+\sum_{\substack{\nu\in\frakd_F^{-1} \\
\nu\gg 0}}\sum_{d|\nu}d^{k-1}\tilde{c}\left(
\frac{p\nu\nu'}{d^2}\right) q_1^{\nu}q_2^{\nu'},
\end{align*}
where $B_k$ denotes the $k$-th Bernoulli number, and $q_j=e^{2\pi i
z_j}$.
\item[(iii)] The lift takes cusp forms to cusp forms.
\end{enumerate}
\hfill $\square$
\end{theorem}

If we define in addition $\Phi(z,f)$ to be identically zero on
$M_k^-(p,\chi_p)$, we obtain the {\em Doi-Naganuma lift} (see
\cite{DN}, \cite{Na}),
\begin{align*}
\DN:M_k(p,\chi_p)\longrightarrow M_k(\Gamma_F).
\end{align*}

It is a fundamental fact that the Doi-Naganuma lift (and theta lifts
in general) behave well with respect to the actions of the Hecke
algebras.

\begin{theorem}
[Doi-Naganuma, Zagier] The Doi-Naganuma lift maps Hecke eigenforms
to Hecke eigenforms. For a normalized Hecke eigenform $f=\sum_n
c(n)q^n\in M_k(p,\chi_p)$ we have
\begin{align*}
L(\DN(f),s)=L(f,s)\cdot L(f^{\rho},s),
\end{align*}
where $L(f,s)$ denotes the Hecke $L$-function of $f$ and
$f^{\rho}=\sum\overline{c(n)}q^n$. \hfill $\square$
\end{theorem}

Let $\Lambda(f,s)=p^{s/2} ( 2\pi)^{-s}\Gamma(s)L(f,s)$ be the
completed Hecke $L$-function of the eigenform $f$. It has the
functional equation
\begin{align*}
\Lambda(f,s)=C \cdot \Lambda(f^\rho,k-s)
\end{align*}
with a non-zero constant $C\in \C$. Therefore
\begin{align*}
R(s)=p^s\left( 2\pi \right)^{-2s}\Gamma(s)^2 L(f,s)L(f^{\rho},s)
\end{align*}
has the functional equation
\begin{align*}
R(s)=R(k-s),
\end{align*}
which looks as the functional equation of the $L$-function of a
Hilbert modular form of weight $k$, see Theorem \ref{thm:L} in
Section \ref{sect:L}. Moreover, all further analytic properties of
$R(s)$ agree with those of $L$-functions of Hilbert modular forms.
Hence, using a converse theorem (similar to Hecke's converse
theorem), one can infer that $R(s)$  really comes from a Hilbert
modular form.

Originally, this argument led Doi and Naganuma to the discovery of
the lifting. Using the converse theorem argument they were able to
prove the existence of the lifting in the few cases where $\calO_F$
is euclidian. Employing a later result of Vaserstein (see \cite{Ge}
Chapter IV.6) on generators of Hilbert modular groups, the proof can
be generalized.

The theta lifting approach came up later, and was suggested by
Eichler and Shintani and worked out by Kudla, Oda, Vigner\'as, and
others.

\subsection{Borcherds products}

Here we consider the Borcherds lift for Hilbert modular surfaces. It
can be viewed as a multiplicative analogue of the Doi-Naganuma lift.
It takes certain weakly holomorphic elliptic modular forms of weight
$0$ to meromorphic Hilbert modular forms which have an infinite
product expansion resembling the Dedekind eta function. The zeros
and poles of such {\em Borcherds products} are supported on
Hirzebruch-Zagier divisors.

\subsubsection{Local Borcherds products}

\label{sect:locbp}

As a warm up,  we study a local analogue of
Borcherds products at the cusps of Hilbert modular surfaces.  This is
a special case of the more general results for $\Orth(2,n)$ of
\cite{BrFr}.

We return to the setup of Section \ref{sect:1}. In particular,
$F\subset \R$ is a real quadratic field of discriminant $D$ and
$\Gamma_F=\Sl_2(\calO_F)$ denotes the Hilbert modular group. We ask
whether the Hirzebruch-Zagier divisors $T_m$ on $X(\Gamma_F)$ are
$\Q$-Cartier. Since the non-compact Hilbert modular surface
$Y(\Gamma_F)$ is non-singular except for the finite quotient
singularities corresponding to the elliptic fixed points, it is
clear that $T_m$ is $\Q$-Cartier on $Y(\Gamma_F)$. We only have to
investigate the behavior at the cusps.

\begin{lemma}
\label{lem:hzlocal} Let $A=(a,b,\lambda)\in L^\vee$ with
$ab-\lambda\lambda'>0$. The closure of the image in $Y_F$ of
\[
\left\{ (z_1,z_2)\in\H^2;\quad az_1z_2+\lambda z_1+\lambda' z_2 +b=0
\right\}
\]
goes through the cusp $\infty$ if and only if $a= 0$.
\end{lemma}

\begin{proof}
This is an immediate consequence of Proposition~\ref{prop:bbtop}
(3).
\end{proof}

Let $m$ be a positive integer. We define the local Hirzebruch-Zagier
divisor at $\infty$ of discriminant $m$
by
\begin{align*}
T_m^\infty
&=\sum_{\substack{\lambda\in \frakd_F^{-1}/\{\pm 1\} \\
-\lambda\lambda'=m/D\\b\in \Z}} \left\{ (z_1,z_2)\in\H^2;\quad
\lambda z_1+\lambda' z_2 +b=0 \right\}\subset \H^2.
\end{align*}
This divisor is invariant under the stabilizer $\Gamma_{F,\infty}$
of $\infty$.

\begin{theorem}
\label{thm:qcart} The Hirzebruch-Zagier divisor $T_m$ on
$X(\Gamma_F)$ is $\Q$-Cartier.
\end{theorem}

\begin{proof}
We have to investigate the behavior at the cusps. Here we only
consider the cusp $\infty$, the other cusps can be treated in the
same way.

We have to show that there is a small open neighborhood $U\subset
X(\Gamma_F)$ of $\infty$ and a holomorphic function $f$ on $U$ such
that
\[
\dv(f)=r\cdot T_m|_U\in \Div(U)
\]
for some positive integer $r$. Here $T_m|_U$ denotes the restriction
of $T_m$ to $U$. In view of Proposition \ref{prop:cuspinfty} and
Lemma~\ref{lem:hzlocal} it suffices to show that there exists a
$\Gamma_{F,\infty}$-invariant holomorphic function $\tilde f:\H^2\to
\C$ such that $\dv(\tilde f)=r \cdot T_m^\infty$. This follows from
Proposition~\ref{localbp} below.
\end{proof}

\begin{remark}
The statement of Theorem \ref{thm:qcart} does {\em not\/} generalize
to Heegner divisors on $\Orth(2,n)$. For instance, for $n>3$ there
are obstructions to the $\Q$-Cartier property at generic boundary
points, which are related to theta series of even definite lattices
of rank $n-2$ with harmonic polynomials of degree 2. (See
\cite{BrFr}, \cite{Lo}.)
\end{remark}

The local Hirzebruch-Zagier divisor $T_m^\infty$ decomposes as a sum
\begin{align}
\label{tmdec}
T_m^\infty = \sum_{\substack{\lambda\in \frakd_F^{-1}/\calO_F^{*,2}\\ -\lambda\lambda'=m/D\\
\lambda>0}} T_\lambda^\infty,
\end{align}
where $\calO_F^{*,2}$ denotes the subgroup of squares in the unit
group $\calO_F^*$, and
\begin{align}
T_\lambda^\infty =
\sum_{\substack{u\in \calO_F^{*,2} \\
b\in \Z}} \left\{ (z_1,z_2)\in\H^2;\quad \lambda u z_1+\lambda' u'
z_2 +b=0 \right\}.
\end{align}
The divisor $T_\lambda$ is invariant under $\Gamma_{F,\infty}$. In
the following, we construct a holomorphic function on
$\H^2/\Gamma_{F,\infty}$ whose divisor is $T_\lambda^\infty$, using
local Borcherds products \cite{BrFr}. We start by introducing some
notation.

The subset
\begin{align}
\label{def:sm}
 S(m)=\bigcup_{\substack{\lambda\in
\frakd_F^{-1}\\-\lambda\lambda'=m/D}} \{ y \in(\R_{>0})^2;\quad
\lambda y_1 +\lambda'y_2=0\}
\end{align}
of $(\R_{>0})^2$ is a union of hyperplanes. It is invariant under
$\Gamma_{F,\infty}$. The complement $(\R_{>0})^2\setminus S(m)$ is
not connected. The connected components are called the {\em Weyl
chambers} (of $\frakd_F^{-1}$) of index $m$.

Let $W$ be a subset of a Weyl chamber of index $m$ and $\lambda\in
\frakd_F^{-1}$ with $-\lambda\lambda'=m/D$. Then $\lambda$ is called
{\em positive} with respect to $W$, if $\tr(\lambda w)>0$ for all
$w\in W$ (which is equivalent to requiring $\tr(\lambda w_0)>0$ for
some $w_0\in W$). In this case we write
\[
(\lambda,W)>0.
\]
Moreover, $\lambda$ is called {\em reduced} with respect to $W$, if
\[
(u\lambda,W)<0  , \quad \text{and}\quad ( \lambda,W )  >0,
\]
for all $u\in \calO_F^{*,2}$ with $u<1$. This condition is
equivalent to
\[
(\eps_0^{-2}\lambda,W)<0  , \quad \text{and}\quad (\lambda,W )  >0.
\]
It implies that $\lambda>0$. We denote by $R(m,W)$ the set of all
$\lambda\in \frakd_F^{-1}$ with $-\lambda\lambda'=m/D$ which are
reduced with respect to $W$.  (Note that this definition slightly
differs from the one in \cite{BB}.) It is a finite set and
\begin{align}
\label{tmdec2} \{\lambda\in \frakd_F^{-1};\quad
-\lambda\lambda'=m/D\}= \{\pm\lambda u; \quad \text{$\lambda\in
R(m,W)$ and $u\in \calO_F^{*,2}$}\}.
\end{align}

Let $W$ be a subset of a Weyl chamber of index $m$ and $\lambda\in
\frakd_F^{-1}$ with $-\lambda\lambda'=m/D$. We define a holomorphic
function $\psi_\lambda^\infty:\H^2\to \C$ by
\[
\psi_\lambda^\infty(z)=\prod_{u\in \calO^{*,2}_F} \big[ 1-e(\sigma_u
\tr(u\lambda z))\big],
\]
where
\[
\sigma_u=\begin{cases}
+1,& \text{if $(u\lambda,W)>0$,}\\
-1,& \text{if $(u\lambda,W)<0$.}
\end{cases}
\]
The sign $\sigma_u$ has to be inserted to obtain a {\em convergent}
infinite product. By construction we have
$\psi_\lambda^\infty=\psi_{-\lambda}^\infty$ and
\[
\dv(\psi_\lambda^\infty) = T_\lambda^\infty.
\]
Moreover,  the product is invariant under translations
$\kzxz{1}{\mu}{0}{1}\in \Gamma_{F,\infty}$. However,
$\psi_\lambda^\infty$ is {\em not\/} invariant under the full
stabilizer of $\infty$. It defines an automorphy factor
\begin{align}
\label{def:j} J(\gamma,z)=\psi_\lambda^\infty(\gamma
z)/\psi_\lambda^\infty( z)
\end{align}
of $\Gamma_{F,\infty}$ acting on $\H^2$, that is, an element of
$H^1(\Gamma_{F,\infty},\calO(\H^2)^*)$. We need to show that this
automorphy factor is trivial up to torsion.
It suffices to consider what happens under the generator
$\kzxz{\eps_0}{0}{0}{\eps_0^{-1}}$ of the subgroup of diagonal
matrices in $\Gamma_{F,\infty}$.  We have
\begin{align*}
\frac{\psi_\lambda^\infty(\eps_0^2 z)}{\psi_\lambda^\infty( z)} &=
\prod_{u\in \calO^{*,2}_F}  \frac{1-e(\sigma_{u/\eps_0^{2}}
\tr(u\lambda z))}{1-e(\sigma_u \tr(u\lambda z))}.
\end{align*}
In this product only one factor is not equal to $1$. If we assume
that $\lambda$ is reduced with respect to $W$, we obtain
\begin{align*}
\frac{\psi_\lambda^\infty(\eps_0^2 z)}{\psi_\lambda^\infty( z)} &=
\frac{1-e(-\tr(\lambda z))}{1-e(\tr(\lambda z))}\\
&=e(1/2-\tr(\lambda z)).
\end{align*}

On the other hand, we consider the invertible holomorphic function
\[
I_\lambda(z)=e\left(\tr\left(\frac{\lambda}{\eps_0^2-1}z\right)\right)
\]
on $\H^2$. It satisfies
\[
\frac{I_\lambda(\eps_0^2 z)}{ I_\lambda(z)}=  e(\tr(\lambda z)).
\]
Moreover, $I_\lambda(z+\mu)=I_\lambda(z)$ for all $\mu\in
(\eps_0^2-1)\calO_F$. Therefore, up to torsion, the automorphy
factor $J(\gamma,z)$ in \eqref{def:j} can be trivialized with
$I_\lambda(z)$. The function
\begin{align}
\Psi_\lambda^\infty(z)=I_\lambda(z) \cdot
\psi_\lambda^\infty(z)=e\left(\tr\left(\frac{\lambda}{\eps_0^2-1}z\right)\right)
\prod_{u\in \calO^{*,2}_F} \big[ 1-e(\sigma_u \tr(u\lambda z))\big]
\end{align}
is holomorphic on $\H^2$, has divisor $T_\lambda^\infty$, and a
power of it is invariant under $\Gamma_{F,\infty}$. Observe that
$\Psi_\lambda^\infty$ does not depend on the choice of the Weyl
chamber $W$, although the factors $I_\lambda$ and
$\psi_\lambda^\infty$ do.

Now it is easy to construct an analogous function for $T_m^\infty$.
We define the Weyl vector of index $m$ for the Weyl chamber $W$ by
\begin{align}
\label{def:weylvector} \rho_{m,W}=\sum_{\lambda\in R(m,W)}
\frac{\lambda}{\eps_0^2-1}.
\end{align}
Moreover, we define the local Borcherds product for $T_m^\infty$ by
\begin{align}
\label{def:localbp} \Psi_m^\infty(z)&=\prod_{
\substack{\lambda\in \frakd_F^{-1}/\calO_F^{*,2}\\ -\lambda\lambda'=m/D\\
\lambda >0}}\Psi_\lambda^\infty(z) = e\big(\tr(\rho_{m,W} z)\big)
\prod_{\substack{\lambda\in \frakd_F^{-1}\\ -\lambda\lambda'=m/D\\
(\lambda,W)>0}} \big[ 1-e(\tr(\lambda z))\big].
\end{align}

\begin{proposition}
\label{localbp} The divisor of $\Psi_m^\infty$ is equal to
$T_m^\infty$. A power of $\Psi_m^\infty$ is invariant under
$\Gamma_{F,\infty}$. \hfill $\square$
\end{proposition}

\begin{example}
\label{ex:wv} We compute $\Psi_1^\infty$ more explicitly. The point
$(1,\eps_0)\in (\R_{>0})^2$ does not belong to $S(1)$. Hence it lies
in a unique Weyl chamber $W$ of index $1$. The set of $\lambda\in
\frakd_F^{-1}$ with $-\lambda\lambda'=1/D$ which are reduced with
respect to $W$ is given by
\[
R(1,W)=\begin{cases} \{\eps_0^2/\sqrt{D}\},& \text{if $\eps_0\eps_0'=-1$,}\\
\{\eps_0/\sqrt{D}, \eps_0^2/\sqrt{D}\},& \text{if
$\eps_0\eps_0'=+1$.}
\end{cases}
\]
The corresponding Weyl vector is equal to
\[
\rho_{1,W}=\begin{cases}
\frac{\eps_0}{\sqrt{D}}\frac{1}{\tr(\eps_0)},& \text{if $\eps\eps_0'=-1$,}\\
\frac{1+\eps_0}{\tr(\sqrt{D}\eps_0)},& \text{if $\eps\eps_0'=+1$.}
\end{cases}
\]
In the case $\eps\eps_0'=-1$, the point $(\eps_0^{-1},\eps_0)$ lies
in the same Weyl chamber $W$. It is often more convenient to work
with this base point. If $\eps\eps_0'=1$, then
$(\eps_0^{-1},\eps_0)\in S(1)$.
\end{example}

\subsubsection{The Borcherds lift}

For the material of the next two sections we also refer to
\cite{Br-survey}. The Doi-Naganuma lift of the Section \ref{sect:DN}
only defines a non-trivial map when $k>0$. (For $k=0$ we have
$M_k(D,\chi_D)=0$.) It is natural to ask if one can also do
something meaningful in the border case $k=0$ where the Siegel theta
function \eqref{def:siegeltheta} reduces to the theta function
$\Theta_0(\tau,Z)$ associated to the standard Gaussian on $V(\R)$.
To get a feeling for this question, one can pretend that there is a
non-trivial element $f=\sum_n c(n)q^n\in M_0^+(p,\chi_p)$ and
formally write down its lifting according to Theorem \ref{theo:ZKB}.
We find that it has the Fourier expansion
\[
\Phi(z,f)=-\frac{B_0}{2k}\tilde{c}(0)+\sum_{\substack{\nu\in\frakd_F^{-1} \\
\nu\gg 0}}\sum_{d|\nu}\frac{1}{d} \tilde{c}\left(
\frac{p\nu\nu'}{d^2}\right) q_1^{\nu}q_2^{\nu'}.
\]
Reordering the summation, this can be written as
\[
\Phi(z,f)=-\frac{B_0}{2k}\tilde{c}(0) -\sum_{\substack{\nu\in\frakd_F^{-1} \\
\nu\gg 0}} \log(1-q_1^{\nu}q_2^{\nu'})^{\tilde{c}(p\nu\nu')} .
\]
Hence, the lifting looks as the logarithm of a ``modular'' infinite
product, resembling the Dedekind eta function.  The idea of
Borcherds, Harvey and Moore was to drop the assumption on $f$ being
holomorphic and to replace it by something weaker \cite{Bo1},
\cite{Bo6}, \cite{Bo2}, \cite{HM}.  They consider a regularized
theta lift for weakly
  holomorphic modular forms. It leads to meromorphic modular forms
with infinite product expansions (roughly of the above type).

This construction works in greater generality for $\Orth(2,n)$. It
yields  a lift from weakly holomorphic modular forms of weight
$1-n/2$ to meromorphic modular forms on $\Orth(2,n)$ with zeros and
poles supported on Heegner divisors. Here we only consider the
$\Orth(2,2)$-case of Hilbert modular surfaces. Moreover, to simplify
the exposition, we assume that the real quadratic field $F$ has
prime discriminant $p$.

Let $\Gamma$ be a subgroup of $\Sl_2(\Q)$ which is commensurable
with $\Sl_2(\Z)$. Recall that a meromorphic modular form of weight
$k$ with respect to $\Gamma$ is called {\em weakly holomorphic} if
it is holomorphic outside the cusps. At the cusp $\infty$ such a
modular form $f$ has a Fourier expansion of the form
\[
f(\tau)=\sum_{\substack{n\in\Z\\n\geq N}} c(n) q^{n/h},
\]
where $N\in \Z$, and $h\in \Z_{>0}$ is the width of the cusp
$\infty$.
By an elementary argument it can be proved that the Fourier
coefficients of $f$ are bounded by
\begin{align}\label{est}
c(n)=O\left(e^{C \sqrt{n}}\right),\qquad n\to \infty,
\end{align}
for some positive constant $C>0$ depending on the order of the poles
at the various cusps of $\Gamma$ (see \cite{BrFu} Section 3).
This estimate is also a consequence of the (much more precise)
Hardy-Rademacher-Ramanujan asymptotic for the coefficients of weakly
holomorphic modular forms.

Let $W_k(p,\chi_p)$ be the space of weakly holomorphic modular forms
of weight $k$ for the group $\Gamma_0(p)$ with character $\chi_p$.
Any modular form $f$ in this space has a Fourier expansion of the
form $f=\sum_{n\gg-\infty}c(n)q^n$. Similarly as in \eqref{pm-space}
we denote by  $W_k^+(p,\chi_p)$ the subspace of those $f\in
W_k(p,\chi_p)$, whose coefficients $c(n)$ satisfy the plus space
condition, that is, $c(n)=0$ whenever $\chi_p(n)=-1$.

\begin{lemma}
\label{lem:pp} A weakly holomorphic modular form $f=\sum_{n} c(n)
 q^n\in W_k^+(p,\chi_p)$ of weight $k\leq 0$ is uniquely determined
by its {\em principal part}
\[
\sum_{n<0} c(n)q^n\in \C[q^{-1}].
\]
\end{lemma}

\begin{proof}
The difference of two elements of $W_k^+(p,\chi_p)$ with the same
principal part is holomorphic at the cups $\infty$. Using the plus
space condition (Lemma 3 of \cite{BB}), one infers that the
difference is also holomorphic at the cusp $0$. Hence, it is a
holomorphic modular form of weight $k\leq 0$ with Nebentypus, and
therefore vanishes identically.
\end{proof}

\begin{corollary}
\label{fcweak} Let $k\leq 0$. Assume that $f\in W_k^+(p,\chi_p)$ has
principal part in $\Q[q^{-1}]$. Then all Fourier coefficients of $f$
are rational with bounded denominators.
\end{corollary}

\begin{proof}
This follows from Lemma \ref{lem:pp} and the properties of the
Galois action on $W_k(p,\chi_p)$.
\end{proof}

Let $f=\sum_{n} c(n) q^n\in W_k^+(p,\chi_p)$. Then
\[
(\R_{>0})^2\setminus \bigcup_{\substack{m>0\\ c(-m)\neq 0}} S(m)
\]
is not connected. The connected components are called the {\em Weyl
chambers} associated to $f$. If $W\subset (\R_{>0})^2$ is such a
Weyl chamber, then the {\em Weyl vector} corresponding to $f$ and
$W$ is defined by
\begin{align}
\rho_{f,W}=\sum_{m>0} \tilde c(-m) \rho_{m,W}\in F.
\end{align}
Here $\rho_{m,W}$ is given by \eqref{def:weylvector} and we have
used the notation \eqref{tildecn}.

We are now ready to state Borcherds' theorem in a formulation that
fits nicely our setting (see \cite{Bo2} Theorem 13.3 and
\cite{BB} Theorem 9).

\begin{theorem}[Borcherds]
\label{thm:borcherds} Let $f=\sum_{n\gg -\infty}c(n)q^n$ be a weakly
holomorphic modular form in $W_0^+(p,\chi_p)$ and assume that
$\tilde c(n)\in \Z$ for all $n<0$. Then there exists a meromorphic
Hilbert modular form $\Psi(z,f)$ for $\Gamma_F$ (with some
multiplier system of finite order) such that:
\begin{enumerate}
\item[(i)]
The weight of $\Psi$ is equal to the constant term $c(0)$ of $f$.
\item[(ii)]
The divisor $Z(f)$ of $\Psi$ is determined by the principal part of
$f$ at the cusp $\infty$. It equals
\[Z(f)=\sum_{n<0} \tilde c(n)T_{-n}.\]
\item[(iii)]
Let $W$ be a Weyl chamber associated to $f$ and put $N=\min \{n;\;
c(n)\neq 0\}$. The function $\Psi$ has the Borcherds product
expansion
\[
\Psi(z,f)=q_1^{\rho} q_2^{\rho'}
\prod_{\substack{\nu\in\frakd_F^{-1}
\\ (\nu,W)>0}} \left(1-q_1^\nu q_2^{\nu'} \right)^{\tilde
c(p\nu\nu')},
\]
which converges normally for all $z$ with $y_1 y_2 > |N|/p$ outside
the set of poles. Here $\rho=\rho_{f,W}$ is the Weyl vector
corresponding to $f$ and $W$, and $q_j^\nu=e^{2\pi i \nu z_j}$ for
$\nu\in F$.
%
%
\end{enumerate}
\end{theorem}

\begin{proof}
We indicate the idea of the proof. We consider the theta lift
(Section \ref{sect:siegeltheta}) for the lattice $L$ in the
quadratic space $V=\Q\oplus \Q\oplus F$ (Section
\ref{sect:hilbert-orth}) and use the accidental isomorphism
$\Gamma_F\cong \Spin_V$. The corresponding Siegel theta function
$\overline{\Theta_0(\tau, z)}$ in weight $0$ transforms as an
element of $M_0^+(p,\chi_p)$ in the variable $\tau$. As a function
of $z$ it is invariant under $\Gamma_F$. The pairing $\langle
f(\tau),\overline{\Theta_0(\tau, z)}\rangle$ (see Proposition
\ref{pairing}) is a $\Sl_2(\Z)$-invariant function in $\tau$.

We consider the theta integral
\begin{align}
\label{theta1}
\int\limits_{\calF} \langle
f(\tau),\overline{\Theta_0(\tau,Z)}\rangle \, \frac{du\,dv}{v^2},
\end{align}
where $\calF=\{\tau\in \H;\; |\tau|\geq 1, \;|u|\leq 1/2\}$ denotes
the standard fundamental domain for $\Sl_2(\Z)$. Formally it defines
a $\Gamma_F$-invariant function on $\H^2$. Unfortunately, because of
the exponential growth of $f$ at the cusps, the integral diverges.
However, Harvey and Moore discovered that it can be regularized as
follows \cite{HM}, \cite{Bo2}, \cite{Kon}: If the constant term
$c(0)$ of $f$ vanishes, one can regularize \eqref{theta1} by taking
\begin{align}
\label{theta2} \lim_{t\to\infty}\int\limits_{\calF_t} \langle
f(\tau),\overline{\Theta_0(\tau,Z)}\rangle \, \frac{du\,dv}{v^2},
\end{align}
where $\calF_t=\{\tau\in \calF; \;|v|\leq t\}$ denotes the truncated
standard fundamental domain.  So the regularization consists in
prescribing the order of integration. We first integrate over $u$
and then over $v$.  If the constant term of $f$ does not vanish, the
limit in \eqref{theta2} still diverges. It can be regularized by
considering
\begin{align}
\label{theta3} \Phi(z,f,s)=\lim_{t\to\infty}\int\limits_{\calF_t}
\langle f(\tau),\overline{\Theta_0(\tau,Z)}\rangle v^{-s}\,
\frac{du\,dv}{v^2}
\end{align}
for $s\in \C$. The limit exists for $\Re(s)$ large enough and has a
meromorphic continuation to the whole complex plane. We define the
regularized theta integral $\Phi(z,f)$ to be the constant term in
the Laurent expansion of $\Phi(z,f,s)$ at $s=0$.

One can show that $\Phi(z,f)$ defines a $\Gamma_F$-invariant real
analytic function on $\H^2\setminus Z(f)$ with a logarithmic
singularity\footnote{If $X$ is a normal complex space, $D\subset X$
a
  Cartier divisor, and $f$ a smooth function on $X\setminus\supp(D)$, then $f$
  has a logarithmic singularity along $D$, if for any local equation
  $g$ for $D$ on an open subset $U\subset X$, the function $f-\log|g|$
  is smooth on $U$.}  along the divisor $-4 Z(f)$ (\cite{Bo2} \S6).
The Fourier expansion of $\Phi(z,f)$ can be computed explicitly by
applying some partial Poisson summation on the theta kernel. It
turns out that
\[
\Phi(z,f)= -4\log\left| \Psi(z,f) (y_1 y_2)^{c(0)/2} \right| -
2c(0)\left(\log(2\pi)+\Gamma'(1)\right),
\]
in the domain of convergence of the infinite product for
$\Psi(z,f)$. Using this identity and the properties of $\Phi(z,f)$,
one can prove that the infinite product has a meromorphic
continuation to $\H^2$ satisfying the hypotheses of the theorem.
\end{proof}

\begin{remark}
The fact that $\Psi(f,z)$ only converges in a sufficiently small
neighborhood of the cusp $\infty$ is due to the rapid growth of the
Fourier coefficients of weakly holomorphic modular forms, see
\eqref{est}.
\end{remark}

Meromorphic Hilbert modular forms that arise  as liftings of weakly
holomorphic modular forms by Theorem \ref{thm:borcherds} are called
{\em Borcherds products}.

The following two propositions highlight the arithmetic nature of
Borcherds products. Via the $q$-expansion principle (see \cite{Ra},
\cite{Ch}) they imply that a suitable power of a Borcherds product
defines a rational section of the line bundle of Hilbert modular
forms over $\Z$.

\begin{proposition}
\label{bpquot} Any meromorphic Borcherds product is the quotient of
two holomorphic Borcherds products.
\end{proposition}

\begin{proof}
See \cite{BBK} Proposition 4.5.
\end{proof}

\begin{proposition}\label{prop:intbp}
For any holomorphic Borcherds product
 $\Psi$ there exists a positive integer $n$ such that:
\begin{enumerate}
\item[(i)]
$\Psi^n$ is a Hilbert modular form for $\Gamma_F$ with trivial
multiplier system.
\item[(ii)]
All Fourier coefficients of $\Psi^n$ are contained in $\Z$.
\item[(iii)]
The greatest common divisor of the Fourier coefficients of $\Psi^n$
is equal to $1$.
\end{enumerate}
\end{proposition}

\begin{proof}
The first assertion is clear. The second and the third follow by
Corollary \ref{fcweak} from the infinite product expansion given in
Theorem \ref{thm:borcherds}(iii).
\end{proof}

\subsubsection{Obstructions}

The Borcherds lift provides an explicit construction of relations
among Hirzebruch-Zagier divisors on a Hilbert modular surface.  It
is natural to seek for a precise description of those linear
combinations of Hirzebruch-Zagier divisors, which are divisors of
Borcherds products. Since the divisor of a Borcherds product
$\Psi(z,f)$ is determined by the principal part of the weakly
holomorphic modular form $f$, it suffices to understand which
Fourier polynomials $\sum_{n<0} c(n)q^{n}\in \C[q^{-1}]$ can occur
as principal parts of elements of $W_0^+(p,\chi_p)$.

A necessary condition is easily obtained. If $f\in W_k^+(p,\chi_p)$
with Fourier coefficients $c(n)$, and $g\in M_{2-k}^+(p,\chi_p)$
with Fourier coefficients $b(n)$, then the pairing $\langle
f,g\rangle$ is a weakly holomorphic modular form of weight $2$ for
$\Sl_2(\Z)$. Thus
\[
\langle f,g\rangle d\tau
\]
is a meromorphic differential on the Riemann sphere whose only pole
is at the cusp $\infty$. By the residue theorem its residue has to
vanish. But the residue is just the constant term in the Fourier
expansion of $\langle f,g\rangle$. We find that
\begin{align}\label{neccond}
\sum_{n\leq 0} \tilde c(n) b(-n)=0.
\end{align}

Applying this condition to the Eisenstein series $E^+_{2-k}(\tau)$,
see \eqref{eis}, one gets a formula for the constant term of $f$.

\begin{proposition}
\label{constterm} Let $k$ be a non-positive integer. Let $f=\sum_n
c(n) q^n\in W_k^+(p,\chi_p)$. Then
\[
c(0)=-\frac{1}{2}\sum_{n<0} \tilde c(n) B^+_{2-k}(-n).
\]
\hfill $\square$
\end{proposition}

Using Serre duality for vector bundles on Riemann surfaces,
Borcherds showed that the necessary  condition is also sufficient
(see \cite{Bo3} and \cite{BB} Theorem 6).

\begin{theorem}\label{serre}
There exists a weakly holomorphic modular form $f\in
W_k^+(p,\chi_p)$ with prescribed principal part $\sum_{n<0}c(n)q^n$
(where $c(n)=0$ if $\chi_p(n)=-1$), if and only if
\[
\sum_{n<0} \tilde c(n)b(-n)=0
\]
for every cusp form $g=\sum_{m>0} b(m)q^m$ in $S_{2-k}^+(p,\chi_p)$.
\hfill $\square$
\end{theorem}

\begin{corollary}\label{serre2}
A formal power series $\sum_{m\geq 0} b(m)q^m\in \C[[q]]^+$ is the
$q$-expansion of a modular form in $M_{2-k}^+(p,\chi_p)$, if and
only if
\[
\sum_{n\leq 0} \tilde c(n)b(-n)=0
\]
for every  $f=\sum_{n} c(n)q^n$ in $W_{k}^+(p,\chi_p)$.
\end{corollary}

\begin{proof}
This follows immediately from Theorem \ref{serre}, see
\cite{Br-survey}, Corollary 4.2.
\end{proof}

If $X$ is a regular projective algebraic variety, we write $\Div(X)$
for the group of divisors of $X$, and $\rat(X)$ for the subgroup
of  divisors of rational functions on $X$. The first Chow group of
$X$ is the quotient
\[
\ch^1(X)=\Div(X)/\rat(X).
\]
Furthermore, we put $\ch^1(X)_\Q=\ch^1(X)\otimes_\Z \Q$.  Recall
that $\ch^1(X)$ is isomorphic to the Picard group of $X$, the group
of isomorphism classes of algebraic line bundles on $X$. The
isomorphism is given by mapping a line bundle $\calL$ to the class
$\cc_1(\calL)$ of the divisor of a rational section of $\calL$.  The
Chow group $\ch^1(X)$ is an important invariant of $X$. It is
finitely generated.

Let $\pi:\widetilde{X}\to X(\Gamma_F)$ be a desingularization. If
$k$ is a positive integer divisible by the order of all elliptic
fixed points of $\Gamma_F$, then $\calM_k:=\pi^*\calM_k(\Gamma_F)$,
the pullback of the line bundle of modular forms of weight $k$,
defines an element of $\Pic(\widetilde{X})$. We consider its class
in $\ch^1(\widetilde{X})$.  More generally, if $k$ is any rational
number, we chose an integer $n$ such that $nk$ is a positive integer
divisible by $n(\Gamma_F)$ and put $\cc_1(\calM_k)=\frac{1}{n}
\cc_1(\calM_{nk})\in \ch^1(\widetilde{X})_\Q$.

The Hirzebruch-Zagier divisors are $\Q$-Cartier on $X(\Gamma_F)$.
Their pullbacks define elements in $\ch^1(\widetilde{X})_\Q$.  We
want to describe their positions in this Chow group.  To this end we
consider the generating series
\begin{align}\label{genseries}
A(\tau) = \cc_1(\calM_{-1/2}) + \sum_{m>0} \pi^*(T_m) q^m \in
\Q[[q]] \otimes_\Q \ch^1(\widetilde{X})_\Q.
\end{align}
%
%
Combining Theorem \ref{thm:borcherds} and Corollary \ref{serre2} one
obtains the following striking application.

\begin{theorem}\label{hirzebruchzagier}
The divisors $\pi^*(T_m)$ generate a subspace of
$\ch^1(\widetilde{X})_\Q$ of dimension $\leq
\dim(M_{2}^+(p,\chi_p))$. The generating series $A(\tau)$ is a
modular form in $M_{2}^+(p,\chi_p)$ with values in
$\ch^1(\widetilde{X})_\Q$, i.e., an element of
$M^+_2(p,\chi_p)\otimes_\Q\ch^1(\widetilde{X})_\Q$.
\end{theorem}

In other words, if $\lambda$ is a linear functional on
$\ch^1(\widetilde{X})_\Q$, then
\begin{align*}
 \lambda\left(\cc_1(\calM_{-1/2})\right) + \sum_{m>0}  \lambda(\pi^*T_m) q^m \in M_{2}^+(p,\chi_p).
\end{align*}
A typical linear functional, one can take for $\lambda$, is given by
the intersection pairing with a fixed divisor on $\widetilde{X}$.
Theorem \ref{hirzebruchzagier} was first proved by Hirzebruch and
Zagier
by computing intersection numbers of Hirzebruch-Zagier divisors with
other such divisors and with the exceptional divisors coming from
Hirzebruch's resolution of the cusp singularities \cite{HZ}. Their
discovery triggered important investigations by several people,
showing that more generally periods of certain special cycles in
arithmetic quotients of orthogonal or unitary type can be viewed as
the coefficients of Siegel modular forms. For instance, Oda
considered cycles on quotients of $\Orth(2,n)$ given by embedded
quotients of $\Orth(1,n)$ \cite{Od}, and Kudla-Millson studied more
general cycles on quotients of $\Orth(p,q)$ and $\Uni(p,q)$ using
the Weil representation and theta functions with values in closed
differential forms \cite{KM1,KM2,KM3}, see also \cite{Fu} for the
case of non-compact quotients. The relationship of the Kudla-Millson
lift and the regularized theta lift is clarified in \cite{BrFu}.

\begin{proof}[Proof of Theorem \ref{hirzebruchzagier}]
Using Borcherds products, Theorem \ref{hirzebruchzagier} can be
proved as follows (see \cite{Bo3}). In view of Corollary
\ref{serre2} it suffices to show that
\[
\tilde c(0)\cc_1(\calM_{-1/2}) + \sum_{n<0} \tilde c(n)
\pi^*(T_{-n})=0 \in \ch^1(\widetilde{X})_\Q
\]
for every $f=\sum_{n} c(n)q^n$ in $W_{0}^+(p,\chi_p)$ with integral
Fourier coefficients. But this is an immediate consequence of Theorem
\ref{thm:borcherds}: Up to torsion, the Borcherds lift of $f$ is a
rational section of $\calM_{c(0)}$ with divisor $\sum_{n<0} \tilde
c(n) \pi^*(T_{-n})$.
\end{proof}

Notice that we have only used (i) and (ii) of Theorem
\ref{thm:borcherds}.  Using the product expansion (iii) in addition,
one can prove an arithmetic version of Theorem
\ref{hirzebruchzagier}, saying that certain arithmetic
Hirzebruch-Zagier divisors are the coefficients of a modular form in
$M_2^+(p,\chi_p)$ with values in an arithmetic Chow group, see
\cite{BBK}, \cite{Br-survey}.  Finally, we mention that this
argument generalizes to Heegner divisors on quotients of
$\Orth(2,n)$.

\begin{remark}
With some further work it can be proved that the dimension of the
subspace of $\ch(\widetilde{X})_\Q$ generated by the
Hirzebruch-Zagier divisors is equal to $\dim M_2^+(p,\chi_p)$, see
Corollary \ref{cor:dim}.
\end{remark}

\subsubsection{Examples}

Recall that $p$ is a prime congruent to $1$ modulo $4$. By a result
due to Hecke \cite{He} the dimension of $S_2^+(p,\chi_p)$ is equal
to $[\frac{p-5}{24}]$. In particular there exist three such primes
for which $S_2^+(p,\chi_p)$ is trivial, namely $p=5,13,17$. In these
cases $W^+_0(p,\chi_p)$ is a free module of rank $\frac{p+1}{2}$
over the ring $\C[j(p\tau)]$. Therefore it is not hard to compute
explicit bases. For any $m\in \Z_{>0}$ with $\chi_p(m)\neq -1$ there is a
unique $f_m=\sum_{n\geq -m} c_m(n)q^n\in W_0^+(p,\chi_p)$ whose
Fourier expansion starts with
\[
f_m=\begin{cases} q^{-m}+c_m(0)+O(q), &\text{if $p\nmid m$,}\\
  \frac{1}{2} q^{-m}+c_m(0)+O(q), &\text{if $p\mid m$.}
    \end{cases}
\]
The $f_m$ ($m\in \Z_{>0}$) form a base of the space $W_0^+(p,\chi_p)$.
The Borcherds lift $\Psi_m$ of $f_m$ is a Hilbert modular form for
$\Gamma_F$ of weight $c_m(0)=-B^+_2(m)/2$ with divisor $T_m$. Here
$B^+_2(m)$ denotes the $m$-th coefficient of the Eisenstein series
$E_2^+(\tau)$ as before.


\label{sect:borcherdsex}

\bigskip

{\bf The case $p=5$. } We consider the real quadratic field
$F=\Q(\sqrt{5})$. The fundamental unit is given by
$\eps_0=\frac{1}{2}(1+\sqrt{5})$. Here the first few $f_m$ were
computed in \cite{BB}. One obtains:
\begin{align*}
f_1&=
q^{-1} + 5 + 11\,q - 54\,q^{4} + 55\,q^{5} + 44\,q^{6} - 395\,q^{9} + 340\,q^{10} + 296\,q^{11} - 1836\,q^{14}+\dots,\\
f_4&= q^{-4} + 15 - 216\,q + 4959\,q^{4} + 22040\,q^{5}
 - 90984\,q^{6} + 409944\,q^{9} + 1388520\,q^{10} +\dots,\\
\displaybreak[2] f_5&=\tfrac{1}{2} \,q^{-5} + 15 + 275\, q +
27550\,q^{4} + 43893\,q^{5} + 255300\,q^{6} + 4173825\,q^{9}
 + \dots,\\
\displaybreak[2] f_6&= q^{-6} + 10 + 264\,q - 136476\,q^{4} +
306360\,q^{
5} + 616220\,q^{6} - 35408776\,q^{9} + \dots,\\
f_9&= q^{-9} + 35 - 3555\,q + 922374\,q^{4} + 7512885\,q
^{5} - 53113164\,q^{6} + 953960075\,q^{9} +  \dots,\\
f_{10}&= \tfrac{1}{2}\,q^{-10} + 10 + 3400\,q + 3471300\,q^{4} +
9614200\,q^{5} + 91620925\,q^{6} + \dots.
\end{align*}
The Eisenstein series $E_2^+(\tau)\in M_2^+(5,\chi_5)$ has the
Fourier expansion
\[
E_2^+(\tau,0)=1-10 q -30 q^4 -30 q^5-20q^6-70 q^9 -20q^{10}-120
q^{11} - 60 q^{14} - 40 q^{15}-\dots.
\]
One easily shows that the weight of any Borcherds product
is divisible by $5$. By a little estimate one concludes that there
is just one holomorphic Borcherds product of weight $5$, namely
$\Psi_1$. There exist precisely $3$ holomorphic Borcherds products
in weight $10$, namely $\Psi_1^2$, $\Psi_6$, and $\Psi_{10}$. In
weight $15$ there are the holomorphic Borcherds products $\Psi_4$,
$\Psi_5$, $\Psi_1^3$, $\Psi_1 \Psi_6$, and $\Psi_1\Psi_{10}$.

It follows from Lemma \ref{lem:hzlocal} that $T_m$ does not go
through the cusp $\infty$ when $m$ is not the norm of some
$\lambda\in \calO_F$.  In particular, $T_6$ and $T_{10}$ do not meet
$\infty$. This also implies that $S(6)=S(10)=\emptyset$. There is
just one Weyl chamber of index $6$ and $10$ (namely $(\R_{>0})^2$)
and the corresponding Weyl vector is $0$.  The divisor $T_1$ does
meet $\infty$. As in Example \ref{ex:wv}, let $W$ be the unique Weyl
chamber of index $1$ containing $(\eps_0^{-1},\eps_0)$.  The
corresponding Weyl vector is
$\rho_1=\frac{\eps_0}{\sqrt{D}}\frac{1}{\tr(\eps_0)}$. We obtain the
Borcherds product expansions
\begin{align*}
\Psi_1 &=q_1^{\rho_1} q_2^{\rho_1'} \prod_{\substack{\nu\in\frakd_F^{-1} \\  \eps_0\nu'-\eps_0'\nu >0}} \left(1-q_1^\nu q_2^{\nu'}\right)^{\tilde c_1(5\nu\nu')},\\
\Psi_6 &= \prod_{\substack{\nu\in\frakd_F^{-1} \\ \nu\gg0}} \left(1-q_1^\nu q_2^{\nu'}\right)^{\tilde c_6(5\nu\nu')},\\
\Psi_{10} &= \prod_{\substack{\nu\in\frakd_F^{-1} \\ \nu\gg0}}
\left(1-q_1^\nu q_2^{\nu'}\right)^{\tilde c_{10}(5\nu\nu')}.
\end{align*}
Gundlach \cite{Gu} constructed a Hilbert modular form $s_5$ for
$\Gamma_F$ with divisor $T_1$ as a product of $10$ theta functions
of weight $1/2$, see Section \ref{sect:gundex}. We have
$s_5=\Psi_1$. Moreover, $s_{15}$, the symmetric cusp form of weight
$15$, is equal to $\Psi_5$.  For further examples we refer to
\cite{Ma}.

\subsection{Automorphic Green functions}

By Theorem \ref{serre} of the previous section we know precisely
which linear combinations of Hirzebruch-Zagier divisors occur as
divisors of Borcherds products on $Y(\Gamma_F)$. It is natural
to ask, whether every Hilbert modular form on $Y(\Gamma_F)$ whose
divisor is a linear combination of Hirzebruch-Zagier divisors is a
Borcherds product, i.e., in the image of the lift of Theorem
\ref{thm:borcherds}. In this section we discuss this question in
some detail. To answer it, we first simplify the problem. We extend
the Borcherds lift to a larger space of ``input modular forms''
and answer the question for this extended lift. In that way we are
led to automorphic Green functions associated with Hirzebruch-Zagier
divisors.

Let $k$ be an integer, let $\Gamma$ be a subgroup of $\Sl_2(\Q)$ which is
commensurable with $\Sl_2(\Z)$, and $\chi$ a character of $\Gamma$. A
twice continuously differentiable function $f:\H\to \C$ is called a
{\em weak Maass form} (of weight $k$ with respect to $\Gamma$ and
$\chi$), if
\begin{enumerate}
\item[(i)]
$f\left(\frac{a\tau+b}{c\tau+d}\right)=\chi(\gamma)(c\tau+d)^{k}f(\tau)$
for all $\kabcd\in \Gamma$;
\item[(ii)]
there is a $C>0$ such that for any cusp $s\in \Q \cup \{\infty\}$ of
$\Gamma$ and $\delta\in \Sl_2(\Z)$ with $\delta\infty=s$ the
function $f_s(\tau) = j(\delta,\tau)^{-k}  f(\delta\tau)$ satisfies
$f_s(\tau)=O(e^{C v})$ as $v\to \infty$;
\item[(iii)]
$\Delta_k f=0$.
\end{enumerate}
Here
\begin{equation}\label{deflap}
\Delta_k = -v^2\left( \frac{\partial^2}{\partial u^2}+
\frac{\partial^2}{\partial v^2}\right) + ikv\left(
\frac{\partial}{\partial u}+i \frac{\partial}{\partial v}\right)
\end{equation}
denotes the usual hyperbolic Laplace operator in weight $k$ and
$\tau=u+iv$.

So if we compare this with the definition of a weakly holomorphic
modular form, we see that we simply replaced the condition that $f$
be holomorphic on $\H$ by the weaker condition that $f$ be
annihilated by $\Delta_k$, and the meromorphicity at the cusps by
the corresponding growth condition. In particular, any weakly
holomorphic modular form is a weak Maass form.  The third condition
implies that $f$ is actually real analytic. Because of the
transformation behavior, it has a Fourier expansion involving
besides the exponential function a second type of Whittaker
function. (See \cite{BrFu} Section 3 for more details.)

There are two fundamental differential operators on modular forms
for $\Gamma$, the Maass raising and lowering operators
\begin{align*} R_k
  =2i\frac{\partial}{\partial\tau} + k v^{-1} \qquad \text{and} \qquad
  L_k = -2i v^2 \frac{\partial}{\partial\bar{\tau}}.
\end{align*}
If $f$ is a differentiable function on $\H$ satisfying the
transformation law (i) in weight $k$, then $L_k f$ transforms in
weight $k-2$, and $R_k f$ in weight $k+2$. It can be shown that the
assignment
\[
f(\tau)\mapsto \xi_k(f)(\tau):=v^{k-2} \overline{L_k f(\tau)} =
R_{-k} v^k\overline{ f(\tau)}
\]
defines an antilinear map $\xi_k$ from weak Maass forms of weight
$k$ to weakly holomorphic modular forms of weight $2-k$. Its kernel
is precisely the space of weakly holomorphic modular forms in weight
$k$.

We write $N_k(p,\chi_p)$ for the space of weak Maass forms of weight
$k$ with respect to $\Gamma_0(p)$ and $\chi_p$. Let us have a closer
look at the map $\xi_k: N_k(p,\chi_p) \to W_{2-k}(p,\chi_p)$. We
denote by $\calN_k(p,\chi_p)$ the inverse image of
$S_{2-k}(p,\chi_p)$ under $\xi_k$, and its plus subspace by
$\calN_k^+(p,\chi_p)$. (Note that our notation is not consistent
with the notation of \cite{BrFu}.)

\begin{theorem}\label{exseq}
We have the following exact sequence:
\[
\xymatrix{ 0\ar[r]& W_k^+(p,\chi_p) \ar[r]& \calN_k^+(p,\chi_p)
\ar[r]^{\xi_k}& S_{2-k}^+(p,\chi_p) \ar[r] & 0. }
\]
\end{theorem}

\begin{proof}
  This can be proved using Serre duality for the Dolbeault resolution
  of the structure sheaf on a modular curve (see \cite{BrFu} Theorem
  3.7) or by means of Hejhal-Poincar\'e series (see \cite{Br2} Chapter
  1).
\end{proof}

Let $k\leq 0$. For every weak Maass form $f\in \calN_k^+(p,\chi_p)$
there is a unique Fourier polynomial $P(f)=\sum_{n<0}c(n)q^n\in
\C[q^{-1}]$ (with $c(n)=0$ if $\chi_p(n)=-1$) such that
$f(\tau)-P(f)(\tau)$ is bounded as $v\to\infty$. It is called the
{\em principal part} of $f$. This generalizes the notion of the
principal part of a weakly holomorphic modular form.

\begin{proposition}
\label{prop:pp} Let $Q=\sum_{n<0}c(n)q^n\in \C[q^{-1}]$ be a Fourier
polynomial satisfying $c(n)=0$ if $\chi_p(n)=-1$. There exists a
unique $f\in\calN_k^+(p,\chi_p)$ whose principal part is equal to
$Q$.
\end{proposition}

\begin{proof}
See \cite{BrFu} Proposition 3.11.
\end{proof}

This Proposition is a key fact, which suggests to study the
regularized theta lift of weak Maass forms. If $f\in
\calN_0^+(p,\chi_p)$, then we define its regularized theta lift
$\Phi(z,f)$ by \eqref{theta3}, in the same way as for weakly
holomorphic modular forms.

\begin{theorem}
\label{green1} Let $f\in \calN_0^+(p,\chi_p)$ be a weak Maass form
with principal part $P(f)=\sum_{n<0}c(n)q^n$ and constant term
$c(0)$.
\begin{enumerate}
\item[(i)]
The regularized theta integral $\Phi(z,f)$ defines a
$\Gamma_F$-invariant function on $\H^2$ with a logarithmic
singularity along $-4Z(f)$, where
\[
Z(f)=\sum_{n<0} \tilde c(n)T_{-n}.
\]
\item[(ii)]
It is a Green function for the divisor $2 Z(f)$ on $Y(\Gamma_F)$ in
the sense of \cite{SABK}, that is, it satisfies the identity of
currents
\begin{align*}
dd^c[ \Phi(z,f) ] +\delta_{2 Z(f)} = [\omega(z,f)]
\end{align*}
on $Y(\Gamma_F)$. Here $\delta_D$ denotes the Dirac current
associated with a divisor $D$ on $Y(\Gamma_F)$ and $\omega(z,f)$ is
a smooth $(1,1)$-form.
\item[(iii)]
If $\Delta^{(j)}=-y_j^2\left( \frac{\partial^2}{\partial x_j^2}+
\frac{\partial^2}{\partial y_j^2}\right)$ denotes the
  $\Sl_2(\R)$-invariant hyperbolic Laplace operator on $\H^2$ in the
  variable $z_j$, then
\[
\Delta^{(j)}\Phi(z,f)=-2c(0).
\]
\end{enumerate}
\end{theorem}

\begin{proof}
See \cite{Br1} and \cite{BBK}.
\end{proof}

In view of Proposition \ref{prop:pp}, for every positive integer $m$
with $\chi_p(m)\neq -1$, there exists a unique weak Maass form
$f_m\in \calN^+_0(p,\chi_p)$, whose principal part is given by
\[
P(f_m)=\begin{cases}
q^{-m}, &\text{if $p\nmid m$,}\\
  \frac{1}{2} q^{-m}, &\text{if $p\mid m$.}
    \end{cases}
\]
Its theta lift
\[
\Phi_m(z)=\frac{1}{2}\Phi(z,f_m)
\]
can be regarded as an {\em automorphic Green function} for $T_m$.

Let $\pi: \widetilde{X}\to X(\Gamma_F)$ be a desingularization. The
Fourier expansion of $\Phi(z,f)$ can be computed explicitly. It can
be used to determine the growth behavior at the boundary of
$Y(\Gamma_F)$ in $\widetilde{X}$. It turns out that the boundary
singularities are of log and log-log type. More precisely, one can
view $\pi^*\Phi(z,f)$ as a pre-log-log Green function for the divisor
$2\pi^*(Z(f))$ on $\widetilde{X}$ in the sense of \cite{BKK} (see
\cite{BBK} Proposition 2.16).  So the current equation in (ii) does
not only hold for test forms with compact support on $Y(\Gamma_F)$,
but also for test forms which are smooth on $\widetilde{X}$.

Moreover, one finds that $\Phi(z,f)$ can be split into a sum
\begin{align}\label{split}
\Phi(z,f)=-2 \log| \Psi(z,f)|^2 + \xi(z,f),
\end{align}
where $\xi(z,f)$ is real analytic on the whole domain $\H^2$ and
$\Psi(z,f)$ is a meromorphic function on $\H^2$ whose divisor equals
$Z(f)$. If $f$ is weakly holomorphic, the function $\xi(z,f)$ is
simply equal to $- 2c(0)\left(\log(
y_1y_2)+\log(2\pi)+\Gamma'(1)\right)$, and we are back in the case
of Borcherds' original lift. However, if $f$ is an honest weak Maass
form, then $\xi$ is a complicated function and $\Psi$ far from being
modular.

The splitting \eqref{split} implies that the smooth form
$\omega(z,f)$ in Theorem \ref{green1} is given by
\[
\omega(z,f)=dd^c \xi(z,f).
\]
By the usual Poincar\'e-Lelong argument, $\frac{1}{2}\omega(z,f)$
represents the Chern class of the divisor $Z(f)$ in the second
cohomology $H^2(Y(\Gamma_F))$.  One can further show that it is a
square integrable harmonic representative. Moreover,
$\frac{1}{2}\pi^*\omega(z,f)$ is a pre-log-log form on
$\widetilde{X}$, representing the class of $\pi^*Z(f)$  in
$H^2(\widetilde{X},\C)$.

We now discuss the relationship between the Borcherds lift (Theorem
\ref{thm:borcherds}) and its generalization in the present section.
For simplicity, we simply write $\calN_k$, $W_k$, $M_k$, $S_k$ for
the spaces $\calN_k^+(p,\chi_p)$, $W_k^+(p,\chi_p)$,
$M_k^+(p,\chi_p)$, $S_k^+(p,\chi_p)$, respectively.  We denote by
$W_{k0}$ the subspace of elements of $W_{k}$ with vanishing constant
term. Moreover, we denote by $M_k^\vee$ the dual of the vector space
$M_k$.

\begin{theorem}
We have the following commutative diagram with exact rows:
\[
\xymatrix{ 0 \ar[r] & W_{00} \ar[r] \ar[d]& \calN_0 \ar[r]\ar[d] &M_2^\vee \ar[r]\ar[d]& 0 \\
0 \ar[r] & \rat(\widetilde{X})_\C \ar[r] & \Div(\widetilde{X})_\C
\ar[r] & \ch^1(\widetilde{X})_\C\ar[r]& 0 }.
\]
Here the map $\calN_0\to M_2^\vee$ is given by $f_m\mapsto a_m$,
where $a_m$ denotes the functional taking a modular form in $M_2$ to
its $m$-th Fourier coefficient.  The map $M_2^\vee \to
\ch^1(\widetilde{X})_\C$ is defined by $a_m\mapsto \pi^*T_m$ for
$m>0$ and $a_0\mapsto \cc_1(\calM_{-1/2})$. The map $\calN_0 \to
\Div(\widetilde{X})_\C$ is defined by $f\mapsto \pi^*Z(f)$.
\end{theorem}

\begin{proof}
The exactness of the first row is an immediate consequence of
Theorem \ref{exseq}. Moreover, by Theorem \ref{thm:borcherds}, if
$f\in W_{00}$, then $Z(f)\in \rat(\widetilde{X})_\C$.
\end{proof}

\begin{remark}
The map $\calN_0\to \Div(\widetilde{X})_\C$ does not really depend
on the analytic properties of the weak Maass forms. In particular
the Green function $\Phi(z,f)$ associated to $f\in \calN_0$ does not
play a role. However, there is an analogue of the above diagram in
Arakelov geometry. If $\widetilde{\calX}$ is a regular model of
$\widetilde{X}$ over an arithmetic ring and $\calT_m$ denotes the
Zariski closure of $\pi^*T_m$ in $\widetilde{\calX}$, then the pair
\[
\widehat{\calT}_m=(\calT_m, \pi^*\Phi_m)
\]
defines an arithmetic divisor in the sense of \cite{BKK}. The map
$\calN_0\to \widehat{\Div}(\widetilde{\calX})$, defined by
$f_m\mapsto\widehat{\calT}_m$, gives rise to a diagram as above for
the first arithmetic Chow group of $\widetilde{\calX}$. So the
generalized Borcherds lift can be viewed as a map to the group of
arithmetic divisors on $\widetilde{\calX}$ (see \cite{BBK},
\cite{Br-survey}).
\end{remark}

\begin{theorem}\label{converse1}
Let $h$ be a meromorphic Hilbert modular form of weight $r$ for
$\Gamma_F$, whose divisor $\dv(h)=\sum_{n<0} \tilde c(n)T_{-n}$ is a
linear combination of Hirzebruch-Zagier divisors. Then
\[
-2 \log | h(z)^2 (y_1y_2)^r| = \Phi(z,f)+\text{constant},
\]
where $f$ is the unique weak Maass form in $\calN_0$ with principal
part $\sum_{n<0}c(n)q^n$.
\end{theorem}

\begin{proof}
(See \cite{Br2} Chapter 5.) Let $f$ be the unique weak Maass form in
$\calN_0$ with principal part $\sum_{n<0}c(n)q^n$. Then $\Phi(z,f)$
is real analytic on $\H^2\setminus Z(f)$ and has a logarithmic
singularity along $-4Z(f)$. Hence
\[
d(z):=\Phi(z,f)+2 \log | h(z)^2 (y_1y_2)^r|
\]
is a smooth function on $Y(\Gamma_F)$. By Theorem \ref{green1}
(iii), it is subharmonic.

One can show that $d(z)$ is in $L^{1+\eps}(Y(\Gamma_F))$ for some
$\eps>0$ (with respect to the invariant measure coming from the Haar
measure). By results of Andreotti-Vesentini and Yau on \mbox{(sub-)}
harmonic functions on complete Riemann manifolds that satisfy such
integrability conditions it follows that $d(z)$ is constant.
\end{proof}

The question regarding the surjectivity of the Borcherds lift raised
at the beginning of this section is therefore reduced to the
question whether the weak Maass form $f$ in the Theorem is actually
weakly holomorphic. It is answered affirmatively in Theorem
\ref{converse2} below.

\begin{corollary}
\label{h11map}
The assignment $\pi^*T_m\mapsto \frac{1}{2}dd^c
\xi(z,f_m)$ defines a linear map
\[
\ch^1_{HZ}(\widetilde{X})_\C\longrightarrow \calH^{1,1}(Y(\Gamma_F))
\]
from the subspace of $\ch^1(\widetilde{X})_\C$ generated by the
Hirzebruch-Zagier divisors to the space of square integrable
harmonic $(1,1)$-forms on $Y(\Gamma_F)$. \hfill $\square$
\end{corollary}

Composing the map $M_2^\vee \to \ch^1(\widetilde{X})_\C$ with the
map $\ch^1_{HZ}(\widetilde{X})_\C\to\calH^{1,1}(Y(\Gamma_F))$ from
Corollary \ref{h11map}, we obtain a linear map
\[
M_2^\vee\longrightarrow \calH^{1,1}(Y(\Gamma_F)).
\]
On the other hand, we have the Doi-Naganuma lift $S_2\to
S_2(\Gamma_F)$, and there is a natural map from Hilbert cusp forms
of weight $2$ to harmonic $(1,1)$-forms on $Y(\Gamma_F)$ (see
e.g.~\cite{Br1} Section 5). Summing up, we get the following
diagram:
\begin{align}
\label{diag} \xymatrix{ M_2^\vee\ar[r] &\ch^1_{HZ}(\widetilde{X})_\C
\ar[r]&
\calH^{1,1}(Y(\Gamma_F))\\
S_2\ar[u]^{f\mapsto (\cdot,f)}\ar[rr]&& S_2(\Gamma_F)\ar[u] }.
\end{align}

\begin{theorem}
The above diagram \eqref{diag} commutes.
\end{theorem}

\begin{proof}
See \cite{Br1} Theorem 8.
\end{proof}

So the above construction can be viewed as a different approach to
the Doi-Naganuma lift, making its geometric properties quite
transparent.

Using, for instance, the description of the Doi-Naganuma lifting in
terms of Fourier expansions, it can be proved that $S_2\to
S_2(\Gamma_F)$ is injective. As a consequence, we obtain the
following {\em converse theorem} for the Borcherds lift (see
\cite{Br1}, \cite{Br2} Chapter 5).

\begin{theorem}\label{converse2}
Let $h$ be a meromorphic Hilbert modular form for $\Gamma_F$, whose
divisor $\dv(F)=\sum_{n<0} \tilde c(n)T_{-n}$ is given by
Hirzebruch-Zagier divisors. Then there is a weakly holomorphic
modular form $f\in W_0$ with principal part $\sum_{n<0}c(n)q^n$,
and, up to a constant multiple, $h$ is equal to the Borcherds lift
of $f$ in the sense of Theorem \ref{thm:borcherds}. \hfill $\square$
\end{theorem}

\begin{corollary}
\label{cor:dim}
The dimension of $\ch^1_{HZ}(\widetilde{X})_\C$ is
equal to $\dim(M_2)$. \hfill $\square$
\end{corollary}

Notice that the analogue of Theorem \ref{converse1} holds for
arbitrary congruence subgroups of $\Gamma_F$ (more generally also
for $\Orth(2,n)$), whereas the analogue of Theorem \ref{converse2}
is related to the injectivity of a theta lift and therefore more
complicated. So far it is only known for particular arithmetic
subgroups of $\Orth(2,n)$, see \cite{Br2}, \cite{Br-survey}. For
example, if we go to congruence subgroups of the Hilbert modular
group $\Gamma_F$, it is not clear whether the analogue of Theorem
\ref{converse2} holds or not. See also \cite{BrFu2} for this question.

\subsubsection{A second approach}
\label{sect:3.3.1}

The regularized theta lift $\Phi_m(z)=\frac{1}{2}\Phi(z,f_m)$ of the
weak Mass form $f_m\in \calN_0$ is real analytic on $\H^2\setminus
T_m$ and has a logarithmic singularity along $-2T_m$.

Here we present a different, more naive, construction of
$\Phi_m(z)$. For details see \cite{Br1}. The idea is to construct
$\Phi_m(z)$ directly as a Poincar\'e series by summing over the
logarithms of the defining equations of $T_m$. We consider the sum
\begin{align}\label{formalsum}
\sum_{\substack{(a,b,\lambda)\in\Z^2\oplus\frakd_F^{-1}\\
ab-\lambda\lambda'=m/D}} \log\left|\frac{az_1 \bar{z}_2 +\lambda
z_1+\lambda'\bar{z}_2+b}{az_1 z_2 +\lambda
z_1+\lambda'z_2+b}\right|.
\end{align}
The denominators of the summands ensure that this function has a
logarithmic singularity along $-2T_m$ in the same way  as
$\Phi_m(z)$. The enumerators are smooth on the whole $\H^2$. They
are included to make the sum formally $\Gamma_F$-invariant.
Unfortunately, the sum diverges. However, it can be regularized in
the following way. If we put
$Q_0(z)=\tfrac{1}{2}\log\left(\tfrac{z+1}{z-1}\right)$, we may
rewrite  the summands as
\[
\log\left|\frac{az_1 \bar{z}_2 +\lambda
z_1+\lambda'\bar{z}_2+b}{az_1 z_2 +\lambda z_1+\lambda'z_2+b}\right|
= Q_{0}\left(1+\frac{ |az_1 z_2 +\lambda z_1+\lambda' z_2+b|^2}{
2y_1 y_2 m/D }\right).
\]
Now we replace $Q_0$ by the $1$-parameter family $Q_{s-1}$ of
Legendre functions of the second kind (cf.~\cite{AS} \S8), defined
by
\begin{equation}\label{legendre}
Q_{s-1}(z)=\int\limits_0^\infty (z+\sqrt{z^2-1} \cosh u)^{-s}du.
\end{equation}
Here $z>1$ and $s\in \C$ with $\Re(s)>0$. If we insert $s=1$, we get
back the above $Q_0$. Hence we consider
\begin{equation}\label{phis}
\phi_m(z,s)  = \sum_{\substack{ a,b\in\Z \\
\lambda\in\frakd_F^{-1}\\ab- \lambda\lambda'=m/D }}
Q_{s-1}\left(1+\frac{ |az_1 z_2 +\lambda z_1+\lambda' z_2+b|^2}{
2y_1 y_2 m/D }\right).
\end{equation}
It is easily seen that this series converges normally for
$z\in\H^2\setminus T_m$ and $\Re(s)>1$ and therefore defines a
$\Gamma_F$-invariant function, which has logarithmic growth along
$-2T_m$. It is an eigenfunction of the hyperbolic Laplacians
$\Delta^{(j)}$ with eigenvalue $s(s-1)$, because of the differential
equation satisfied by $Q_{s-1}$. Notice that for $D=m=1$ the
function $\Phi_m(z,s)$ is simply the classical resolvent kernel for
$\Sl_2(\Z)$ (cf.~\cite{Hej}, \cite{Ni}). One can compute the Fourier
expansion of $\phi_m(z,s)$ explicitly and use it to obtain a
meromorphic continuation  to $s\in \C$. At $s=1$ there is a simple
pole, reflecting the divergence of the formal sum \eqref{formalsum}.
We define the regularization $\phi_m(z)$ of \eqref{formalsum} to be
the constant term in the Laurent expansion of $\phi_m(z,s)$ at
$s=1$.

It turns out that $\phi_m$ is up to an additive constant equal to
$\Phi_m$. The Green functions $\phi_m$ can be used to give different
proofs of the results of the previous section and of
Theorem~\ref{thm:borcherds}. Similar Green functions on $\Orth(2,n)$
are investigated in the context of the theory of spherical functions
on real Lie groups in \cite{OT}.

\subsection{CM values of Hilbert modular functions}

In this section we consider the values of Borcherds products on
Hilbert modular surfaces at certain CM cycles. We report on some
results obtained in joint work with T. Yang, see \cite{BY}. This
generalizes work of Gross and Zagier on CM values of the
$j$-function \cite{GZ}.

\subsubsection{Singular moduli}

We review  some of the results of Gross and Zagier on the
$j$-function. We begin by recalling some background material.

Let $k$ be a field and $E/k$ an elliptic curve, that is, a
non-singular projective curve over $k$ of genus $1$ together with a
$k$-rational point. If $\operatorname{char} k\neq 2,3$, then by the
Riemann-Roch theorem one finds that $E$ has a Weierstrass equation
of the form
\[
y^2=4x^3-g_2 x- g_3,
\]
with $g_2,g_3\in k$ and $g_2^3-27 g_3^2\neq 0$. The $j$-invariant of
$E$ is defined by
\[
j(E)=1728\frac{g_2^3}{g_2^3-27 g_3^2}.
\]

A basic result of the theory of elliptic curves says that if $k$ is
algebraically closed then two elliptic curves over $k$ are
isomorphic if and only if they have the same $j$-invariant.
Moreover, for every given $a\in k$ there is an elliptic curve with
$j$-invariant $a$. So the assignment $E\mapsto j(E)$ defines a
bijection
\[
\{\text{elliptic curves over $k$}\}/\sim\;\longrightarrow k.
\]

Over $\C$, the theory of the elliptic functions implies that any
elliptic curve is complex analytically isomorphic to a complex torus
$\C/L$, where $L\subset \C$ is a lattice. (Here $g_2=60
G_4(L)$ and $g_3=140 G_6(L)$ where $G_4,G_6$ are the usual
Eisenstein series of weight $4$ and $6$.) Two elliptic curves $E,
E'$ over $\C$ are isomorphic if and only if the corresponding
lattices $L$, $L'$ satisfy
\[
L=a L'
\]
for some $a\in \C^*$. On the other hand it is easily seen that we
have a bijection
\[
\Sl_2(\Z)\bs\H\longrightarrow \{ \text{lattices in $\C$}\}
/\C^*,\qquad [\tau]\mapsto [\Z\tau+\Z].
\]
Summing up, we obtain a bijection
\begin{align}
\label{moduli:C} \Sl_2(\Z)\bs\H\longrightarrow \{\text{elliptic
curves over $\C$}\}/\sim,\qquad [\tau]\mapsto [\C/(\Z\tau+\Z)].
\end{align}
Hence, the $j$-invariant induces a function on $Y(1):=\Sl_2(\Z)\bs\H$.
A more detailed examination of the map in \eqref{moduli:C} shows that
$j$ is a holomorphic function on $Y(1)$ with the Fourier expansion
$j(\tau)=q^{-1}+744+196884 q+\dots$ at the cusp $\infty$.

So we may view the $j$-function as a function on the coarse moduli
space of isomorphism classes of elliptic curves over $\C$. There are
special points on $Y(1)$ which correspond to special elliptic
curves, namely to elliptic curves with complex multiplication.

Let $K/\Q$ be an imaginary quadratic field with ring of integers
$\calO_K$. A point $\tau\in \H$ is called a CM point of type
$\calO_K$ if the corresponding elliptic curve $E_\tau =
\C/(\Z\tau+\Z)$ has complex multiplication $\calO_K\hookrightarrow
\End(E_\tau)$, or equivalently if $\Z\tau+\Z\subset K$ is a
fractional ideal. We may consider the $0$-cycle $\CM(K)\subset Y(1)$
given by the points $\tau$ for which $E_\tau$ has complex
multiplication by $\calO_K$.

The values of the $j$-function at CM points are classically known as
{\em singular moduli}. If $\tau_0$ is a CM point of type $\calO_K$,
then, by the theory of complex multiplication, $j(\tau_0)$ is an
algebraic integer generating the Hilbert class field of $K$.
Moreover, the Galois group $\Gal(H/K)$ acts transitively on
$\CM(K)\subset Y(1)$. This implies that
\[
j(\CM(K))=\prod_{[\tau]\in \CM(K)} j(\tau)
\]
is an integer. It is a natural question to ask for the shape of this
number. At the beginning of the 20-th century, Berwick made
extensive computations of these numbers and conjectured various
congruences \cite{Be}. We listed some values in Table \ref{j:table}.

\begin{table}[h]
\caption{\label{j:table} Some CM values of the $j$-function}
\begin{center}
\begin{tabular}{|r|r|l| }
\hline \rule[-3mm]{0mm}{8mm}
$|\operatorname{disc}(K)|$ & $h(K)$ &$(j(\CM(K)))^{1/3}$\\
\hline
3 & 1 & 0\\
4 & 1 & $2^2\cdot 3$\\
7 & 1 & $3\cdot 5$\\
8 & 1 & $2^2\cdot 5$\\
11 & 1 & $2^5$ \\
19 & 1 & $2^5\cdot 3$ \\
23 & 3 & $5^3\cdot 11 \cdot 17$\\
31 & 3 & $3^3 \cdot 11\cdot 17 \cdot 23$\\
43 & 1 & $2^6 \cdot 3 \cdot 5$\\
47 & 5 & $5^5\cdot 11^2\cdot 23\cdot 29$\\
59 & 3 & $2^{16}\cdot 11$\\
67 & 1 & $2^5 \cdot 3 \cdot 5\cdot 11$\\
71 & 7 & $11^3 \cdot 17^2 \cdot 23\cdot 41\cdot 47\cdot 53$\\
\hline
\end{tabular}
\end{center}
\end{table}

In \cite{GZ}, Gross and Zagier found an explicit formula for the
prime factorization of $j(\CM(K))$ and proved Berwick's conjectures.
More precisely, they considered the function $j(z_1)-j(z_2)$ on
$Y(1)\times Y(1)$.

Let $K_1$ and $K_2$ be two imaginary quadratic fields of
discriminants $d_1$ and $d_2$, respectively. Assume $(d_1,d_2)=1$,
and put $D=d_1 d_2$. We consider the CM cycle $\CM(K_1)\times
\CM(K_2)$ on $Y(1)\times Y(1)$ and put
\begin{align*}
J(d_1,d_2)=\prod_{\substack{[\tau_1]\in \CM(K_1)\\
[\tau_2]\in \CM(K_2)}}
\left(j(\tau_1)-j(\tau_2)\right)^{\frac{4}{w_1
w_2 }},
\end{align*}
where $w_i$ is the number of units in $K_i$.

\begin{theorem}[Gross, Zagier]
\label{GZmain}
We have
\begin{equation}\label{GZformula}
   J(d_1, d_2)^2 =\pm \prod_{\substack{ x, n, n' \in \mathbb
   Z, \\n, n'>0 \\ x^2 +4 n n' =D }} n^{\epsilon(n')}.
\end{equation}
Here  $\epsilon$ is the genus character defined  as follows:
$\epsilon(n)= \prod \epsilon(l_i)^{a_i}$ if $n$ has the prime
factorization $n =\prod l_i^{a_i}$, and
$$
\epsilon(l) =\begin{cases}
       (\frac{d_1}l) &\text{if $l \nmid d_1$,}
       \\
        (\frac{d_2}l) &\text{if $l \nmid d_2$,}
        \end{cases}
$$
for primes $l$ with $(\frac{D}l) \ne -1$.
\end{theorem}

In particular,  this result implies that the prime factors of
$J(d_1,d_2)$ are bounded by $ D/4$. Since $j(\CM(\Q(\sqrt{-3})))=
j(e^{2\pi i/3})=0$,
we obtain an explicit formula for the CM values of $j$ as a special
case. It leads to the values in Table \ref{j:table}.

The surface $Y(1)\times Y(1)$ can be viewed as the Hilbert modular
surface corresponding to the real quadratic ``field'' $\Q\oplus\Q$
of discriminant $1$. Moreover, $j(z_1)-j(z_2)$ is a Borcherds
product on this surface given by
\begin{align}\label{prodjj}
j(z_1)-j(z_2)=q_1^{-1}\prod_{\substack{m>0\\n\in \Z}}(1-q_1^m
q_2^n)^{c(mn)}.
\end{align}
Here $q_j=e^{2\pi i z_j}$, and $c(n)$ is the $n$-th Fourier
coefficient of $j(\tau)-744$. In fact, this is the celebrated
denominator identity of the monster Lie algebra, which is crucial in
Borcherds' proof of the moonshine conjecture. From this viewpoint it
is natural to ask if the formula of Gross and Zagier has a
generalization to Hilbert modular surfaces. In the rest of this
section we report on joint work with T.~Yang on this problem
\cite{BY}. See also \cite{Ya} for further motivation and background
information.

\subsubsection{CM extensions}
\label{sect:3.4.2}

As before, let $F\subset\R$ be a real quadratic field. Let $K$ be a CM
extension of $F$, that is, $K=F(\sqrt{\Delta})$, where $\Delta\in F$
is totally negative. We view both $K$ and $F(\sqrt{\Delta'})$ as
subfields of $\mathbb C$ with $\sqrt{\Delta}, \sqrt{\Delta'} \in
\mathbb \H$.  The field $M =F(\sqrt\Delta, \sqrt{\Delta'})$ is Galois
over $\Q$. There are three possibilities for the Galois group
$\Gal(M/\Q)$ of $M$ over $\Q$:
\begin{align*}
\Gal(M/\Q)=\begin{cases} \Z/2\Z\times\Z/2\Z, & \text{if $K/\Q$ is biquadratic},
\\ \Z/4\Z, & \text{ if $K/\Q$ is cyclic}, \\
D_4, & \text{ if $K/\Q$ is non Galois}.
\end{cases}
\end{align*}

\begin{lemma} Let the notation be as above, and let $\tilde F=\mathbb{Q}(\sqrt{\Delta \Delta'})$.
\begin{enumerate}
\item[(i)]    $K/\Q$ is biquadratic if and only if $\tilde F
   =\mathbb Q$.
\item[(ii)]    $K/\Q$ is cyclic if and only if $\tilde F = F$.
\item[(iii)]    $K/\Q$ is non-Galois if and only if $\tilde F \ne F$
    is a real quadratic field.
    \hfill$\square$
\end{enumerate}
 \end{lemma}

Gross and Zagier considered a biquadratic case. Here we assume that
$K$ is non-biquadratic, i.e., $\tilde F$ is a real quadratic field.
Then $M/\Q$ has an automorphism $\sigma$ of order $4$ such that
 \begin{equation} \label{eq1.2}
  \sigma(\sqrt{\Delta}) =\sqrt{\Delta'}, \quad
  \sigma(\sqrt{\Delta'})=-\sqrt\Delta.
 \end{equation}
 Notice that $K$ has four CM types, i.e., pairs of non complex conjugate complex embeddings: $\Phi=\{1, \sigma \}$, $\sigma\Phi
 =\{\sigma, \sigma^2\}$, $\sigma^2\Phi$, and $\sigma^3 \Phi$. Since
 $K$ is not biquadratic, these CM types are primitive. We write
 $(\tilde K, \tilde \Phi)$ for the reflex of $(K,\Phi)$. Then $\tilde
 K = \Q(\sqrt{\Delta }+\sqrt{\Delta'} )$ and
 $\tilde F$ is the real quadratic subfield of $\tilde K$. We refer to
 \cite{ShimuraCM} for details about CM types and reflex fields.
%

For the rest of this section we assume that the discriminant of $F$ is a
prime $p\equiv 1\pmod{4}$. Moreover, we suppose that the discriminant
$d_K$ of $K$ is given by $d_K=p^2q$ for a prime $q\equiv 1\pmod{4}$.
This assumption guarantees that the class number of $K$ is odd, which
is crucial in the argument of $\cite{BY}$. It implies that $\tilde
F=\Q(\sqrt{q})$ and $d_{\tilde{K}}=q^2p$.  In Table \ref{cmfielddata}
we listed a few CM extensions of $F=\Q(\sqrt{5})$ satisfying the
assumption, including the class number $h_K$, and a system of
representatives for the ideal class group of $K$.

\begin{table}[h]
\caption{\label{cmfielddata} CM extensions of $\Q(\sqrt{5})$}
\begin{center}
\begin{tabular}{|r|l|l|l|}
\hline
\rule[-3mm]{0mm}{8mm}
$q$ & $K=F(\sqrt{\Delta})$ & $h_K$ & $\CL(K)$  \\
\hline
\rule[-3mm]{0mm}{8mm}
$5$ &$\Delta=-\frac{5+\sqrt{5}}{2}$& 1& $\mathcal{O}_K = \darstell{}{\sqrt{\Delta}}$\\
\hline
\rule[-3mm]{0mm}{8mm}
41 & $\Delta=- \frac{13+\sqrt{5}}{2}$ & 1& $\mathcal{O}_K =\darstell{}{\frac{1}{2}(\sqrt{\Delta}+\frac{3+\sqrt{5}}{2})}$\\
\hline
\rule[-3mm]{0mm}{8mm}
61 & $\Delta = - (9+2 \sqrt{5})$ & 1& $\mathcal{O}_K =\darstell{}{\frac{1}{2}(\sqrt{\Delta}+1)}$\\
\hline
\rule[-3mm]{0mm}{8mm}
109 & $\Delta = -\frac{21+\sqrt{5}}{2} $& 1& $\mathcal{O}_K =\darstell{}{\frac{1}{2}(\sqrt{\Delta}+\frac{3+\sqrt{5}}{2})}$\\
\hline
\rule[-3mm]{0mm}{8mm}
241 & $\Delta=- \frac{33+5 \sqrt{5}}{2}$& 3& $\mathcal{O}_K = \darstell{}{\frac{1}{2}(\sqrt{\Delta}+\frac{3+\sqrt{5}}{2})}$, \\
\rule[-3mm]{0mm}{8mm}&&&
$\mathfrak{A} =\darstell{2}{\frac{1}{2}(\sqrt{\Delta}+\frac{9+3 \sqrt{5}}{2})}$,\\
\rule[-3mm]{0mm}{8mm}&&&
$\mathfrak{B} =\darstell{4}{\frac{1}{2}(\sqrt{\Delta}+\frac{9+3 \sqrt{5}}{2})}$\\
\hline
\rule[-3mm]{0mm}{8mm}
281 & $\Delta=- \frac{37+7 \sqrt{5}}{2}$ & 3 & $\mathcal{O}_K = \darstell{}{\frac{1}{2}(\sqrt{\Delta}+\frac{1+\sqrt{5}}{2})}$,\\
\rule[-3mm]{0mm}{8mm}&&&
$\mathfrak{A} = \darstell{2}{\frac{1}{2}(\sqrt{\Delta}+\frac{1+\sqrt{5}}{2})}$,\\
\rule[-3mm]{0mm}{8mm}&&&
$\mathfrak{B} = \darstell{4}{\frac{1}{2}(\sqrt{\Delta}+\frac{9+\sqrt{5}}{2})}$\\
\hline
\rule[-3mm]{0mm}{8mm}
409 & $\Delta=- \frac{41+3 \sqrt{5}}{2}$ & 3&
$\mathcal{O}_K = \darstell{}{\frac{1}{2}(\sqrt{\Delta}+\frac{1+\sqrt{5}}{2})}$,\\
\rule[-3mm]{0mm}{8mm}&&&
$\mathfrak{A} = \darstell{2}{\frac{1}{2}(\sqrt{\Delta}+\frac{7+ 3 \sqrt{5}}{2})}$,\\
\rule[-3mm]{0mm}{8mm}&&&
$\mathfrak{B} = \darstell{4}{\frac{1}{2}(\sqrt{\Delta}+\frac{-1+3 \sqrt{5}}{2})}$\\
\hline
\end{tabular}
\end{center}
\end{table}

\subsubsection{CM cycles}

We now define CM points on Hilbert modular surfaces analogously
to the CM points on the modular curve $Y(1)$ above.
Recall that the Hilbert modular surface $Y(\Gamma_F)$ corresponding to
$\Gamma_F=\Sl_2(\calO_F)$ parameterizes isomorphism classes of triples
$(A, \imath, m)$, where
\begin{enumerate}
\item[(i)]
$A$ is an abelian surface over $\C$,
\item[(ii)]
$ \imath: \calO_F \to \End(A)$ is a real multiplication by $\calO_F$,
\item[(iii)] and
$
m:(P_A, P_A^+) \to\left(
\frakd_F^{-1}, \frakd_F^{-1, +}\right)
$
is an $\calO_F$-linear isomorphism between the polarization module
$P_A$ of $A$ and $\frakd_F ^{-1}$, taking the subset
of polarizations to totally positive elements of
$\frakd_F^{-1}$. 
\end{enumerate}
(See e.g.~\cite{Go}, Theorem 2.17
and \cite{BY} Section 3.)
The moduli interpretation can be used to construct a model of the
Hilbert modular surface $Y(\Gamma_F)$ over $\Q$, see \cite{Ra},
\cite{DePa}, \cite{Ch}.

Let $\Phi=(\sigma_1, \sigma_2)$ be a CM type of $K$. A point $z=(A,
\imath, m) \in Y(\Gamma_F)$ is said to be a CM point of type $(K,
\Phi)$ if one of the following equivalent conditions holds (see
\cite{BY} Section 3 for details):

\begin{enumerate}
\item[(i)]
As a point $z \in \mathbb H^2$, there is $\tau \in K$ such that
$\Phi(\tau) =(\sigma_1(\tau),  \sigma_2(\tau)) =z$ and such that
$\Lambda_\tau=\calO_F \tau + \mathfrak \calO_F$ is a fractional ideal of $K$.

\item[(ii)] $(A, \imath)$ is a CM abelian variety of type $(K, \Phi)$
  with complex multiplication $\imath': \calO_K
  \hookrightarrow \End(A)$ such that $\imath=\imath'|_{\calO_F}$.
\end{enumerate}

Let $\Phi=\{1,\sigma\}$ be the CM type of $K$ defined in Section
\ref{sect:3.4.2}.  Let $\CM(K, \Phi, \calO_F)$ be the CM $0$-cycle in
$Y(\Gamma_F)$ of CM abelian surfaces of type $(K, \Phi)$.  By the
theory of complex multiplication \cite{ShimuraCM}, the field of moduli
for $\CM(K, \Phi,\calO_F)$ is the reflex field $\tilde K$ of $(K,
\Phi)$.  In fact, one can show that the field of moduli for
\[
\CM(K)=\CM(K, \Phi, \calO_F) + \CM(K, \sigma^3\Phi, \calO_F)
\]
is $\mathbb Q$ (see \cite{BY}, Remark 3.5). Therefore, if $\Psi$ is a rational
function on $Y(\Gamma_F)$, i.e., a Hilbert modular function for $\Gamma_F$
over
$\mathbb Q$, then $ \Psi(\CM(K))$ is a rational number.
The purpose of the following section is to find
a formula for this number, when $\Psi$
is given by a Borcherds product.

\subsubsection{CM values of Borcherds products}

We keep the above assumptions on $F$ and $K$. We denote by $W_{K}$
the number of roots of unity in $K$. For an ideal $\fraka$ of
$\tilde F$ we consider the representation number
\[
\rho(\mathfrak a) =\#\{ \mathfrak A \subset \calO_{\tilde K};\;
     N_{\tilde K/\tilde F} \mathfrak A = \mathfrak a\}
\]
of $\fraka$ by integral ideals of $\tilde K$.  We briefly write
$|\fraka|$ for the norm of $\fraka$. For a non-zero element $t\in
d_{\tilde K/\tilde F}^{-1}$ and a prime ideal $\mathfrak l$ of
$\tilde F$, we put
\begin{align*}
B_t(\mathfrak l) &= \begin{cases}
     (\ord_{\mathfrak l} t +1) \rho(t d_{\tilde K/\tilde F}  \mathfrak l^{-1}) \log
     | \mathfrak l| &\text{ if $\mathfrak l$ is non-split in $\tilde
     K$,}\\
     0  &\text{ if $\mathfrak l$ is split in $\tilde K$,}
     \end{cases}
\end{align*}
and \begin{align*}
B_t & =\sum_{\mathfrak l} B_t(\mathfrak l).
\end{align*}
We remark  that $\rho(\mathfrak a) =0$ for a non-integral ideal
$\mathfrak a$, and that for every $t \ne 0$, there are at most
finitely many  prime ideals $\mathfrak l$ such that $B_t(\mathfrak
l)\ne 0$. In fact, when $t >0 > t'$, then $B_t =0$ unless there is
{\it exactly one} prime ideal $\mathfrak l$ such that
$\chi_{\mathfrak l}(t)=-1$, in which case $B_t=B_t(\mathfrak l)$
(see \cite{BY}, Remark 7.3). Here $\chi =\prod_{\mathfrak l}
\chi_{\mathfrak l}$ is the quadratic Hecke character of $\tilde F$
associated to $\tilde K/\tilde F$. The following formula for the CM
values of Borcherds products is proved in \cite{BY}.

\begin{theorem}
\label{BYmain} 
Let $f=\sum_{n\gg -\infty}c(n)q^n\in W_0^+(p,\chi_p)$, and assume that
$\tilde c(n)\in \Z$ for all $n<0$, and $c(0)=0$.  Then the Borcherds
lift $\Psi=\Psi(z,f)$ (see Theorem \ref{thm:borcherds}) is a rational
function on $Y(\Gamma_F)$, whose value at the CM cycle $\CM(K)$
satisfies $$
\log |\Psi(\CM(K))| = \frac{W_{\tilde K}}4 \sum_{m > 0}
\tilde c(-m) b_m, $$
where
\[
b_m =\sum_{\substack{ t =\frac{n+m\sqrt q}{2p} \in d_{\tilde
K/\tilde F}^{-1}\\ |n| < m \sqrt q }} B_t.
\]
\end{theorem}

Observe that the number of roots of unity $W_{\tilde K}$ is equal to
$2$ unless $p=q=5$, in which case $W_{\tilde K}=10$. The theorem shows
that the prime factorization of $\Psi(\CM(K))$ is determined by the
arithmetic of the reflex field $\tilde K$.

\begin{corollary}
\label{cor1.2} Let the notation be as in Theorem
\ref{BYmain}. Then
\begin{equation} \label{eq1.9}
 \Psi(\CM(K))
 =\pm \prod_{\text{$l$  rational prime}}
 l^{ e_l },
\end{equation}
where
$$
e_l = \frac{W_{\tilde K}}4 \sum_{m >0} \tilde c(-m) b_m(l),
$$
and
$$
b_m(l) \log l= \sum_{\mathfrak l |l} \sum_{\substack{ t
=\frac{n+m\sqrt q}{2p} \in d_{\tilde K/\tilde
    F}^{-1}\\ |n| < m \sqrt q }} B_t(\mathfrak l).
$$
Moreover, when $K/\mathbb Q$ is cyclic, the sign in $(\ref{eq1.9})$
is positive.
\end{corollary}

As in the case that Gross and Zagier considered, see Theorem
\ref{GZmain}, we find that the prime factors of the CM value are
small.

\begin{corollary}
\label{cor1.3}
Let the notation and assumption be as in Corollary
\ref{cor1.2}. Then $e_l=0$ unless   $4 p l | m^2 q -n^2$ for some $m
\in M:=\{ m \in \Z_{>0};\; \tilde c(-m) \ne 0\}$ and some integer
$|n| < m \sqrt q$. 
\end{corollary}

\begin{corollary}
\label{cor1.4}
Let the notation and assumption be as in Corollary
\ref{cor1.2}.
Every prime factor   of
$\Psi(\CM(K)) $ is less than or equal to $\frac{N^2 q}{4p}$, where
$N=\max(M)$.
\end{corollary}

We now indicate the idea of the proof of Theorem \ref{BYmain}.  It
roughly follows the analytic proof of Theorem \ref{GZmain} given in
\cite{GZ}, although each step requires some new ideas.   
By the construction of the Borcherds lift and by the results of 
Section \ref{sect:3.3.1}, 
we have
\[
-4\log |\Psi(z,f)| = \Phi(z,f)=\sum_{m > 0}
\tilde c(-m) \phi_m(z),
\]
where $\phi_m(z)$ denotes the automorphic Green function for $T_m$.
Consequently, it suffices to
compute $\phi_m(\CM(K))$.  Using a CM point, the lattice $\Z^2\oplus
\frakd_F^{-1}$ defining the automorphic Green function can be related
to some ideal of the reflex field $\tilde K$ of $(K,\Phi)$.  In that
way, one derives an expression for $\phi_m(\CM(K))$ as an infinite sum
involving arithmetic data of $\tilde K /\tilde F$.

To come up with a finite sum for the CM value $\phi_m(\CM(K))$, we
consider an auxiliary function.  It is constructed using an incoherent
Eisenstein series (see e.g.~\cite{Ku1}) 
of weight $1$ on $\tilde F$ associated to
$\tilde K /\tilde F$.  We consider the central derivative of this
Eisenstein series, take its restriction to $\Q$, and compute its
holomorphic projection.  

In that way we obtain a holomorphic cusp form $h\in S_2^+(p,\chi_p)$
of weight $2$. Its $m$-th Fourier coefficient is the sum of two parts.
One part is the infinite sum for $\phi_m(\CM(K))$, the other part is a
linear combination of the quantity $b_m$ (what we want) and the
logarithmic derivative of the Hecke $L$-series of $\tilde K /\tilde
F$.  Finally, the duality between $W_0^+(p,\chi_p)$ and $S_2^+(p,\chi_p)$ 
of Theorem \ref{serre}, applied to $f$ and $h$, implies a relation for the
Fourier coefficients of $h$, which leads to the claimed formula.

Notice that the assumption in Theorem \ref{BYmain}
that the constant term of $f$ vanishes can
be dropped. Then the Borcherds lift of $f$ is a meromorphic modular
form of non-zero weight, and one can prove a formula for
$\log\|\Psi(\CM(K),f)\|_{\Pet}$, where $\|\cdot \|_{\Pet}$ denotes the
Petersson metric on the line bundle of modular forms (see \cite{BY}
Theorem 1.4).

In a recent preprint \cite{Scho}, Schofer obtained a formula for the
evaluation of Borcherds products on $\Orth(2,n)$ at CM $0$-cycles
associated with biquadratic CM fields by means of a different method.
It would be interesting to use his results to derive explicit formulas
as in Theorem \ref{BYmain} for the values of Hilbert modular functions
at CM cycles associated to  biquadratic CM fields.  
Finally, notice that Goren and Lauter have
recently proved results on the CM values of Igusa genus two invariants
using arithmetic methods \cite{GL}.

\subsubsection{Examples}

We first consider the real quadratic field $F=\Q(\sqrt{5})$ and the
cyclic CM extension $K=\Q(\sqrt{\zeta_5})$, where $\zeta_5=e^{2 \pi i
  /5}$. So $p=q=5$.  If $\sigma$ denotes the complex embedding of $K$
taking $\zeta_5$ to $\zeta_5^2$ then $\Phi=\{1,\sigma\}$ is a CM type
of $K$. We have $\calO_K=\calO_F+\calO_F\zeta_5$, and the corresponding
CM cycle $\CM(K,\Phi)$ is represented by the point
$(\zeta_5,\zeta_5^2)\in \h^2$.

In Section \ref{sect:borcherdsex} we constructed some Borcherds products
for $\Gamma_F$.
Using the basis $(f_m)$ of $W_0^+(p,\chi_p)$ we see that the Borcherds products 
\begin{align*}
R_1(z)&=
\Psi(z,f_6-2f_1)=\frac{\Psi_6}{\Psi_1^2},\\
R_2(z)&=
\Psi(z,f_{10}-2f_1)=\frac{\Psi_{10}}{\Psi_1^2}
\end{align*} 
are rational functions on $Y(\Gamma_F)$ with divisors $T_6-2T_1$ and
$T_{10}-2T_1$, respectively.  Let us see what the above results say
about $R_1(\CM(K))$.  We have $M=\{1,6\}$ and $N=6$.  According to
Corollary \ref{cor1.4}, the prime divisors of $R_1(\CM(K))$ are bounded
by $9$. Consequently, only the primes $2,3,5,7$ can occur in the
factorization. The divisibility criterion given in Corollary
\ref{cor1.3} actually shows that only $2,3,5$ can occur.
The exact value is given by Corollary \ref{cor1.2}. It is equal to 
$R_1(\CM(K)) = 2^{20}\cdot 3^{10}$.

In Table \ref{valuesp=5} we listed some further CM values of $R_1$ and 
$R_2$.

\begin{table}[h]
\caption{\label{valuesp=5} The case $F=\mathbb Q(\sqrt{5})$}
\begin{center}
\begin{tabular}{|c|c|c| }
  \hline \rule[-3mm]{0mm}{8mm}
  $q$ & $ R_1 (\CM(K))$ & $ R_2(\CM(K))$  \\
  \hline $5$ (cyclic) &   $2^{20}\cdot 3^{10}$ &   $2^{20}\cdot 5^{10}$ \\
  \hline  $41$  & $2^{14}\cdot 3^{10}\cdot 61\cdot 73$ &  $2^{14}\cdot 5^9\cdot 37\cdot 41$\\
  \hline $61$ & $2^{20}\cdot 3^6\cdot 13\cdot 97\cdot 109$ &  $2^{20}\cdot 5^9\cdot 61$
  \\
  \hline $109$ &  $2^{20}\cdot 3^8\cdot 61\cdot 157\cdot 193$ &
  $ 2^{20}\cdot 5^{12}\cdot 73$\\
  \hline $149$ &  $2^{20}\cdot 3^{10} \cdot 31^2\cdot 37\cdot 229$ &
  $ 2^{20}\cdot 5^{12} \cdot 17\cdot 113$ \\
  \hline $269$ &  $2^{20}\cdot 3^{10} \cdot 13^{-2}\cdot 37^2\cdot 61\cdot 97 \cdot 349\cdot 433$ &
  $2^{20}\cdot 5^{14}\cdot 13^{-1}\cdot 53\cdot 73\cdot 233$\\
  \hline
\end{tabular}
\end{center}
\end{table}

\end{document}